%% file: 2001-16.tex
\documentclass{gtart}
\input gtoutput

\lognumber{169}
\volumenumber{5}\papernumber{16}\volumeyear{2001}
\pagenumbers{441}{519}
\accepted{5 May 2001}
\received{16 November 2000}
\published{5 May 2001}
\proposed{Robion Kirby}
\seconded{Ronald Fintushel, Tomasz Mrowka}

\usepackage{amsmath,amssymb}
\usepackage[mathscr]{eucal}
\newtheorem{thm}{Theorem}[section]
\newtheorem{prop}[thm]{Proposition}
\newtheorem{lem}[thm]{Lemma}
\DeclareMathOperator{\Spin}{Spin}
\DeclareMathOperator{\sw}{sw}
\DeclareMathOperator{\SW}{SW}
\DeclareMathOperator{\Gr}{Gr}
\DeclareMathOperator{\kernel}{kernel}
\DeclareMathOperator{\cokernel}{cokernel}
\DeclareMathOperator{\End}{End}
\DeclareMathOperator{\clif}{cl}
\DeclareMathOperator{\Hom}{Hom}
\DeclareMathOperator{\dist}{dist}
\DeclareMathOperator{\Ima}{Im}
\DeclareMathOperator{\Rea}{Re}
\newcommand*{\C}{\ensuremath{\mathbb C}}
\newcommand*{\Z}{\ensuremath{\mathbb Z}}
\newcommand*{\R}{\ensuremath{\mathbb R}}
\newcommand*{\mm}{\ensuremath{\mathscr M}}
\newcommand*{\dd}{\ensuremath{\mathscr D}}
\newcommand*{\gaga}{\ensuremath{\mathscr G}}
\newcommand*{\aaa}{\ensuremath{\mathscr A}}
\newcommand*{\sss}{\ensuremath{\mathscr S}}
\newcommand*{\ww}{\ensuremath{\mathscr W}}
\newcommand*{\ccc}{\ensuremath{\mathscr C}}
\newcommand*{\zz}{\ensuremath{\mathscr Z}}
\newcommand*{\ee}{\ensuremath{\mathscr E}}
\newcommand*{\oo}{\ensuremath{\mathscr O}}
\newcommand*{\ooo}{\ensuremath{\mathscr O}}
\newcommand*{\hh}{\ensuremath{\mathbf H}}
\newcommand*{\qq}{\ensuremath{\mathscr Q}}
\newcommand*{\uu}{\ensuremath{\mathscr U}}

\title{The Seiberg--Witten invariants and $4$--manifolds\\with essential tori}
\asciititle{The Seiberg-Witten invariants and 4-manifolds with essential tori}
\author{Clifford Henry Taubes}
\address{Department of Mathematics\\Harvard University\\Cambridge, 
MA 02138, USA}
\email{chtaubes@math.harvard.edu}
\begin{abstract}
A formula is given for the Seiberg--Witten invariants of a $4$--manifold that 
is cut along certain kinds of $3$--dimensional tori. The formula involves a 
Seiberg--Witten invariant for each of the resulting pieces.
\end{abstract}
\asciiabstract{
A formula is given for the Seiberg-Witten invariants of a 4-manifold that 
is cut along certain kinds of 3-dimensional tori. The formula involves a 
Seiberg-Witten invariant for each of the resulting pieces.}

\asciikeywords{Seiberg-Witten invariants, gluing theorems}

\begin{document}
\maketitlepage

\section{Introduction}\label{sec:1}

This article presents a proof of a cited but previously unpublished 
Mayer--Vietoris theorem for the Seiberg--Witten invariants of four
dimensional 
manifolds which contain certain embedded $3$--dimensional tori. This 
Mayer--Vietoris result is stated in a simple, but less than general form in 
Theorem \ref{thm:1.1}, below. The more general statement is given in Theorem
\ref{thm:2.7}. 
Morgan, Mrowka and Szab\'o have a different (as yet unpublished) proof of
this 
theorem. The various versions of the statement of Theorem \ref{thm:2.7} have
been 
invoked by certain authors over the past few years (for example \cite{MT}, 
\cite{FS1,FS2}, \cite{McT}) and so it is past time for the
appearance of its
proof. 

Note that when the $4$--manifold in question is the product of the circle with 
a $3$--manifold, then Theorem \ref{thm:2.7} implies a Mayer--Vietoris theorem 
for
the 
$3$--dimensional Seiberg--Witten invariants. A proof of the latter along 
different lines has been given by Lim \cite{L}. 

The formulation given here of Theorems \ref{thm:1.1} and \ref{thm:2.7} is 
directly
a consequence 
of conversations with Guowu Meng whose conceptual contributions deserve 
this special acknowledgment here at the very outset.

By way of background for the statement of Theorem \ref{thm:1.1}, consider a
compact, 
connected, oriented $4$--manifold with $b^{2+}\ge 1$. Here, $b^{2+}$ denotes
the 
dimension of any maximal subspace of $H^{2}(X;\mathbb{R})$ on which the cup 
product 
pairing is positive definite. (Such a subspace will be denoted here by 
$H^{2+}(X; \mathbb{R})$.)

When $b^{2+} > 1$, then $X$ has an unambiguous Seiberg--Witten
invariant. In its simplest incarnation, the latter is a map from the
set of $\Spin^{\mathbb{C}}$ structures on $X$ to $\mathbb{Z}$ which is
defined up
to $\pm 1$. Moreover, the sign is pinned with the choice of an
orientation for the real line $L_{X}$ which is the product of the top
exterior power of $H^{1}(X; \mathbb{R})$ with that of $H^{2+}(X;
\mathbb{R})$. The Seiberg--Witten invariants in the case where $b^{2+}
= 1$ can also be defined, but with the extra choice of an orientation
for $H^{2+}(X; \mathbb{R})$. In either case, the Seiberg--Witten
invariants are defined via an algebraic count of solutions to a
certain geometrically natural differential equation on $X$. (See
\cite{W1},
\cite{M}, \cite{KKM}.)

Now, imagine that $M  \subset  X$ is a compact, oriented $3$--dimensional 
submanifold. Supposing that $M$ splits $X$ into two manifolds with boundary, 
$X_{+}$ and $X_{-}$, the problem at hand is to compute the Seiberg--Witten 
invariants for $X$ in a Mayer--Vietoris like way in terms of certain
invariants 
for $X_{+}$ , $X_{-}$ and $M$. Such a formula exists in many cases (see,
eg,
\cite{KM}, \cite{MST}, \cite{MMS}, \cite{OS}.) Theorem \ref{thm:1.1} addresses
this problem in the case 
where: 

\begin{itemize}
\item  $M \subset  X$ is a $3$--dimensional torus.

\item  There is a class in $H^{2}(X; \mathbb{Z})$ with 
non-trivial restriction to $H^{2}(M; \mathbb{Z})$.
\end{itemize}

To consider the solution to this problem, invest a moment to discuss the 
structure of the set, 
$\sss(X)$,
of 
$\Spin^{\C}$ 
structures on
$X$.  In particular, 
remark that for any oriented $4$--manifold $X$, this set 
$\sss(X)$
is defined as the  
set of equivalence classes of pairs 
$(Fr, F)$,       
where 
$Fr  \to  X$ 
is a 
principal                      
$SO(4)$ 
reduction of the oriented, general linear frame bundle for 
$TX$, while $F$ is a lift of 
$Fr$ 
to a principal
$\Spin^{\C}(4)$ 
bundle. In this 
regard, remember that 
$SO(4)$ 
can be identified with 
$(SU(2)\times SU(2))/\{\pm 1\}$ 
in which case 
$\Spin^{\C}(4)$ 
appears as
$(SU(2)\times SU(2)\times U(1))/\{\pm 1\}$. 
Here, $SU(2)$ is the
group of 
$2\times 2$, unitary matrices with determinant $1$ and $U(1)$ is the
circle, the group of 
unit length complex numbers. In any event, since 
$\Spin^{\C}(4)$ 
is an 
extension of $SO(4)$ by the circle, any lift, $F$, of 
$Fr$, 
projects back to 
$Fr$ 
as 
a particularly homogeneous principal $U(1)$ bundle over 
$Fr$.

One can deduce from the preceding description of 
$\sss(X)$ 
that
the latter can be  
viewed in a canonical way as a principal
$H^{2}(X; \mathbb{Z})$ 
bundle over a point. 
In particular, 
$\sss(X)$ 
can be put in 1-1 correspondence with
$H^{2}(X; \mathbb{Z})$, 
but  
no such correspondence is natural without choosing first a fiducial element 
in 
$\sss       (X)$. 
However, there is the canonical `first Chern class' map 
\begin{equation}\label{eq:1.2}
c\colon\sss(X)  \to  H^{2}(X; \mathbb{Z}),
\end{equation}
which is induced by the homomorphism from Spin$^{C}(4)$ to $U(1)$ which 
forgets 
the $SU(2)$ factors. With respect to the $H^{2}(X; \mathbb{Z})$ action
on $\sss(X)$, the map  $c$ obeys 
\begin{equation}\label{eq:1.3}
c(es) = e^{2} c(s),
\end{equation}
for any 
$e  \in  H^{2}(X; \mathbb{Z})$ and 
$s  \in  \sss(X)$. Here, and 
below, the cohomology is viewed as a multiplicative group. Note that 
\eqref{eq:1.3}
implies that 
$c$ is never onto,
and not injective when there is 
$2$--torsion in the second cohomology. 
By the way, $c$'s image in the 
$\mod{2}$
cohomology is the second 
Stiefel--Whitney class of $TX$. 

In any event, if $X$ is a compact, oriented $4$--manifold with $b^{2+} > 1$, 
then 
the Seiberg--Witten invariants define, via the map $c$ in \eqref{eq:1.2}, a map
\begin{equation*}
\underline{\sw}\colon H^{2}(X; \mathbb{Z}) \to  \mathbb{Z},
\end{equation*}
which is defined up to $\pm 1$ without any additional choices. That is, 
$\underline{\sw}(z) \equiv \sum_{s:c(s)=z} \sw({s})$, where 
$\sw(s)$ denotes the value of the Seiberg--Witten invariant on the 
class $s  \in  \sss(X)$. Note that $\underline{\sw} = 0$ but for 
finitely many classes in $H^{2}(X; \mathbb{Z})$ if $b^{2+} > 1$.

Now, for a variety of reasons, it proves useful to package the map 
$\underline{\sw}$ in a manner which will now be described. To start, introduce 
$\Z H^{2}(X; \mathbb{Z})$, 
the free $\mathbb{Z}$ module generated by
the elements in the second  
cohomology. The notation in this regard is such that the abelian group 
structure on $H^{2}(X; \mathbb{Z})$ (as a vector space) is represented in a 
multiplicative fashion. For example, the identity element, $1$, corresponds to 
the trivial class, and more generally, the vector space sum of two classes 
is represented as their product. With the preceding notation understood, a 
typical element in 
$\mathbb{Z}H^{2}(X; \mathbb{Z})$ 
consists of a
formal sum  $\sum a(z)z$, where  
the sum is over the classes  $z \in  H^{2}(X; \mathbb{Z})$ with 
$a(z)\in\mathbb{Z}$
being zero but for finitely many classes. Thus, a choice of basis over
$\mathbb{Z}$ for  
$H^{2}(X; \mathbb{Z})$, makes 
elements of $\mathbb{Z}H^{2}(X;
\mathbb{Z})$ into finite Laurent series.         

With $\mathbb{Z}H^{2}(X; \mathbb{Z})$ understood, the invariant
$\underline {\sw}$ in the $b^{2+} > 1$ case can be packaged neatly as an
element in $\mathbb{Z}H^{2}(X; \mathbb{Z})$, namely
\begin{equation}\label{eq:1.5}
\underline{\SW}_{X}  \equiv \sum_{z} \underline{\sw}(z)z.
\end{equation}

In the case where $b^{2+}$ = 1, a choice of orientation for 
$H^{2+}(X; \mathbb{R})$ 
is needed to define $\underline {\sw}$, and in this case the analog of 
\eqref{eq:1.5} is 
a `semi-infinite' power series rather than a finite Laurent series. In this 
regard, a power series such as $\sum_{z} a(z) z$ is termed semi-infinite with 
respect to a given generator of $H^{2+}(X; \mathbb{Z})$ when the following is 
true: 
For any real number $m$, only a finite set of classes 
$z \in H^{2}(X; \mathbb{Z})$ 
have both $a(z) \neq 0$ and cup product pairing less than $m$ with the 
generator. 
In the case of $\underline {\SW}_{X}$, the choice of a Riemannian metric and 
an orientation for $H^{2+}(X; \mathbb{Z})$ determines the generator in 
question. 
However, $\underline {\SW}_{X}$ does not depend on the metric, it depends 
only on the chosen orientation of $H^{2+}(X; \mathbb{R})$. The associated, 
extended 
version of $\mathbb{Z}H^{2}(X; \mathbb{Z})$ which admits such power
series will not be  
notationally distinguished from the original. In any event, when $b^{2+} = 
1$, the extra choice of an orientation for $H^{2}(X; \mathbb{R})$ yields a 
natural 
definition of $\underline {\sw}$ so that \eqref{eq:1.5} makes good sense as 
an element 
in the extended $\mathbb{Z}H^{2}(X; \mathbb{Z})$.

Now, suppose that $X$ is a compact, connected $4$--manifold with boundary, 
$\partial X$, 
with each component of the latter being a $3$--torus. Assume, in addition, 
that 
there is a fiducial class, $\varpi $, in $H^{2}(X; \mathbb{Z})$ whose 
pull-back is 
non-zero in the cohomology of each component of $\partial X$ . Theorem 
\ref{thm:2.5} to come 
implies that such a manifold also has a Seiberg--Witten invariant, $\underline 
{\SW}_{X}$, which lies either in 
$\mathbb{Z}H^{2}(X,\partial X;\Z)$ 
or, in certain cases, a  
particular extension of this group ring which allows semi-infinite series. 
In this case, the extension in question consists of formal power series such 
as $\sum_{z} a(z) z$ where, for 
any given real number 
$m$, only a finite set of $z$'s
have both $a(z) \neq 0$ 
and cup product pairing less than $m$ with $\varpi $.
(This 
extension will not be notationally distinguished from
$\mathbb{Z}H^{2}(X, \partial X ; \Z)$.)  
In any event, $\underline {\SW}$ is defined, as in the no boundary case, via 
an 
algebraic count of the solutions to a version of the Seiberg--Witten 
equations. This invariant is defined up to a sign with the choice of $\varpi$ 
and the sign is fixed with the choice of an orientation for the line 
$L_{X}$ which is the product of the top exterior power or $H^{1}(X,
\partial X ; \mathbb{R})$  
with that of $H^{2+}(X, \partial X ; \mathbb{R})$. Even in the
non-empty boundary case,  $\underline {\SW}$ is a diffeomorphism invariant.

By way of an example, the invariant $\underline {\SW}$ for the product, 
${D}^{2} \times T^{2}$, of the closed, $2$--dimensional disk with the 
torus is $t(1-t^{2})^{-1} = t + t^{3} +\cdots$, 
where $t$ is Poincar\'e dual to the class of the torus. For another example, 
take $n$ to be a positive integer and let $E(n)$ denote the simply connected, 
minimal elliptic surface with no multiple fibers and holomorphic Euler 
characteristic $n$. The invariant $\underline {\SW}$ for the complement in 
$E(n)$ 
of an open, tubular neighborhood of a generic fiber is 
$(t-t^{-1})^{n-1}$, where $t$ is the Poincar\'e dual of a fiber.

With the description of $\underline {\SW}$ in hand, a simple version of
the promised Mayer--Vietoris formula can be stated. In this regard,
mind that certain pairs of elements in $\mathbb{Z}H^{2}(X, \partial
X ; \Z)$ can be multiplied together as formal power series. The
multiplication rule used here is the evident one where $(\sum_{z} a(z)
z)\cdot (\sum_{z} a'(z) z) \equiv \sum_{z} [\sum_{(w,x): wx=z} a(x)
a'(w)] z$.

\begin{thm}\label{thm:1.1} 
Let $X$ be a compact, connected, 
oriented, $4$--manifold with $b^{2+} = 1$ and with boundary 
consisting of a disjoint union of $3$--dimensional tori. Let $M  \subset  
X$ be an embedded, $3$--dimensional torus and suppose that 
there is a fiducial class $\varpi\in H^{2}(X; \mathbb{R})$ whose 
pull-back is non-zero in the cohomology of $M$ and in that of each 
component of the boundary of $X$. 

\begin{itemize}
\item If $M$ splits $X$ as a pair, 
$X_{+} \cup X_{-}$, of $4$--manifolds with boundary, let $j_{\pm}$ 
denote the natural, $\mathbb{Z}$--linear extensions of the canonical 
homomorphisms from  $H^{2}(X_{\pm},\partial X_{\pm};\mathbb{Z})$ to 
$H^{2}(X, \partial X ; \mathbb{Z})$ which arise by coupling the
excision isomorphism  
with those from the
long exact cohomology sequences of the pairs 
$X_{-}, X_{+}\subset  X$. Then $j_{-}(\underline {\SW}_{X_{-}})$ 
and $j_{+}(\underline{\SW}_{X_{+}})$ can be 
multiplied together in $\mathbb{Z}H^{2}(X, \partial X ; \Z)$ and
\begin{equation*}
\underline {\SW}_{X} = j_{-}(\underline {\SW}_{X_{-}})
j_{+}(\underline {\SW}_{X_{+}}).
\end{equation*}
Here, the orientation for the line $L_{X}$ is induced by 
chosen orientations for the analogous lines for $X_{+}$ and 
$X_{-}$. Also, if $X$ is compact and $b^{2+} = 1$, 
then $\varpi $ naturally defines the required orientation of 
$H^{2+}(X; \mathbb{R})$.  

\item If $M$ does not split $X$, introduce 
$X_{1}$ to denote the complement of a tubular neighborhood of $M$ 
in $X$. In this case, 
\begin{equation*}
\underline {\SW}_{X} = j(\underline {\SW}_{X_{1}}),
\end{equation*}
where $j$ is the $\mathbb{Z}$--linear extension of the map from 
$H^{2}(X_{1},\partial X_{1}; \mathbb{Z})$ to $H^{2}(X, \partial
X ; \mathbb{Z})$ which  
arises by coupling the excision isomorphism with the natural homomorphism 
from the long exact cohomology sequence of the pair $M  \subset  X$. 
In addition the orientation for the line $L_{X}$ is induced 
by a chosen orientation of the analogous line for $X_{1}$. Finally, 
if $X$ is compact and $b^{2+} = 1$, then $\varpi $ 
naturally defines the needed orientation for $H^{2+}(X; \mathbb{R})$.
\end{itemize}
\end{thm}

As remarked at the outset, this theorem has a somewhat more general version 
which is given as Theorem \ref{thm:2.7}. The latter differs from Theorem 
\ref{thm:1.1} in that 
it discusses the Seiberg--Witten invariants proper   
rather than their averages 
over $2$--torsion classes. In any event, Theorem \ref{thm:1.1} follows 
directly 
as a  corollary to Theorem \ref{thm:2.7}.

When $X$ in Theorem \ref{thm:1.1} has the form $S^{1} \times  Y$, where 
$Y$ is a 
$3$--manifold, then the Seiberg--Witten invariants of $X$ are the same 
as those 
that are defined for $Y$ by counting solutions of a 
$3$--dimensional version of 
the Seiberg--Witten equations. In this case, Theorems \ref{thm:1.1} and 
\ref{thm:2.7} imply 
Mayer--Vietoris theorems for the $3$--dimensional Seiberg--Witten invariants. 
In 
particular, the $3$--dimensional version of Theorem \ref{thm:1.1} is stated as 
Theorem 
5.2 in \cite{MT}. As noted above, Lim \cite{L} has a proof of the 
$3$--dimensional 
version of Theorem \ref{thm:2.7}.

Before ending this introduction, a two part elipsis   
is in order which may or 
may not (depending on the reader) put some perspective on the subsequent 
arguments which lead back to Theorem \ref{thm:1.1}. 

Part one of this elipsis addresses, in a sense, a raison d'etre for Theorem 
\ref{thm:1.1}. To start, remark that the Seiberg--Witten invariants, like the 
Donaldson 
invariants \cite{A}, \cite{W2}, follow the `topological field theory' 
paradigm where 
Mayer--Vietoris like results are concerned. To elaborate: According to the 
topological field theory paradigm, the solutions to the $3$--dimensional 
version of the Seiberg--Witten equations associate a vector space with inner 
product to each $3$--manifold; and then the Seiberg--Witten equations on a 
$4$--manifold with boundary are expected to supply a vector in the boundary 
vector space. Moreover, when two $4$--manifolds with boundary are glued 
together across identical boundaries to make a compact, boundary free 
$4$--manifold, the field theory paradigm has the Seiberg--Witten invariants of 
the latter equal to the inner product of the corresponding vectors in the 
boundary vector space. 

Now, the fact is that the topological field theory paradigm is stretched 
somewhat when boundary $3$--tori are present. However, in the situation at 
hand, which is to say when there is a $2$--dimensional cohomology class with 
non-zero restriction in the cohomology of each boundary component, the 
paradigm is not unreasonable. In particular, the relevant boundary vector 
space is $1$--dimensional and so the topological field theory paradigm 
predicts 
that the Seiberg--Witten invariants of the $4$--manifolds with boundary under 
consideration here are simply numbers. And, when two such $4$--manifolds are 
glued across identical boundaries, then the Seiberg--Witten invariants of the 
result should be the product of the invariants for the pair. This last 
conclusion is, more or less, exactly what Theorem \ref{thm:1.1} states.

By the way, the essentially multiplicative form of the Mayer--Vietoris gluing 
theorems in \cite{MST} have an identical topological field theoretic 
`explanation.' 

Part two of this elipsis concerns the just mentioned \cite{MST} paper. The 
latter 
produced a Mayer--Vietoris gluing theorem for certain Seiberg--Witten 
invariants of a $4$--manifold cut along the product of a genus two or more 
surface and a circle. In particular, \cite{MST} considers only $Spin^{\C}$
structures $s$ whose class $c(s)$ evaluates on the surface  
to give two less than twice the genus; and \cite{MST} states a gluing theorem 
which is the genus greater than one analog of those given here. The case of 
genus one was not treated in \cite{MST} because the genus one case requires 
some 
special arguments. This paper gives a part of the genus one story. 
Meanwhile, other aspects of the genus one cases, Dehn surgery like gluings 
in particular, can be handled using the results in \cite{MMS}. 

By the way, a version of the 
gluing theorem in the surface genus greater than one context 
of \cite{MST} is 
used in \cite{MS}. 

The introduction ends with the list that follows of the section headings.

\begin{enumerate}\small
\item Introduction

\item The Seiberg--Witten invariants
	\begin{enumerate}
	\item The differential equation
	\item A topology on the set of solutions
	\item The structure of $\mm$
	\item Compactness properties   
	\item The definition of the Seiberg--Witten invariants
	\item Invariance of the Seiberg--Witten invariants
	\item The Mayer--Vietoris gluing theorems
	\end{enumerate}

\item Preliminary analysis
	\begin{enumerate}
	\item Moduli spaces for $T^{3}$
	\item Fundamental lemmas
	\item Immediate applications to the structure of $\mm$
	\item The family version of Proposition \ref{thm:2.4}
	\item Gluing moduli spaces
	\item Implications from gluing moduli spaces
	\end{enumerate}

\item Energy and compactness
	\begin{enumerate}
	\item The first energy bound
	\item Uniform asymptotics of $(A,\psi)$
	\item Refinements for the cylinder
	\item Vortices on the cylinder
	\item The moduli space for $\R\times T^{3}$
	\item Compactness in some special cases
	\end{enumerate}

\item Refinements for the cylinder
	\begin{enumerate}
	\item The operator $\dd_{c}$ when $X = \R \times  T^{3}$
	\item Decay bounds for $\kernel(\dd_{c})$ when $c \in  \mm_{P}$
	\item More asymptotics for solutions on a cylinder
	\item The distance to a non-trivial vortex
	\end{enumerate}

\item Compactness
	\begin{enumerate}
	\item Proof of Proposition \ref{thm:2.4}
	\item Proof of Proposition \ref{thm:3.7}
	\item Proof of Proposition \ref{thm:3.9}
	\end{enumerate}

\item $3$--dimensional implications

\end{enumerate}

This work was supported in part by the National Science Foundation.

\section{The Seiberg--Witten invariants}\label{sec:2}

This section provides a review of the definition of the Seiberg--Witten 
invariants for compact $4$--manifolds, and then extends the definition in 
Theorem \ref{thm:2.5} to cover the cases which are described in the 
introduction. It 
ends with the statement of Theorem \ref{thm:2.7}, which is the principle 
result of 
this article.

\sh{a) The differential equation}

In what follows, $X$ is an oriented, Riemannian $4$--manifold which can be 
non-compact. But, if the latter is the case, assume that there is a compact 
$4$--manifold with boundary $X_{0} \subset X$ whose 
boundary, $\partial X_{0}$, is 
a disjoint union of $3$--tori and whose complement is isometric to the half 
infinite cylinder $[0, \infty)  \times  \partial X_{0}$. 
To be precise, the metric on 
the $[0, \infty)$  factor should 
be the standard metric and the metric on $\partial X_{0}$ 
should be a flat metric. (Unless stated to the contrary, all metrics under 
consideration on $T^{3}$ will be flat.) The letter `$s$' is used below to 
denote a fixed function on $X$ which restricts 
to $[0, \infty)  \times  \partial X_{0}$ 
as the standard Euclidean coordinate on the factor $[0, \infty)$.

With the metric given, a $\Spin^{\C}$ structure is simply a lift (up to 
obvious equivalences) to a $\Spin^{\C}(4)$ principal bundle of the bundle 
$Fr\to X$ of oriented, orthonormal frames in the tangent bundle to $X$. Let 
$\sss_{0}(X_{0})\subset\sss(X)$ denote the subset of 
$\Spin^{\C}$ structures 
$s$ with $c({s}) = 0$ on $\partial X_{0}$. Choose a $\Spin^{\C}$ structure 
$s\in\sss_{0}(X_{0})$.

Associated to $s$'s principal $\Spin^{\C}(4)$ bundle $F\to X$ are a 
pair of the $\C^{2}$ vector bundles $S_{\pm}\to X$ as well as the 
complex line bundle $K\to X$. Here, $S_{+}$ arises from the group 
homomorphism which sends $\Spin^{\C}(4) =
(SU(2)\times SU(2) \times 
U(1))/\{\pm 1\}$ to $U(2) = (SU(2)\times U(1))/\{\pm 1\}$ by 
forgetting the first $SU(2)$ factor. Meanwhile, $S_{-}$ arises from the 
homomorphism to $U(2)$ which forgets the second factor; and $K$ arises by 
forgetting both factors. Note that $K = \det(S_{+}) = \det(S_{-})$ and the 
first Chern class of $K$ is the class $c({s})$.

Fix a self-dual $2$--form $\omega $ on $X$ which is non-zero and covariantly 
constant on each component of $[0, \infty) \times \partial X_{0}$. 
This is to say 
that the restriction of $\omega $ to such a component has the form 
$\omega=ds\wedge\theta+\omega_{0}$, where $\omega_{0}$ is a 
non-zero, covariantly constant $2$--form on $T^{3}$ and $\theta $ 
is the metric 
dual to $\omega_{0}$. 

Consider now the set of smooth configurations $(A,\psi)$ consisting of a 
connection, $A$, on $\det(S_{+})$ and a section $\psi $ of $S_{+}$ which solve 
the equations
\begin{equation}
\begin{split}
\bullet\text{ }&P_{+}F_{A}=\tau (\psi \otimes \psi^{\dagger}) - i\cdot 
\omega;\\
\bullet\text{ }&D_{A}\psi=0\label{eq:2.1};\\
\bullet\text{ }&\int_{X} |F_{A}|^{2}<\infty.
\end{split}
\end{equation}
The notation in \eqref{eq:2.1} is as follows: 

\begin{itemize}
\item $F_{A}$ denotes the curvature $2$--form of the connection $A$.

\item $P_{+}$ denotes the metric's orthogonal projection from the 
bundle of $2$--forms to the bundle, $\Lambda_{+}$, of self dual $2$--forms. 
(The latter is associated to the bundle $Fr\to X$ via the representation 
from $SO(4) = (SU(2) \times SU(2))/\{\pm 1\}$ to $SO(3) = SU(2)/\{\pm 
1\}$ which forgets the first $SU(2)$ factor. There is, of course, the bundle, 
$\Lambda_{-}$, of anti-self dual $2$--forms that is obtained via the 
representation to $SO(3)$ which forgets the second $SU(2)$ factor.) 

\item  $\tau $ denotes the homomorphism from $\End(S_{+}) = S_{+}  \otimes
S_{+}^{*}$ which is the hermitian adjoint to the Clifford multiplication 
homomorphism from $\Lambda_{+}$ into $\End(S_{+})$. 

\item  $D_{A}$ denotes a version of the Dirac operator. In particular, 
$D_{A}$ is the first order, elliptic operator which sends a section of 
$S_{+}$ to one of $S_{-}$ by composing a certain $A$--dependent covariant 
derivative on $S_{+}$ with the Clifford multiplication endomorphism from 
$S_{+} \otimes T^*X$ to $S_{-}$. Here, the covariant 
derivative is defined from the 
connection on $F$ which is obtained 
by coupling the connection $A$ with the 
pull-back from $Fr$ of the metric's Levi--Civita     
connection. 

\item In the last point of \eqref{eq:2.1}, the norm and the implicit volume 
form are defined by the given Riemannian metric.
\end{itemize}

A certain algebraic count of the solutions to \eqref{eq:2.1} gives the 
Seiberg--Witten 
invariants.

\sh{b) A topology on the set of solutions}

The set of solutions to \eqref{eq:2.1} is topologized as follows: Fix a base 
connection $A_{b}$ on $\det(S_{+})$ which is flat on $[0, \infty) \times 
\partial X_{0}$. With $A_{b}$ fixed, the set of connections on $\det(S_{+})$ 
can be 
identified with the space of smooth, imaginary valued $1$--forms, $i\cdot 
\Omega^{1}$. With the preceding understood, 
the space of solutions to 
\eqref{eq:2.1} 
is topologized by its embedding in the Fr\'echet 
space $i\cdot\Omega^{1}\oplus 
C^{\infty}(S_{+})\oplus\mathbb{R}$ 
which sends $(A, \psi )$
to $(A - A_{b}, 
\psi,\int_{X}|F_{A}|^{2})$. In this regard, the vector 
spaces $\Omega ^{1}$ and $C^{\infty}(S_{+})$ are topologized by the weak 
$C^{\infty}$ topology in which a typical neighborhood of $0$ is the space of 
sections which are small in the $C^{k}$ topology for some finite $k$ on some 
compact subset of $X$.

Note that the group $C^{\infty}(X;S^{1})$ acts continuously on the space of 
solutions to \eqref{eq:2.1} if this group is given 
the weak $C^{\infty}$ Fr\'echet structure 
in which a pair of maps are close if they are $C^{k}$ close for some finite 
$k$ on some compact subset of $X$. Here, 
$\varphi\in C^{\infty}(X;S^{1})$ 
sends a pair $(A,\psi)$ to $(A - 2\varphi^{-1}d\varphi,\varphi\psi)$. 
The quotient of the space of solutions by this action 
(with the quotient topology) will be called the \textit{moduli space} 
of solutions to \eqref{eq:2.1}. An orbit of $C^{\infty}(X;S^{1})$ 
will be called a `gauge orbit' and two solutions on the same gauge orbit 
will be deemed `gauge equivalent.' Except where confusion appears likely, a 
pair $(A,\psi)$ and its gauge orbit will not be notationally 
distinguished.

Before embarking on a detailed discussion of the structure of the moduli 
space of solutions to \eqref{eq:2.1}, 
some remarks are in order which concern an 
important consequence of the constraint given by the third point in 
\eqref{eq:2.1}. 
In particular, if $A$ is any connection on $K$, 
then up to factors of $2\pi$ and $i=\sqrt{-1}$, the 
curvature, $F_{A}$, is a closed $2$--form on $X$ whose cohomology class 
gives $c({s})$, the first Chern class of $K$. Now, by assumption, 
$c({s})$ restricts to zero on the ends of $X$ and thus lies in the image 
of the natural homomorphism 
from $H^{2}(X_{0}, \partial X_{0}; \Z)$ in $H^{2}(X;\Z)$. 
And, if $F_{A}$ is square integrable on $X$, arguments to follow prove 
that $F_{A}$ canonically defines a preimage, $c_{A}$, of $c({s})$ in 
$H^{2}(X_{0}, \partial X_{0}$; \Z).

The construction of $c_{A}$ employs an application of the abelian version of 
Uhlenbeck's compactness theorem \cite{U} as follows: Use $s$ to denote the 
Euclidean coordinate on the half line factor of 
$[0,\infty)\times\partial X_{0}$. 
Then, for large $s_{0}\in[0, \infty)$, Uhlenbeck's theorem insures that for 
any $s > s_{0}$, the connection $A$ restricts to the cylinder $[s, s + 1] 
\times\partial X_{0}$ as $A^{s}_{f}+a_{s}$, where $A^{s}_{f}$ is a flat 
connection on $\partial X_{0}$ and where $a_{s}$ is an imaginary $1$--form 
on $[s,s+1]\times\partial X_{0}$. Moreover, 
according to Uhlenbeck's theorem the 
sequence, indexed by  $s\in[s_{0},\infty)$, of the $L^{2}_{1}$ norms of 
$a_{s}$ over the defining domains, $[s, s +1]\times\partial X_{0}$, converges 
to zero as $s$ tends to infinity. Now, with the preceding understood, fix a 
sequence, $\{\beta_{s}\}$ of `cut-off' functions on 
$[0, \infty)$ , indexed by $s\in[s_{0},\infty)$, 
with $\beta_{s}=1$ on $[0,s]$, $\beta_{s}=0$ 
on $[s+1,\infty)$ and $|\beta_{s}'| < 2$ everywhere. Then, for $s > s_{0}$, 
introduce the connection $A^{s}$ on $K$ that equals 
$A^{s}_{f} + \beta_{s} a_{s}$ on $[s, \infty)\times\partial X_{0}$ 
and $A$ everywhere else. By 
construction, the curvature $2$--form of $A^{s}$ is zero on 
$[s+1,\infty) \times\partial X_{0}$. 
Thus, this $2$--form gives a bona fide 
class in the relative 
cohomology group $H^{2}(X,[s+1,\infty)\times\partial X_{0};\Z)$. And, as the 
latter group is canonically isomorphic to $H^{2}(X_{0},\partial X_{0};\Z)$, 
the 
curvature $2$--form of $A^{s}$ defines a class 
in this last group as well. When 
$s$ is sufficiently large, the latter class is the desired $c_{A}$. 

Of course, this definition of $c_{A}$ makes sense provided that the $s\in 
[s_{0}, \infty)$ indexed set of curvature $2$--forms 
$\{F_{A^s}\}$ give identical 
classes in $H^{2}(X_{0},\partial X_{0}$;\Z) when $s$ is large. 
To prove that such 
is the case, consider the classes indexed by some pair $s$ and $s +\delta$ 
with $\delta\in(0, 1)$. The corresponding curvatures both define 
integral classes in the relative group 
$H^{2}(X,[s+\delta+1)\times\partial X_{0};\Z)$ and it is 
sufficient to prove that these two classes 
agree. In particular, since $A^{s}$ and $A^{s+\delta}$ agree on 
$X-([s,\infty)\times\partial X_{0})$ and since $H^*(T^{3};\Z)$ 
has no torsion, such will 
be the case if the composition of exterior product and then integration over 
$X$ pairs both curvature $2$--forms identically with a given set of closed 
$2$--forms on $X$ that generate $H^{2}(X;\mathbb{Z})$. 
And, for this purpose, it is 
enough to take a generating set of forms which are covariantly constant on 
$[0,\infty)\times\partial X_{0}$. 

With these last points understood, it then follows that the relevant 
cohomology classes agree if the curvature forms in question are close in the 
$L^{2}$ sense on\linebreak
$[s,s+2]\times\partial X_{0}$ since these curvature forms 
agree on the complement of $[s,\infty)\times\partial X_{0}$. Of course, these 
forms are $L^{2}$ close (when $s$ is large) since each separately has small 
$L^{2}$ norm on 
$[s,s+2]\times\partial X_{0}$ by virtue of the fact that 
the $1$--forms $a_{s}$ and $a_{s+\delta}$ have small $L^{2}_{1}$ norms on 
this same cylinder.

With only minor modifications, the preceding argument that $c_{A}$ is well 
defined yields the following:

\begin{lem}\label{thm:2.1} 
Let $\aaa$ denote the space of smooth 
connections on $K$ whose curvature $2$--form is square integrable, here 
endowed with the smallest topology for which the assignment of $A$ 
to $\int_{X}|F_{A}|^{2}$ is continuous and 
which allows $C^{\infty}$ convergence on compact sets. Then, 
the assignment of $c_{A}$ to $A\in\aaa$ defines a 
locally constant function on $\aaa$. In particular, the connected 
components of the moduli space $\mm({s})$ are labeled, in part, 
by the set of elements in $H^{2}(X_{0},\partial X_{0};\Z)$ which map 
to $c({s})$ under the natural homomorphism to $H^{2}(X;\mathbb{Z})$. 
\end{lem}

With Lemma \ref{thm:2.1} in mind and supposing that $s\in\sss_{0}(X)$ 
and 
a preimage, $z$, of $c({s})$ in $H^{2}(X_{0},
\partial X_{0};\Z)$ have been 
given, introduce 
$\mm\equiv\mm({s},z)$ to denote the subspace of    
pairs 
$(A,\psi)$ in the moduli space 
of solutions to the $s$ version of 
\eqref{eq:2.1} for which $c_{A} = z$. 

\sh{\boldmath c) The structure of $\mm$}

The local structure of $\mm$ 
is described in the next proposition. However, the 
statement of this proposition requires a preliminary digression to point out 
certain topological features of $X$. To start the digression, note that there 
is an integer valued, bilinear pairing on $H^{2}(X_{0},\partial X_{0};\Z)$ 
which is obtained by composing the cup-product map to $H^{4}(X_{0}, 
\partial X_{0};\Z)$ 
with evaluation on the fundamental class. In contrast to the 
case where $\partial X_{0} = \varnothing$, 
this form has a null space in the non-empty 
boundary case, that being the image of $H^{1}(\partial X_{0})$ via the natural 
connecting homomorphism of the long exact sequence for the pair $(X_{0}, 
\partial X_{0})$. Thus, the cup product pairing is both well defined and 
non-degenerate on the image in $H^{2}(X_{0})$ of $H^{2}(X_{0}, 
\partial X_{0})$. 
Use $z\bullet z'$ to denote the cup product pairing between 
classes $z$ and $z'$. Meanwhile, use 
$H^{2+}(X_{0},\partial X_{0};\mathbb{R})\subset
H^{2}(X_{0},\partial X_{0};\mathbb{R})$ to denote a maximal dimensional vector 
subspace on which this cup product pairing is positive definite and use 
$b^{2+}(X_{0})$ to denote the dimension of $H^{2+}(X_{0},\partial X_{0}; 
\mathbb{R})$. 
Also, use $\tau$ to denote the signature of the cup product pairing. Thus, 
$\tau = b^{2+}-b^{2-}$, where $b^{2-}$ is the dimension of the 
maximal vector subspace in $H^{2}(X_{0},\partial X_{0})$ where the cup product 
pairing is negative definite. Finally, the digression ends by introducing 
$b^{1}_{0}$ to denote the dimension of $H^{1}(X_{0},\partial X_{0}; 
\mathbb{R})$. 

Here is the promised local structure result: 

\begin{prop}\label{thm:2.2} 
Let $c \equiv (A, \psi) \in\mm$. 
Then, there exists a Fredholm operator $\dd_{c}$ of index 
$d\equiv b^{1}_{0}-1-b^{2+}+4^{-1}(c({s})\bullet 
c({s})-\tau)$; a real analytic map $f$, from a 
ball in the kernel of $\dd_{c}$ to the cokernel of $\dd_{c}$ 
mapping the origin in the ball to the origin in the cokernel of 
$\dd_{c}$; and, provided that $\psi$ is not identically 
zero, a homeomorphism, $\varphi $, from $f^{-1}(0)$ onto 
an open neighborhood of $c$ in $\mm$. 
\end{prop}
Note that when $\partial X_{0}\ne\varnothing$, 
there are no solutions to \eqref{eq:2.1} where $\psi$ 
is identically zero. 

For $c=(A,\psi)$, the operator $\dd_{c}$ is a differential operator which 
is initially defined to send $(b,\eta)\in iC^{\infty}(T^*X)\oplus 
C^{\infty}(S_{+})$ to the element in  $iC^{\infty}(\R)\oplus 
iC^{\infty}(\Lambda_{+})\oplus C^{\infty}(S_{-})$ 
whose components in the three summands are as 
follows:
\begin{equation}
\begin{split}
&d^*b-2(\psi^{\dagger}\eta - \eta^{\dagger}\psi);\\
&P_{+}db - \tau'(\eta\oplus\psi^{\dagger}+\psi\oplus\eta^{\dagger});\\
&D_{A}\eta+\clif(b)\psi.\label{eq:2.2}
\end{split}
\end{equation}
Here, $\tau'$ denotes the polarization of the bilinear form $\tau $ which 
appears in \eqref{eq:2.1} 
and $\clif(\cdot)$ denotes the Clifford multiplication 
endomorphism from $T^*X$ to $\Hom(S_{+},S_{-})$. To make $\dd_{c}$ Fredholm, a 
preliminary domain and range are defined to allow only sections with compact 
support. This preliminary domain is then completed using the Sobolev 
$L^{2}_{1}$ norm, while the preliminary range is completed using the 
Sobolev $L^{2}$ norm. 

Under favorable conditions, the local neighborhoods described in Proposition 
\ref{thm:2.2} 
fit nicely together to give $\mm$ the structure of a smooth $d$--dimensional 
manifold. The following proposition elaborates:

\begin{prop}\label{thm:2.3} 
With reference to the previous 
proposition, the set of points in $\mm$ where the cokernel of the 
operator $\dd_{c}$ is $\{0\}$ has the structure of a smooth, 
$d$--dimensional manifold whereby the homeomorphism $\varphi$ 
is a smooth coordinate chart. In addition, this last portion of $\mm$ 
is orientable and canonically so with a choice of orientation for 
the line 
$\Lambda^{\mathrm{top}}H^{1}(X_{0},\partial X_{0};\mathbb{R})\otimes\Lambda 
^{\mathrm{top}}H^{2+}(X_{0},\partial X_{0};\mathbb{R})$. Finally, if 
$\partial X_{0}\ne\varnothing$ 
or if $b^{2+} > 0$, then there exists a Baire subset of 
choices for the $2$--form $\omega$ in \eqref{eq:2.1} for which the 
corresponding $\mm$ is everywhere a smooth manifold. In fact, given an 
open set $U$ in $X_{0}$ whose closure is disjoint from 
$\partial X_{0}$, and a self-dual form $\omega'$ on 
$X_{0}$ that has the required structure near $\partial X_{0}$, 
there is a Baire set of smooth, self-dual extensions, $\omega$, of 
$\omega'$ from $X_{0} - U$ to the whole of 
$X_{0}$ for which this same conclusion holds. For such $\omega$,  
$\mm = \varnothing$ when $d < 0$; and when $d \ge 0$, 
then the cokernel of $\dd_{c}$ is trivial for every $c\in\mm$.
\end{prop}
The proofs of these last two propositions are given in Section \ref{sec:3} 
of this 
paper.

In what follows, a point $c\in\mm$ will be called a `smooth point' when the 
cokernel of the operator $\dd_{c}$ is trivial.

\sh{d) Compactness properties}

In the case where $X_{0}$ has no boundary, the moduli space $\mm$ is compact; 
this compactness is one of the remarkable features of the Seiberg--Witten 
equations. However, even the simplest example with non-empty boundary, 
$T^{2}\times D^{2}$, can yield non-compact moduli spaces. Even so, 
certain zero and $1$--dimensional subspaces of $\mm$ are compact if the form 
$\omega $ is suitably chosen. 

A two part digression follows as a preliminary to the specification of the 
constraints on $\omega$.

\textbf {Part 1}\qua Remember that $X_{0}$ is assumed to have a 
class $\varpi\in H^{2}(X_{0};\mathbb{R})$ which is non-zero in the 
cohomology of each component of $\partial X_{0}$. 
Meanwhile, as the chosen $2$--form 
$\omega$ 
is constant on each component of $\partial X_{0}$, it defines a cohomology 
class, $[\omega]\in H^{2}(\partial X_{0};\mathbb{R})$. 
With this understood, say 
that $\omega$ is \textit{tamed} by $\varpi$ when $[\omega] = \varpi$ 
in $H^{2}(\partial X_{0};\mathbb{R})$. 

\textbf {Part 2}\qua Introduce $\varsigma({s}) 
\subset H^{2}(X_{0},\partial X_{0};\Z)$ to denote the set of elements which 
map to $c({s})$ in $H^{2}(X_{0};\Z)$. Next, introduce the sets
\begin{equation}\begin{split}\label{eq:2.3}
&\mm_{s} \equiv \bigcup_{z\in\varsigma(s)}\mm({s},z)\text{ and}\\
&\mm_{s,m} \equiv \bigcup_{z\in\varsigma(s):z\bullet\varpi\le m} 
\mm({s}, z)\subset\mm_{s}.
\end{split}\end{equation}
Endow $\mm_{s}$ with the topology of $C^{\infty}$ 
convergence on compact subsets of 
$X$ and give $\mm_{s,m}$ the subspace topology. This is the topology which 
arises by embedding the space of solutions to \eqref{eq:2.1} 
in the Fr\'echet space 
$i\cdot\Omega^{1}\oplus C^{\infty}(S_{+})$ with the latter given the 
$C^{\infty}$ weak topology. Any $\mm({s},z)\subset\mm_{s}$ will be 
called a \textit{stratum} of $\mm_{s}$.

With $\mm_{s}$ and each $\mm_{s,m}$ understood, here is the most that can be 
said at this point about compactness:

\begin{prop}\label{thm:2.4} 
Let $\varpi\in H^{2}(X_{0};\mathbb{R})$ 
be a class with non-zero pull-back to the cohomology of each 
component of $\partial X_{0}$. With $\varpi$ given, use a form 
$\omega$ in \eqref{eq:2.1} that is tamed by $\varpi$. 
Then each $\mm_{s,m}\subset\mm_{s}$ is compact and 
contains only a finite number of strata. Moreover, fix a 
self-dual form $\omega'$ that is non-zero and covariantly constant 
on each component of $[0,\infty)\times\partial X_{0}$ and that is tamed 
by $\varpi$; and fix a non-empty, open set $U\subset X_{0}$. 
Then, there is a Baire set of smooth, self-dual forms $\omega$ that 
agree with $\omega'$ on $X - U$ and have the 
following properties:
\begin{itemize}
\item As in Proposition \ref{thm:2.3}, each stratum of $\mm_{s}$ 
is a smooth manifold of dimension 
$d$ given in Proposition \ref{thm:2.2}.
Moreover, the cokernel of the operator
$\dd_{c}$ vanishes for each $c\in\mm_{s}$.
\item The boundary of the closure in $\mm_{s}$ of any 
stratum intersects the remaining strata as a codimension $2$ submanifold.
\end{itemize}
\end{prop}
Roughly said, Proposition \ref{thm:2.4} 
guarantees the compactness of the zero set in 
$\mm({s},z)$ of a reasonably 
chosen section of a $d$ or $(d-1)$--dimensional    
vector bundle over $\mm_{s,m}$. 

By the way, it turns out that an extra cohomology condition on the class 
$\varpi$ guarantees the compactness of the whole of each $\mm({s},z)$. 
Indeed, this pleasant situation arises when the restriction of $\omega $ to 
each component of $\partial X_{0}$ defines a cohomology 
class which is not a linear 
multiple of an integral class. Proposition \ref{thm:4.6} 
gives the formal statement.

Proposition \ref{thm:2.4} is proved in Section \ref{sec:6}a.

\sh{e) The definition of the Seiberg--Witten invariants}

The simplest version of the Seiberg--Witten invariant for $X_{0}$ associates 
an integer to a pair $s\in\sss_{0}(X)$ and $z\in H^{2}(X_{0}, 
\partial X_{0};\Z)$ mapping to $c({s})$. This integer will be denoted by 
$\sw({s},z)$. As in the empty boundary case, it is obtained via an 
algebraic count of the elements in $\mm$. However, there are some additional 
subtleties when $\partial X_{0}\ne\varnothing$ because $\mm$ 
need not be compact. 

In what follows, $X_{0}$ is as described above except that positivity of 
$b^{2+}$ will be implicitly assumed when 
$\partial X_{0}=\varnothing$. Fix a $\Spin^{\C}$ 
structure $s\in\sss_{0}(X_{0})$ and a class 
$z\in\varsigma({s})\subset H^{2}(X_{0},\partial X_{0};\Z)$. Note 
that $s$ provides the integer $d$ in Proposition \ref{thm:2.2}. 
Also, fix a class 
$\varpi\in H^{2}(X_{0};\mathbb{R})$ which is non-zero in the cohomology of 
each component of $\partial X_{0}$. Finally, orient 
$L\equiv\Lambda^{\mathrm{top}}H^{1}(X_{0},\partial X_{0};\mathbb{R})\otimes 
\Lambda^{\mathrm{top}}H^{2+}(X_{0},\partial X_{0};\mathbb{R})$ and, when 
$\partial X_{0}=\varnothing$ and $b^{2+} = 1$, orient 
$H^{2+}(X_{0};\mathbb{R})$.

What follows is the definition of $\sw$.

\textbf{Case 1}\qua This case has either $d < 0$ or $d$ odd. Set 
$\sw({s}, z) = 0$ in this case.

\textbf{Case 2}\qua This case has $d = 0$. Choose a form 
$\omega$ 
in \eqref{eq:2.1} which is tamed by $\varpi$ 
and which is such that each stratum 
of $\mm_{s}$ has the structure described in Propositions \ref{thm:2.3} 
and \ref{thm:2.4}. The 
latter insure that $\mm = \mm({s},z)$ is a finite set of points. In 
addition, each point $c\in\mm$ comes with a sign, 
$\varepsilon(c)\in\{\pm 1\}$, 
from the orientation. With these points understood, set 
\begin{equation}\label{eq:2.4}
\sw({s},z)\equiv\sum_{c}\varepsilon(c),
\end{equation}
where the sum is taken over all $c\in\mm$.

\textbf{Case 3}\qua This case has $d > 0$ and even. Once again, 
choose a form $\omega$ in \eqref{eq:2.1} 
which is tamed by $\varpi$ and which is 
such that $\mm_{s}$ has the structure described in Propositions \ref{thm:2.3} 
and \ref{thm:2.4}. 
Thus, each stratum of $\mm_{s}$ is an oriented, $d$--dimensional manifold. 

Next, choose a set, $\Lambda\subset X$, of $d/2$ distinct 
points, and for each $x\in\Lambda$, specify a $\C$--linear surjection 
$\Phi_{x}\colon S_{+}\big|_{x}\to\C$. Use $\underline\Lambda$ to 
denote the resulting set of $d/2$ pairs $(x,\Phi_{x})$. With 
$\underline\Lambda$ understood, set 
\begin{equation}\label{eq:2.5}
\mm^{\underline\Lambda}\equiv\{c=(A,\psi)\in\mm : 
\Phi_{z}(\psi(x))=0\text{ for each }x\in\underline\Lambda\}.
\end{equation}
Note that $\mm^{\underline\Lambda}$ can 
be viewed as the zero set of a smooth section 
of a $d/2$--dimensional complex vector bundle over $\mm$. 
This understood, Sard's 
theorem guarantees that $\mm^{\underline\Lambda}$ is discrete 
for a Baire set of data 
$\underline\Lambda$, and each $c\in\mm^{\underline\Lambda}$ comes with a 
sign, $\varepsilon(c)\in\{\pm 1\}$. 
Moreover, Proposition \ref{thm:2.4} 
guarantees that this Baire set can be found so that
the corresponding 
$\mm^{\underline\Lambda}$ is a finite set.

With $\underline\Lambda$ now chosen from the afore mentioned Baire set 
of possibilities, define $\sw({s},z)$ by \eqref{eq:2.4} but with the sum 
restricted to those $c$ in the set $\mm^{\underline\Lambda}$.

When $X$ is compact, there also exists an extension of $\sw$ whose image is in 
$\Lambda^*H^{1}(X_{0};\Z)\equiv\Z\oplus 
H^{1}\oplus\Lambda^{2}H^{1}\oplus\cdots$. 
The latter is described in \cite{T3} 
and the definition there can be readily adapted to the non-compact setting 
described here. Theorems \ref{thm:2.5} and \ref{thm:2.7} 
below have reasonably self-evident 
analogs which apply to this extended $\sw$. 
However, to prevent an already long 
paper from getting longer, the extended version of $\sw$ will not be discussed 
further here. Thus, the statements of the versions of Theorems 
\ref{thm:2.5} and \ref{thm:2.7} 
that apply to the extended $\sw$ are left to the reader to supply. 

\sh{f) Invariance of the Seiberg--Witten invariants}

With $\sw(\cdot)$ so defined, there is an obvious question to address: To 
what extent does $\sw(\cdot)$ depend on the various choices that enter its 
definition? 

In the case where $X_{0}$ is compact, the following answer is well known 
(the arguments are given in \cite{W1}, but see also \cite{M}, \cite{KKM}):
{\sl
\begin{itemize}
\item If $b^{2+} > 1$, then $\sw$ is 
independent of the choice of Riemannian metric and form $\omega$; 
its absolute value depends only on the $\Spin^{\C}$ 
structure, and the sign is determined by the orientation of the line 
$L\equiv\Lambda^{\mathrm{top}}H^{1}(X_{0};\mathbb{R})\otimes 
\Lambda^{\mathrm{top}}H^{2+}(X_{0};\mathbb{R})$. Moreover, 
$\sw(\varphi^*s)=\sw(s)$ when $\varphi$ is a diffeomorphism of 
$X_{0}$ which preserves the orientation of the line $L$.

\item  If $b^{2+}=1$, first specify an orientation 
of $H^{2+}(X_{0};\mathbb{R})$. Then, $\sw$ is independent of the 
choice of Riemannian metric and form $\omega$ provided that the 
integral over $X_{0}$ of the wedge of $\omega $ with an 
oriented harmonic representative of $H^{2+}(X_{0};\mathbb{R})$ is 
sufficiently large and positive. So defined, the absolute value of $\sw$ 
only depends on the $\Spin^{\C}$ structure and the 
orientation of $H^{2+}(X_{0},\R)$, and its sign is determined by 
the orientation of the line $L$. Furthermore, $\sw(\varphi^*{s}) = \sw({s})$ 
when $\varphi$ is a 
diffeomorphism of $X_{0}$ which preserves the orientations of the 
line $L$ and $H^{2+}(X_{0};\mathbb{R})$. 
\end{itemize}
}

The next result provides an answer to the opening question in the case where 
$\partial X_{0}$ is not empty.

\begin{thm}\label{thm:2.5} 
Suppose that $\partial X_{0}\ne\varnothing$. First, 
choose a class $\varpi\in H^{2}(X_{0};\mathbb{R})$ which is 
non-zero in the cohomology of each component of $\partial X_{0}$. Next, 
define $\sw$ on a pair $(s,z)$ using $\varpi$ 
as described in the preceding subsection. Then, the 
result is independent of the chosen metric and the form $\omega$ 
provided that the latter is tamed by $\varpi$. Here, the absolute 
value of $\sw$ is determined solely by the triple $(s,z,\varpi)$ 
and the sign is determined by the chosen orientation for 
the line $L\equiv\Lambda^{\mathrm{top}}H^{1}(X_{0},\partial X_{0};\mathbb{R})
\otimes 
\Lambda^{\mathrm{top}}H^{2+}(X_{0},\partial X_{0};\mathbb{R})$. Moreover, if 
$\varphi$ is a diffeomorphism of $X_{0}$ which fixes the 
orientation of the line $L$, then the value of $\sw$ on 
$\varphi^*(s,z,\varpi)$ is the same as its value on 
$(s,z,\varpi)$. Finally, $\sw$ is insensitive to 
continuous deformation of $\varpi$ in $H^{2}(X_{0};\mathbb{R})$ 
through classes with non-zero restriction in the cohomology of each 
component of $\partial X_{0}$. 
\end{thm}

The proof of Theorem \ref{thm:2.5} is provided in Section \ref{sec:3}d.

\sh{g) The Mayer--Vietoris gluing theorems}

The purpose of this subsection is to state the advertised generalization of 
the Mayer--Vietoris gluing result given by Theorem \ref{thm:1.1}. 
This generalization 
is summarized in Theorem \ref{thm:2.7}, 
below, but a four-part digression comes first 
to set the stage.

\textbf{Part 1}\qua In what follows, $X_{0}$ is a compact, 
oriented $4$--manifold 
with boundary such that each component of $\partial X_{0}$ is a 
$3$--torus. Suppose next that there is an embedded $3$--torus 
$M\subset X_{0}$ which separates $X_{0}$ so that $X_{0} = X_{+}\cup X_{-}$, 
where $X_{\pm}$ are $4$--manifolds with boundary 
embedded in $X$ which intersect 
in $M$.

With the preceding set up understood, introduce the lines $L_{0}$, $L_{+}$ 
and $L_{-}$ via
\begin{equation}\label{eq:2.7}
L_{\diamondsuit}\equiv\Lambda^{\mathrm{top}}H^{1}(X_{\diamondsuit}, 
\partial X_{\diamondsuit};\mathbb{R})\otimes\Lambda^{\mathrm{top}} 
H^{2+}(X_{\diamondsuit},\partial X_{\diamondsuit};\mathbb{R}), 
\end{equation}
where $\diamondsuit$ is a stand in for $0$, $+$ or $-$. An argument from 
\cite{MST} 
can be adapted almost verbatim to establish the existence of a canonical 
isomorphism
\begin{equation}\label{eq:2.8}
L_{0}\approx L_{+}\otimes L_{-}.
\end{equation}
Thus, orientations of $L_{+}$ and $L_{-}$ canonically induce an orientation 
of $L_{0}$.

If $M\subset X_{0}$ is non-separating, introduce $X_{1}$ to denote the 
complement in $X_{0}$ of a tubular neighborhood of $M$. Then, the 
afore-mentioned argument from \cite{MST} adapts readily to establish the 
existence of a canonical isomorphism between $L_{0}$ from \eqref{eq:2.7} 
and $L_{1}\equiv\Lambda^{\mathrm{top}}H^{1}(X_{1},\partial X_{1};\mathbb{R})
\otimes\Lambda 
^{\mathrm{top}}H^{2+}(X_{1},\partial X_{1};\mathbb{R})$. 

\textbf{Part 2}\qua By assumption, there is a class 
$\varpi\in H^{2}(X_{0};\mathbb{R})$ whose 
pull-back is not zero in the cohomology of 
each component of $\partial X_{0}$. Theorem \ref{thm:2.7}, 
below, will assume that the 
pull-back of $\varpi$ to the cohomology of $M$ is also non-zero. With this 
understood, then the pull-back of $\varpi$ will be non-zero in the 
cohomology of each component of $\partial X_{+}$ 
and each component of $\partial X_{-}$ in 
the case when $M\subset X$ is separating. Likewise, when $M$ is not 
separating, then the pull-back of $\varpi $ in the cohomology of each 
component of $\partial X_{1}$ will be non-zero.

By the way, in the case when $X_{0}$ is compact and has $b^{2+} = 1$, the 
choice of a class $\varpi\in H^{2}(X_{0};\mathbb{R})$ whose pull-back to 
the cohomology of $M$ is non-zero supplies an orientation for 
$H^{2+}(X_{0};\Z)$. Indeed, because $\varpi\ne0\in H^{2}(M;\Z)$, 
there is a class in $H_{2}(M;\Z)$ whose push-forward in $H_{2}(X_{0};\Z)$ 
is non-zero and which pairs positively with $\varpi $. This homology class 
has self-intersection number zero, so its image in $H^{2}(X_{0};\Z)$ lies 
on the `light cone.' Thus, the latter's direction specifies an orientation 
to any line in $H^{2}(X_{0};\mathbb{R})$ on which the cup-product pairing is 
positive definite.

\textbf{Part 3}\qua This part of the digression contains the 
instructions for the construction of a $\Spin^{\C}$ structure on $X$ from what 
is given on $X_{\pm} $ or $X_{1}$. To start, consider a somewhat abstract 
situation where $Y$ is a smooth, oriented $4$--manifold and $U\subset Y$ is 
any set. Having defined $\sss(Y)$ 
as in the introduction, define $\sss(U)$ to denote 
the equivalence class of pairs $(Fr\big|_{U},F_{U})$, where $Fr\big|_{U}$ is a 
principal $SO(4)$ reduction of the restriction of the oriented, general linear 
frame bundle of $Y$ to $U$, and where $F_{U}$ is a lift of $Fr\big|_{U}$ to a 
principal $\Spin^{\C}(4)$ bundle. This definition provides a tautological 
pull-back map $\sss(Y)\to\sss(U)$ 
which intertwines the action of $H^{2}(Y;\Z)$ 
with that of its image in $H^{2}(U;\Z)$. 

With the preceding understood, let $\sss_{0M}(X_{0})\subset 
\sss_{0}(X_{0})$ denote the set of $\Spin^{\C}$ 
structures whose image under $c$ 
in \eqref{eq:1.2} 
is zero under pull-back to the cohomology of $M$. When $M$ separates 
$X_{0}$, then the pull-back map from the preceding paragraph defines a map 
$\wp^{0}\colon\sss_{0M}(X_{0})\to\sss_{0}(X_{-})\times 
\sss_{0}(X_{+})$. Meanwhile, in the case where $M$ is non-separating, there is 
the analogous $\wp^{0}\colon\sss_{0M}(X_{0})\to\sss_{0}(X_{1})$. In 
this regard, note that $\sss_{0M}(X_{0})$ is a principal
homogeneous space 
for the image in $H^{2}(X_{0};\Z)$ of 
$H^{2}(X_{0},M;\Z)$ and the map 
$\wp$        
intertwines the action of the latter 
group with its image in either 
$H^{2}(X_{-};\Z)\times H^{2}(X_{+};\Z)$ or $H^{2}(X_{1};\Z)$ as 
the case may be.

With $\wp^{0}$ understood, the question arises as to the sense in which 
it can be inverted. The answer requires the introduction of some 
additional terminology. For this purpose, let $Y$ be a compact, oriented 
$4$--manifold with boundary which is a disjoint union of tori. Introduce 
$\sss_{0}(Y,\partial Y)$ to denote the set of pairs $(s,z)$ where ${s} 
\in\sss_{0}(Y)$ and where $z\in H^{2}(Y,\partial Y;\Z)$ maps to $c({s}) 
\in H^{2}(Y;\Z)$ under the long exact sequence homomorphism. Note that 
$\sss_{0}(Y,\partial Y)$ is a principal 
homogeneous space for the group $H^{2}(Y,\partial Y; 
\Z)$. Perhaps it is needless to say that there is a tautological `forgetful' 
map from $\sss_{0}(Y,\partial Y)$ 
to $\sss_{0}(Y)$ which intertwines the action of 
$H^{2}(Y,\partial Y;\Z)$ with that of its image in $H^{2}(Y;\Z)$.

With the new terminology in hand, consider:

\begin{lem}\label{thm:2.6} 
Depending on whether $M$ does or does not 
separate $X_{0}$, there is a canonical map, $\wp$, from 
$\sss_{0}(X_{-},\partial X_{-})\times\sss_{0}(X_{+},\partial X_{+})$ 
or $\sss_{0}(X_{1},\partial X_{1})$ into $\sss_{0}(X_{0},\partial X_{0})$ 
respectively, which has the following properties:
\begin{itemize}
\item The image of $\wp$ lies in $\sss_{0M}(X_{0},\partial X_{0})$
\item $\wp$ either intertwines the action of 
$H^{2}(X_{-},\partial X_{-};\Z)\times H^{2}(X_{+},\partial X_{+};\Z)$ 
or that of $H^{2}(X_{1},\partial X_{1};\Z)$, as the case may 
be, with their images in $H^{2}(X_{0},\partial X_{0};\Z)$.
\item The composition of $\wp$ and then $\wp^{0}$ 
gives the canonical forgetful map. 
\end{itemize}
\end{lem}

\begin{proof}[Proof of Lemma \ref{thm:2.6}] 
What follows is the argument for the case when 
$M$ separates $X$. The argument for the other case is analogous and is left to 
the reader.

To start, remark that the given $\Spin^{\C}$ structures ${s}_{\pm}$ 
can be patched together over $M$ with the specification of an isomorphism over 
$M$ between the corresponding lifts, $F_{\pm}\to Fr$. In this regard, 
note that the choice of a Riemannian metric on $X$ which is a product flat 
metric on a tubular neighborhood, $U\approx I\times M$, of $M\subset X$ 
determines principal $SO(4)$ reductions of the general linear 
frame bundles of $X_{\pm}$ which are consistent with the inclusions of 
$X_{\pm}$ in $X$. 

Having digested the preceding, note next that the space of isomorphisms 
between $F_{+}\big|_{M}$ and itself which cover the identity on 
$Fr\big|_{M}$ has 
a canonical identification with the space of maps from $M$ to the circle; thus 
the space of isomorphisms $\varphi\colon F_{+}\big|_{M}\to F_{-}\big|_{M}$ 
which cover the identity on $Fr\big|_{M}$ has a non-canonical identification 
with $C^{\infty}(M;S^{1})$. This 
implies that the set of homotopy classes of 
such maps is a principal bundle over a point for the the group $H^{1}(M;\Z)$. 
In this regard, note that a pair of isomorphisms between $F_{+}\big|_{M}$ 
and $F_{-}\big|_{M}$ yield the same $\Spin^{\C}$ structures over $X$ 
if and only 
if they differ by a map to $S^{1}$ which extends over either $X_{+}$ or 
$X_{-}$. 

In any event, a choice of isomorphism from $F_{+}\big|_{M}$ to 
$F_{-}\big|_{M}$ 
covering $Fr\big|_{M}$ is canonically equivalent to a choice of isomorphism 
between the restrictions to $M$ of the associated $U(1)$ line bundles 
$K_{\pm}$. Meanwhile, as $c({s}_{\pm}) = 0$, the data $z_{+}\in 
H^{2}(X_{+},\partial X_{+}; \Z)$ mapping to $c({s}_{+})$ canonically 
determines a homotopy class of isomorphisms from $K_{+}\big|_{M}$ to 
$M\times\C$. Likewise, $z_{-}$ determines a homotopy class of isomorphisms 
from $K_{-}\big|_{M}$ to $M\times\C$. With the preceding understood, use 
the composition of an isomorphism 
$K_{+}\big|_{M}\approx M\times\C$ 
in the $z_{+}$ determined                 
class with the inverse of one between 
$K_{-}\big|_{M}$ to $M\times\C$ from the $z_{-}$ determined class to 
construct the required isomorphism between $F_{+}\big|_{M}$ and 
$F_{-}\big|_{M}$. 
\end{proof}                       

\textbf{Part 4}\qua Lemma \ref{thm:2.6} makes the point that the image 
of $\wp$ contains only those $\Spin^{\C}$ structures on $X$ whose image under 
the map $c$ pulls back as zero to $H^{2}(M;\Z)$. There may well be other 
$\Spin^{\C}$ structures on $X$. Even so, a case of the main theorem in 
\cite{KM} 
asserts that $\sw({s}) = 0$ if $c({s})$ does not pull back as zero 
in $H^{2}(M;\Z)$.

The digression is now over, and so the stage is set for the main theorem:

\begin{thm}\label{thm:2.7} 
Let $X_{0}$ be a compact, connected, 
oriented $4$--manifold with (possibly empty) boundary consisting of a 
disjoint union of $3$--dimensional tori such that restriction to each 
boundary component induces a non-zero pull-back map on the second 
cohomology. If the boundary is empty, require that $b^{2+}\ge1$. 
Let $M\subset X_{0}$ be an embedded $3$--dimensional 
torus for which the restriction induced pull-back homomorphism from 
$H^{2}(X_{0};\mathbb{R})$ to $H^{2}(M;\mathbb{R})$ is non-zero. 
Choose a class $\varpi\in H^{2}(X_{0};\mathbb{R})$ whose 
pull-back in the cohomology of $M$ and in that of every component of 
$\partial X_{0}$ is non-zero. If $M$ splits $X_{0}$ as a 
pair, $X_{-}\cup X_{+}$, of $4$--manifolds with boundary, then 
orient the lines $L_{\pm}$ and then orient the corresponding line 
$L_{0}$ via \eqref{eq:2.8}. Otherwise, orient the line $L_{1}$ and 
use the isomorphism $L_{1}\approx L_{0}$ to orient the 
latter. If $X_{0}$ has empty boundary and $b^{2+} > 1$, 
use the orientation for $L_{0}$ to define the map $\sw$ on 
$\sss(X_{0})$. If $X_{0}$ has empty boundary and $b^{2+} = 
1$, use the orientation for $L_{0}$ and that for 
$H^{2+}(X_{0};\mathbb{R})$ as defined by $\varpi$ to define $\sw$ 
on $\sss(X_{0})$. Finally, if $X_{0}$ has non-trivial 
boundary, use the orientation for $L_{0}$ and the class $\varpi$ 
to define $\sw$ on $\sss(X_{0},\partial X_{0})$.  
\begin{itemize}
\item Suppose that $M$ splits $X_{0}$ as a 
pair, $X_{-}\cup X_{+}$, of $4$--manifolds with boundary. 
Use the chosen orientations for the lines $L_{\pm}$ and 
the restriction of $\varpi$ to $X_{\pm}$ to define the 
corresponding maps $\sw\colon\sss_{0}(X_{\pm},\partial X_{\pm})\to\Z$. 
Then, for all $({s}, z)\in\sss_{0}(X_{0},\partial X_{0})$, 
there are just finitely many pairs $(({s}_{-},z_{-}), 
({s}_{+},z_{+}))\in\wp^{-1}(({s},z))$ 
with either $\sw(({s}_{-},z_{-}))$ or 
$\sw(({s}_{+},z_{+}))$ non-zero; and with this 
last fact understood,
\begin{equation*}
\sw(({s}, z)) = \sum_{((s_{-},z_-),(s_{+},z_{+}))\in{\wp}^{-1}((s,z))} 
\sw(({s}_{-},z_{-}))\sw(({s}_{+},z_{+})).
\end{equation*}
\sloppy
\item If $M$ does not split $X$, use 
the chosen orientations for the line
$L_{1}$ and the restriction 
of $\varpi$ to $X_{1}$ to define the corresponding map 
$\sw\colon\sss_{0}(X_{1},\partial X_{1})\to\Z$. Then, for each 
$({s},z)\in\sss_{0}(X_{0},\partial X_{0})$, there are just 
a finite number of $(s_{1},z_{1})\in\wp^{-1}(({s},z))$ with 
$\sw(({s}_{1},z_{1}))\ne0$; and with this 
last point understood,
\begin{equation*}
\sw(({s}, z)) = \sum_{((s_{1},z_{1}))\in\wp^{-1}((s,z))}
\sw(({s}_{1},z_{1})).
\end{equation*}
\fussy
\end{itemize}
\end{thm}

The proof of Theorem \ref{thm:2.7} is given in Section \ref{sec:3}f.

\section{Preliminary analysis}\label{sec:3}

The proofs of Theorems \ref{thm:2.5} and \ref{thm:2.7} 
use many of the ideas from \cite{MST}, but new 
techniques are also involved. The new techniques enter into the proof of 
Proposition \ref{thm:2.4} 
and into the proofs of related compactness assertions which 
concern the moduli spaces for manifolds with long cylinders that are 
products of $3$--tori with intervals. These related compactness assertions are 
summarized below in Propositions \ref{thm:3.7} and \ref{thm:3.9}. 

This section sees to the separation of these compactness aspects of the 
proofs of Theorems \ref{thm:2.5} and 
\ref{thm:2.7} from the more well known techniques. In so 
doing, it explains how these theorems follow from Propositions 
\ref{thm:2.4}, \ref{thm:3.7} and 
\ref{thm:3.9} 
while leaving the proofs of the latter to the subsequent sections of 
this paper. The details start in Subsection \ref{sec:3}a below. A guide to the
analytical points that arise in this and the subsequent sections immediately 
follows. 

Any attempt to define an `invariant' via an integer weighted count of 
solutions to an equation must deal with the following two absolutely central 
issues: First, there must be some 
guarantee of a finite count. Second, to insure 
the invariance of the count, solution 
appearance and disappearance with change of movable parameters must occur in 
groups with zero aggregate count. Given a suitable topology on the solution 
set, both of these issues are issues of compactness. Indeed the former 
concerns the compactness of the solution set for some fixed parameter value 
while the latter concerns the compactness of a family of solutions spaces as 
determined by a corresponding family of parameter values. 

The investigation of this compactness issue starts with the next two points. 
The first is a standard application of 
elliptic regularity theory and the 
second follows from a Bochner--Weitzenb\"ock 
formula for the Seiberg--Witten 
equations.

{\sl
\begin{itemize}
\item A bound on the $L^{2}$ norms of $F_{A}$ and 
$\psi$ on a ball implies uniform $C^{\infty}$ estimates for 
some pair on the gauge orbit of $(A,\psi)$ on the 
concentric, half-radius ball. Hence, 
the space of gauge orbits of solutions 
that satisfy an a priori $L^{2}$ norm 
on a ball is precompact in the 
concentric, half-radius ball.

\item For the toroidal end manifolds under consideration, the 
$\Spin^{\C}$ structure determines a uniform $L^{2}$ bound 
for both $F_{A}$ and $\psi$ on any ball when the 
perturbing form $\omega$ in \eqref{eq:2.1} has the properties stated in 
Proposition \ref{thm:2.4}.
\end{itemize}}

Although quite powerful, the preceding points are not powerful enough to 
imply compactness for $\mm$ 
since they do not foreclose 
leakage down the ends of 
$X$. The characterization of this leakage requires an investigation of 
solutions to \eqref{eq:2.1} 
on finite and infinite cylinders. This investigation 
begins with the derivation of the following key result: 

{\sl
\begin{itemize}
\item If $(A, \psi)$ solves \eqref{eq:2.1} on a given cylindrical 
portion of $\R\times T^{3}$ where $\omega$ is 
non-zero and constant, and if $F_{A}$ has small $L^{2}$ 
norm on each subcylinder of length $4$ in the given cylinder, then $|F_{A}|$ 
decays exponentially from both ends of the given cylinder.
\end{itemize}}
This last fact, as confirmed in the initial subsections below, is ultimately 
a consequence of the structure of a moduli space of a certain version of the 
Seiberg--Witten equations on $T^{3}$. 

As argued in Section \ref{sec:4}, 
the preceding point with the first two implies: 

{\sl
\begin{itemize}
\item All limits of non-convergent sequences in $\mm$ are 
described by data sets that consist of a second solution to \eqref{eq:2.1} 
on $X$ 
and also a finite set (with topological bound on its size) of 
solutions on $\R\times T^{3}$. 
\end{itemize}}
Given an assertion that precludes the extra solutions on $\R\times T^{3}$, 
this last point implies that $\mm$ is compact. Moreover, certain 
geometric assumptions about the tubular end $4$--manifold and the form 
$\omega$ actually do imply the absence of the relevant $\R\times T^{3}$ 
solutions.

When these just mentioned geometric assumptions are not met, then the 
argument for compactness employs the following additional observations:

{\sl
\begin{itemize}
\item The relevant moduli space of solutions to \eqref{eq:2.1} on 
$\R\times T^{3}$ is a smooth manifold with a vector bundle whose 
fiber at each point is the cokernel of a Fredholm operator associated to the 
solution in question. This vector bundle has positive dimension and comes 
with a canonical, nowhere zero section.
\item This canonical section defines a canonical, 
$\R^{2}$--valued function on every moduli space of solutions to \eqref{eq:2.1} 
over $X$. There is one such function for each end of $X$. 
\item One or more of these canonical functions vanishes at any 
solution on $X$ that appears in the data set for a non-convergent 
sequence in $\mm$.
\item If $\omega$ is suitably generic on the 
interior of $X_{0}$, then the zero set in each moduli space of each 
of these canonical functions is a codimension $2$ submanifold.
\end{itemize}}

These last points are used in Section \ref{sec:6} 
to establish the assertions of 
Proposition \ref{thm:2.4} 
and those of the forthcoming Propositions \ref{thm:3.7} and 
\ref{thm:3.9}. These points are established in Section \ref{sec:5} using 
a detailed analysis of the first order operator that is obtained by 
linearizing the Seiberg--Witten equations about a solution on the cylinder 
$\R\times T^3$.

As illustrated by the preceding discussion, the moduli spaces on $T^{3}$ and 
$\R\times T^{3}$ are central to a significant chunk of the compactness 
story and so they are the focus of much of the subsequent discussion.

\sh{\boldmath a) Moduli spaces for $T^{3}$}

The story behind Theorems 
\ref{thm:2.5} and \ref{thm:2.7} starts here with a description of the 
moduli spaces of a version of the Seiberg--Witten equations on a flat, 
oriented $3$--torus. These moduli spaces enter into the story because they are 
naturally identified with the moduli spaces of translationally invariant 
solutions to a version of \eqref{eq:2.1} 
on $\R\times T^{3}$. In any event, the 
equations considered here require the choice of a flat metric on $T^{3}$ 
plus a lift of the resulting $SO(3)$ frame bundle to a principal 
$\Spin^{\C}(3) = U(2)$ 
bundle. Such lifts are classified up to isomorphism by the first 
Chern class of the complex line bundle, $K$, which is associated to the 
determinant representation of $U(2)$ on $\C$. 
Note that this class, $c$, is an even 
class in $H^{2}(T^{3};\Z)$. The equations will also require the choice of 
a covariantly constant $2$--form, $\omega_{0}$, on $T^{3}$.

It is worth digressing here momentarily to comment some on the relationship 
between the $3$ and $4$--dimensional stories. To start the digression, keep in 
mind that the tori that arise in this article come as constant `time' slices 
of an oriented $\R\times T^{3}$. In addition, the $\R$ factor will come 
oriented and thus induce an orientation on $T^{3}$; this will be the 
implicit orientation of choice.

With compatible orientations for $\R$, $T^{3}$ and $\R\times T^{3}$ 
understood, the set of $\Spin^{\C}$ structures on $\R\times T^{3}$ has a 
canonical, 1-1 correspondence with the set of lifts of the $SO(3)$ principal 
frame bundle to a $U(2)$ bundle. Indeed, this correspondence comes about via a 
natural map from the set of isomorphism classes of $U(2)$ lifts of the $SO(3)$ 
frame bundle of $T^{3}$ to $\sss(\R\times T^{3})$. What follows is the 
definition of this map.

To describe the aforementioned map, note first that a $U(2)$ lift, $P$, of the 
$SO(3)$ frame bundle comes with an associated, principal $\Spin^{\C}(4)$ 
bundle, 
$F_{P}$, which will be viewed both as a bundle over 
$T^{3}$ and, via pull-back, 
as one over $\R\times T^{3}$. In this regard, $F_{P}$ is defined using 
the representation which sends a pair $(h,\lambda)\in U(2) = (SU(2) 
\times S^{1})/\{\pm 1\}$ to $(h,h,\lambda)\in\Spin^{\C}(4)$. 
This associated $\Spin^{\C}(4)$ bundle is a lift of the pull-back to $\R 
\times T^{3}$ of the analogously defined, associated principal $SO(4)$ 
bundle to the frame bundle of $T^{3}$. Meanwhile this last $SO(4)$ bundle is 
canonically isomorphic to the frame bundle of $\R\times T^{3}$. Indeed, 
a choice of oriented, unit length tangent vector field to the $\R$ factor 
induces just this isomorphism. 

Thus, in the manner just described, an isomorphism class of principal $U(2)$ 
lifts of the $SO(3)$ frame bundle of $T^{3}$ canonically determines an element 
${s}_{P}\in\linebreak\sss(\R\times T^{3})$. The inverse of this map 
starts with ${s}\in\sss(\R\times T^{3})$ and the corresponding 
$\Spin^{\C}(4)$ bundle $F$. The oriented, unit length tangent vector to $\R$ 
defines the reduction of $F$ to a principal $U(2)$ 
bundle, $P$, as follows: First, 
define $Fr^{3}\subset Fr$ as the subset of frames whose first basis 
element is metrically dual to the chosen vector on the $\R$ factor. Then, view 
$F$ as a bundle over $Fr$ 
(with fiber $S^{1}$) and use $P\subset F$ to denote 
restriction of this bundle to $Fr^{3}$. This $P$ is a principal $U(2)$ bundle 
which covers the $SO(3)$ frame bundle of $T^{3}$. 

In any event, let $P\to T^{3}$ denote a principal $U(2)$ lift of the 
$SO(3)$ frame bundle. Introduce $S$ to denote the complex $2$--plane 
bundle $P\times_{U(2)}\C^{2}$. The Seiberg--Witten equations on $T^{3}$ are 
equations for a pair $(A,\psi)$, where $A$ is a connection on the complex 
line bundle $\Lambda^{2}S = P\times_{U(2)}\C$ and where $\psi $ 
is a section of $S$. These equations read:
\begin{equation}\begin{split}\label{eq:3.1}
\bullet\text{ }&F_{A} = 
\tau (\psi \otimes \psi^{\dagger}) - i\cdot \omega_{0};\\
\bullet\text{ }&D_{A}\psi=0. 
\end{split}\end{equation}
In this last equation, $\tau$ denotes the homomorphism from 
$\End(S) = S \otimes S^*$
which is the hermitian adjoint to the Clifford multiplication homomorphism 
from $\Lambda^{2}(T^{3})$ into $\End(S)$ while $D_{A}$ denotes a version 
of the Dirac operator. In particular, $D_{A}$ is the first order, elliptic 
operator which sends a section of $S$ to another section of $S$ by composing a 
certain $A$--dependent covariant derivative on $S$ with the Clifford 
multiplication endomorphism from $S \otimes T^*X$ to $S$. 
Here, (and below) the 
covariant derivative is defined from the connection on $F$ which is obtained 
by coupling the connection $A$ with the pull-back from $Fr$ of the metric's 
Levi--Civita connection. 

A second digression is in order here concerning the relationship between the 
$3$ and $4$ 
dimensional versions of Clifford multiplication. To start, introduce 
the $\Spin^{\C}(4)$ principal bundle $F_{P}\to\R\times T^{3}$ 
which corresponds to $P$, and then introduce the associated $C^{2}$ bundles 
$S_{\pm}\to\R\times T^{3}$. Both are canonically isomorphic to $S = 
P\times_{U(2)}\C^{2}$. In this way, the $4$--dimensional Clifford 
multiplication endomorphism from $T^*(\R\times T^{3})\times S_{+}\to 
S_{-}$ induces a Clifford multiplication map $T^*(T^{3})\times S 
\to S$. Here, $T^*(T^{3})$ is viewed as a summand in $T^*(\R\times T^{3})$ 
via the pull-back monomorphism from the projection map to 
$T^{3}$. By the way, with $S_{+}$ and $S_{-}$ identified as $S$, Clifford 
multiplication by the $1$--form $dt$ is just multiplication by $i$.

With the digression now over, introduce $\mm_{P}$ to denote the moduli space 
of solutions to \eqref{eq:3.1} 
for a given flat metric, covariantly constant form 
$\omega_{0}$ and lift $P$ of the $SO(3)$ frame bundle. Thus, $\mm_{P}$ is the 
quotient of the space of smooth pairs $(A,\psi)$ which solve \eqref{eq:3.1} 
by the 
action of the group $C^{\infty}(T^{3};S^{1})$. Here, the action is the same 
as in the $4$--dimensional case.  Note that $\mm_{P}$ is given the quotient 
topology.

The lemma that follows describes the salient features of $\mm_{P}$: 

\begin{lem}\label{thm:3.1} 
Given the flat metric, there exists 
$\delta > 0$ such that if $|\omega_{0}| < \delta$, 
then the space $\mm_{P}$ is empty unless 
$c_{1}(S) = 0$. In addition, given that $c_{1}(S) = 0$ 
and $\omega_{0}\ne0$, then $\mm_{P}$ is a single 
point, the orbit of a pair $(A_{0},\psi_{0})$, where $A_{0}$ 
is a trivial connection, $\psi_{0}$ is covariantly 
constant and $\tau(\psi_{0} \otimes \psi_{0}^{\dagger}) = i 
\omega_{0}$. When $\omega_{0} = 0$, then $\mm_{P}$ 
consists of the orbits of those pairs $(A, 0)$, where 
A is a flat connection; thus $\mm_{P} = H^{1}(T^{3}; 
\R^{3})/H^{1}(T^{3};\Z^{3})$. 
\end{lem}

\begin{proof}[Proof of Lemma \ref{thm:3.1}] 
The proof sticks closely to a well worn trail 
initially blazed by Witten in \cite{W1} 
and translated to the $3$--manifold context 
at the start of Section 5 of \cite{MST}. In fact, the argument follows almost 
verbatim the discussion in the proof of \cite{MST}'s 
Lemma 5.1. Here is a brief 
synopsis of the argument: First of all, in the case where $\omega_{0}= 
0$, use the Bochner--Weitzenb\"ock formula 
$D_{A}^{2}=\nabla_{A}^*\nabla_{A}            
\psi+2^{-1}\clif(F_{A})\cdot\psi$ on 
$T^{3}$, and the two Seiberg--Witten equations to conclude that $\nabla 
_{A}\psi=0$. Here, $\nabla_{A}$ denotes the covariant derivative 
on sections of $S$ which is defined by $A$ and the Levi--Civita 
connection on the 
$SO(3)$ frame bundle. As $\psi$ is covariantly constant, either $F_{A} = 0$ 
or $\psi = 0$, or both. In particular, the first point in \eqref{eq:3.1} 
is only 
consistent with both vanishing.

In the case where $\omega_{0}\ne0$, Clifford multiplication on $S$ by 
$\omega_{0}$ defines a skew Hermitian, covariantly constant endomorphism 
of $S$. Decompose $S$ into the eigenbundles for this endomorphism and write 
$\psi $ with respect to this decomposition as $(\alpha,\beta)$. In so 
doing, follow the steps in \cite{W1} 
or at the beginning of Section 5 in \cite{MST}. 
Here, it may be useful to consider lifting the story to $S^{1}\times 
T^{3}$ with $S^{1}$ viewed as $\R/\Z$ to make the connection with the 
$4$--dimensional framework in these references. In any event, with $\psi $ 
written as $(\alpha,\beta)$, the second equation in \eqref{eq:3.1} becomes a 
coupled system of equations for the pair $(\alpha,\beta)$. Continuing 
the analysis in either \cite{W1} or \cite{MST} 
then leads directly to the conclusion 
that $\beta= 0$ and that when $|\omega|$ is small, then $\alpha $ is 
constant and $F_{A} = 0$.
\end{proof}                       

An orbit in $\mm_{P}$ is termed either a smooth point or not. These are 
technical terms which correspond to whether or not there are non-trivial 
deformations of pairs on the given orbit which solve \eqref{eq:3.1} 
to first order. A 
precise definition is given in \cite{MST} 
subsequent to Lemma 5.2 and continued 
in \cite{MST}'s 
section 5.1. As is standard in gauge theory problems, the notion 
of being smooth or not is based on whether or not a certain Fredholm 
operator has vanishing cokernel. Here, the operator in question is obtained 
by first linearizing the equations in \eqref{eq:3.1} 
about a given solution, then 
restricting the domain to the $L^{2}$--orthogonal complement of the tangents 
to the orbit of the $C^{\infty}(T^{3};S^{1})$ action, and finally projecting 
the resulting expression onto an isomorphic image of this same orthogonal 
complement. In this regard, the tangents to the orbit here appear as the 
image of the operator denoted by $D_{x}$ in Section 5 of \cite{MST}. 
Alternately, 
one can view the equations in \eqref{eq:3.1} 
as an $S^{1}$ invariant version of \eqref{eq:2.1} 
on $S^{1}\times T^{3}$ in which case the operator under scrutiny here 
is the $S^{1}$ invariant version of the operator in \eqref{eq:2.2}. 
In any event, 
consider:

\begin{lem}\label{thm:3.2}
Suppose that $P$ is a principal $U(2)$ 
lift of the frame bundle of $T^{3}$ for which $c_{1}(S) = 
0$. If $\omega_{0}\ne0$, then the point in 
$\mm_{P}$ is a smooth point.
\end{lem}

\begin{proof}[Proof of Lemma \ref{thm:3.2}] 
This is a straightforward computation since the 
operators involved have constant coefficients. The details are left to the 
reader.
\end{proof}

Remark that $\mm_{P}$ has no smooth points when $\omega_{0} = 0$. Moreover, 
in this case, the solution which corresponds to the trivial connection has a 
larger space of first order deformations than do the others.

The remainder of this paper considers only the case where $P\to T^{3}$ 
is the `trivial lift,' that is, the lift of the $SO(3)$ frame bundle for which 
the associated $C^{2}$ bundle $S$ is topologically trivial. This assumption 
about P is made implicitly from here to the end of the paper.

\sh{b) Fundamental lemmas}

The preceding lemmas on the structure of $\mm_{P}$ can be used to deduce 
certain key facts about solutions to the Seiberg--Witten equations on the 
product of $T^{3}$ with an interval. The latter are summarized by the lemmas 
in this subsection. Here is the first: 

\begin{lem}\label{thm:3.3}
Fix a flat metric on $T^{3}$ and then 
fix a non-zero, constant $2$--form $\omega$ on $\R\times T^{3}$. 
The metric on $T^{3}$, the form $\omega$, a 
non-negative integer $k$ and a choice of $\varepsilon > 0$ 
determine $\delta > 0$ which has the following 
significance: Suppose that $(A,\psi)$ is a solution to 
\eqref{eq:2.1} on $Y\equiv(2, -2)\times T^{3}$ as 
defined by the product flat metric and $\omega$. If $\int 
_{Y}|F_{A}|^{2} < \delta$, then $(A,\psi)$ has 
$C^{k}$ distance $\varepsilon$ or less on the 
subinterval $[-1, 1] \times T^{3}$ from a point in the gauge 
orbit on $Y$ of $(A_{0},\psi_{0})\in\mm_{P}$.
\end{lem}
By way of explaining terminology, the statement that some $(A_{1},\psi_{1})$ 
is in the `gauge orbit on $Y$ of $(A_{0},\psi_{0})\in\mm_{P}$' 
means only that $(A,\psi)$ is gauge equivalent via some element 
in $C^{\infty}(Y;S^{1})$ to a pair of connection and section of $S_{+}$ which 
are the pull-backs from $T^{3}$ of a pair which solves \eqref{eq:3.1}. 

The second fundamental lemma can be viewed as a corollary to Lemmas 
\ref{thm:3.2} and 
\ref{thm:3.3}. Here is this second lemma:

\begin{lem}\label{thm:3.4}
Fix a flat metric on $T^{3}$ and then 
fix a non-zero, self-dual, constant $2$--form $\omega$ on $\R\times 
T^{3}$. The metric and form $\omega$ determine a 
constant $\delta > 0$ and a set of constants $\{\zeta 
_{k}\}_{k=0,1,\ldots}$ with the following significance: 
Suppose that $R\ge0$ and that $(A,\psi)$ is a 
solution to \eqref{eq:2.1} on $Y\equiv(-R-2,R+2)\times T^{3}$ 
as defined by the product metric and $\omega$. 
Suppose as well that $\int_{U}|F_{A}|^{2}<\delta$ 
for every length $2$ cylinder $U=(t -1, t + 1)\times 
T^{3}\subset Y$. Then, there is a pair $(A_{1},\psi_{1})$ 
in the gauge orbit on $Y$ of $(A_{0},\psi_{0})\in\mm_{P}$ such that
\begin{equation}\label{eq:3.2}
|\nabla^{k}(A - A_{1})| + |(\nabla_{A_{1}})^{k}(\psi-\psi_{1})| 
\le \zeta_{k}(e^{-\delta(R-t)}+e^{-\delta(R+t)})
\end{equation}
at any point $(t, x)\in[-R, R]\times T^{3}$.
\end{lem}

These two lemmas are proved shortly so accept them for now to consider one 
of their more immediate consequences, that elements in the moduli space $\mm$ 
of Sections \ref{sec:2}a--c
decay exponentially fast along the cylindrical ends of $X$.

\begin{lem}\label{thm:3.5} 
Let $X_{0}$ be the usual $4$--manifold 
with toroidal boundary components, and let $X = X_{0}\cup([0,\infty)  
\times\partial X_{0})$ denote the corresponding non-compact manifold 
with a metric which restricts to $[0,\infty)\times\partial X_{0}$ as a 
flat, product metric. Let $\omega$ denote a self-dual form on 
$X$ which is constant and non-zero on each component of $[0,\infty)  
\times\partial X_{0}$. This data determines $\delta > 0$, a 
sequence of constants $\{\zeta_{k}\}_{k=0,1,2,\ldots}$ and, with 
the choice of $r\ge1$, a constant $R$; and these constants 
have the following significance: Let $(s, z)\in\sss_{0}(X_{0}, 
\partial X_{0})$, let $\mm$ denote the 
resulting moduli space of solutions to \eqref{eq:2.1}, and let $(A, 
\psi)\in\mm$ obey $\int_{[r,\infty)\times\partial X_0} 
|F_A|^2 < \delta$. 
Then, on each component, $Y\subset[R,\infty)\times
\partial X_{0}$ there is a point $(A_{1},\psi_{1})$ on the 
gauge orbit on $Y$ of $(A_{0},\psi_{0})\in\mm_{P}$ for which
\begin{equation}\label{eq:3.3}
|\nabla^{k}(A-A_{1})| + |(\nabla_{A_{1}})^{k}(\psi - 
\psi_{1})|\le\zeta_{k}e^{-\delta(t-R)}
\end{equation}
at any point $(t, x)\in Y$.
\end{lem}

\begin{proof}[Proof of Lemma \ref{thm:3.5}] 
Write $t$ in \eqref{eq:3.2} as $t' + R$ and then take $R$ to 
infinity in the resulting equation to obtain bounds on $(A-A_{1},
\psi-\psi_{1})$ at points $(t',x)\in[0,\infty)\times T^{3}$. The 
latter are identical to those in \eqref{eq:3.3} after shifting $t$ 
in \eqref{eq:3.3} to $t' - R$.
\end{proof}

This subsection ends with the proofs of Lemmas \ref{thm:3.3} and 
\ref{thm:3.4}.

\begin{proof}[Proof of Lemma \ref{thm:3.3}] 
The $L^{2}$ bound on $P_{+}F_{A}$ by 
any $\delta\ge0$ immediately yields an $L^{2}$ bound on $|\psi|$ since
\begin{equation}\label{eq:3.4}
|\tau(\psi\otimes\psi^{\dagger})|=z|\psi|^{2} 
\end{equation}
with $z$ a universal constant. Thus, since $|\omega|$ is a constant, the 
first point in \eqref{eq:2.1} 
implies that there is a constant $z_{\omega}$ and a 
bound of the form
\begin{equation}\label{eq:3.5}
\int_{Y}\left||\psi|-z_{\omega}\right|^{4}\le z_{1}\delta,
\end{equation}
with $z_{1}$ depending only on $|\omega|$.

The next step obtains bounds on the $L^{2}$ norm of $\nabla_{A}\psi $, 
and the Bochner--Weitzenb\"ock formula for the Dirac operator is the principle 
tool for doing so. Without assumptions on the Riemannian metric, the 
connection $A$ and the section $\psi $, this formula reads
\begin{equation}\label{eq:3.6}
D_{A}^*D_{A}\psi=\nabla_{A}^*\nabla_{A}\psi+4^{-1}s\psi+2^{-1}\clif_{+}
(P_{+}F_{A})\cdot\psi.
\end{equation}
Here, $s$ denotes the metric's scalar curvature while $D_{A}^*$ and $\nabla 
_{A}^*$ denote the formal, $L^{2}$ adjoints of the Dirac operator $D_{A}$ 
and the covariant derivative $\nabla_{A}$. Also, $\clif_{+}(\cdot)$ 
denotes the Clifford multiplication induced homomorphism from $\Lambda_{+}$ 
into $\End(S_{+})$.

As indicated first by Witten in \cite{W1}, 
this formula in conjunction with the 
Seiberg--Witten equations can be used with great effect to analyze the 
behavior of solutions to \eqref{eq:2.1}. 
In particular, the left hand side of \eqref{eq:3.6} 
is zero when the second line in \eqref{eq:2.1} 
holds, while the first line can be 
used to control the term with $P_{+}F_{A}$. 

In any event, for the purposes at hand, take the inner product of both sides 
of \eqref{eq:3.6} 
with $\psi$ and use the equations in \eqref{eq:2.1} to rewrite the result 
as
\begin{equation}\label{eq:3.7}
2^{-1}d^*d|\psi|^{2}+|\nabla_{A}\psi|^{2}+2^{-1}|P_{+}F_{A}|^{2}
+2^{-1}\langle P_{+}F_{A},i\omega\rangle=0.
\end{equation}
Here, $d^*$ denotes the formal $L^{2}$ adjoint of the exterior derivative $d$. 

Equation \ref{thm:3.7} implies that the square of the $L^{2}$ norm of 
$|\nabla_{A}\psi|$ over $[-3/2, 3/2]\times T^{3}$ is bounded by $z_{2} 
\delta^{1/2}$, with $z_{2}$ only dependent on $|\omega|$. Indeed, to 
obtain such a bound, first replace $|\psi|^{2}$ by $(|\psi^{2}|-z_{\omega})$ 
in the first term on the left side. Then, multiply both sides 
of the resulting equation by a smooth, non-negative function which equals $1$ 
on $[-3/2, 3/2]\times T^{3}$ and vanishes near the boundary of $Y$. Next, 
integrate the result over $Y$ and then integrate by parts to remove the 
derivatives in $d^*d$ from $|\psi|^{2}-z_{\omega}$. Finally, an appeal 
to \eqref{eq:3.5} 
and a suitable application of the inequality $2|ab|\le\delta^{-1/2}
|a|^{2}+\delta^{1/2}|b|^{2}$ produces the asserted bound.

With the $L^{2}$ norms of $\psi$, $\nabla_{A}\psi $ and $F_{A}$ 
bounded in terms of $\delta$ over the interior cylinder $[-3/2, 3/2] 
\times  T^{3}$, the next step proves that there is, given $\varepsilon  > 
0$, a value for $\delta $ which implies that the $(A,\psi)$ has 
$L^{2}_{1}$ distance $\varepsilon $ or less from a pair $(A_{1},\psi 
_{1})$ in the gauge orbit on $Y$ of $(A_{0},\psi_{0})\in\mm_{P}$. 
This step is a straightforward argument by contradiction which 
invokes fairly standard elliptic techniques. In this regard, the only 
novelty is that the action of 
$C^{\infty}(Y;S^{1})$ must be used to make \eqref{eq:2.1} 
an elliptic system. Indeed, ellipticity can be achieved by using the 
$C^{\infty}(Y;S^{1})$ action to write 
$(A = A_{0} + b, \psi  = \psi_{0} + \eta)$,     
where $(a,\eta)$ are constrained 
so that the expression 
in the top line of \eqref{eq:2.2} 
holds. The details of all of this are left to the 
reader. 

Finally, given that $(A,\psi)$ is $L^{2}_{1}$ close to a point in the 
gauge orbit on $Y$ of $(A_{0},\psi_{0})$, the final step proves that $(A, 
\psi)$ is $C^{k}$ close to $(A_{0},\psi_{0})$. This last part of the 
proof is also left to the reader as it constitutes a direct application of 
standard procedures in elliptic regularity theory.
\end{proof}

\begin{proof}[Proof of Lemma \ref{thm:3.4}] 
The lemma can be proved by invoking, with only 
minor changes, the argument which \cite{MST} 
uses to prove its Corollary 6.17. 
However, a slightly more direct argument can made by filling out the sketch 
that follows. The sketch starts with the remark that $\varepsilon  > 0$ 
provides a positive upper bound for the $L^{2}$ norm over $U$ of $F_{A}$
that has the following significance: When the $L^{2}$ norm of $F_{A}$ is 
less than this bound, then $(A,\psi)$ is gauge equivalent to $(A_{0} + b, 
\psi_{0} + \eta)$, where the $C^{1}$ norms of $(b,\eta)$ are bounded 
by $\varepsilon $ on the sub-cylinder $Y'\equiv[-R-1, R+1]\times
T^{3}$; and where the Seiberg--Witten equations in terms of $\lambda   
\equiv(a,\eta)$ have the form
\begin{equation}\label{eq:3.8}
\partial_{t}\lambda + L_{0}\lambda + r(\lambda) = 0.
\end{equation}
Here, $\partial_{t}$ 
is the tangent vector field to the line segment factor in $Y$, 
and $L_{0}$ is a linear, symmetric, first order differential operator. 
Meanwhile, $r(\lambda)$ in \eqref{eq:3.8} 
is a `remainder' term which is formally 
second order in $\lambda $. To be precise here, the $L^{2}$ norm of 
$r(\lambda)$ on each constant $t\in[-R-1,R+1]$ slice of $Y$ is 
bounded by the product of $\varepsilon $ and the    
$L^{2}$ norm of $\lambda $ on 
the same constant $t$ slice. 

What follows are some important points about $L_{0}$. First, $L_{0}$ is 
determined solely by the metric on $T^{3}$, $\omega $ and $(A_{0},\psi_{0})$, 
and thus only by the metric and $\omega $. Second, the $L^{2}$ 
spectrum of $L_{0}$ is discrete, real, and lacks accumulation points. Third, 
$0$ is not in the spectrum of $L_{0}$. This last conclusion is essentially the 
statement of Lemma \ref{thm:3.2}
and plays the starring role in the subsequent part of 
the argument. 

With $L_{0}$ introduced, let $\lambda_{\pm} $ denote the projections of 
$\lambda $ onto the respective eigenspaces of $L_{0}$ with positive ($+$) and 
negative ($-$) eigenvalues. Next, introduce functions $f_{\pm} $ on the 
interval $[-R - 1, R + 1]$ whose values at a point $t$ are the respective 
$L^{2}$ norms of $\lambda_{\pm} $ on the corresponding constant $t$ slice 
of $Y$. Then \eqref{eq:3.8} yields the differential inequalities
\begin{equation}\begin{split}\label{eq:3.9}
&\partial_{t} f_{+} + (E - \varepsilon) f_{+} - \varepsilon f_{-} \le 0;\\
&\partial_{t} f_{-} - (E - \varepsilon) f_{-} + \varepsilon f_{+} \ge 0.
\end{split}\end{equation}
Here, $E > 0$ is the distance between $0$ 
and the spectrum of $L_{0}$.

With \eqref{eq:3.9} understood, a simple comparison argument establishes the 
following: When $\varepsilon\ll E$, then the inequalities in \eqref{eq:3.9} 
require 
both $f_{\pm} $ to decay exponentially from the ends of $[-R - 1, R + 1]$. 
This exponential decay for the $L^{2}$ norms of $(a,\eta)$ on the constant 
$t$ slices of $Y'$ can then be 
bootstrapped to give \eqref{eq:3.2}
using standard 
elliptic regularity techniques.
\end{proof}

\sh{\boldmath c) Immediate applications to the structure of $\mm$}

Lemmas \ref{thm:3.1}--\ref{thm:3.5} 
have certain automatic consequences with regard to the 
moduli spaces which are considered in Section \ref{sec:2}c. 
In particular, 
Propositions \ref{thm:2.2} and \ref{thm:2.3} follow from these lemmas.

\begin{proof}[Proof of Proposition \ref{thm:2.2}] 
The argument here for the structure of $\mm$ 
is completely analogous to that derived in \cite{T1} for the $SU(2)$ self-dual 
moduli spaces on manifolds with cylindrical ends. In this regard, observe 
that the lower two components of the image of $\dd_{c}$ in \eqref{eq:2.2} 
are nothing 
more than the linearization of the equations in the first two points of 
\eqref{eq:2.1}. 
Meanwhile, the vanishing of the first component of the image of 
$\dd_{c}$ in \eqref{eq:2.2} 
only asserts that the given section of $i\cdot T^*X\oplus 
S_{+}$ is $L^{2}$--orthogonal to the space of tangents to the orbit of the 
gauge group. The fact that $\dd_{c}$ is Fredholm with the $L^{2}_{1}$ domain 
and $L^{2}$ range follows from Lemma \ref{thm:3.2} 
by standard arguments; for example, 
by invoking Lemma \ref{thm:3.5} 
to control the behavior of $(A,\psi)$ on $[0, \infty)  
\times  \partial X_{0}$, the fact that $\dd_{c}$ 
is Fredholm follows almost directly 
from results in \cite{APS}.

The formula in Proposition \ref{thm:2.2} 
for the index of $\dd_{c}$ can be derived with 
the help of the excision properties of the index from the following input: 
First, the formula in Proposition \ref{thm:2.2} 
holds when $X$ is compact, see eg\
\cite{W1}. 
Second, take $X = \R \times T^{3}$ and $\omega $ to be a constant, 
non-zero self-dual $2$--form. 
Then, take the solution $c$ to be the pull-back via 
the projection to $T^{3}$ of a solution in Lemma \ref{thm:3.1}'s 
space $\mm_{P}$. Here, 
the form $\omega_{0}$ in \eqref{eq:3.1} is the pull-back to $\{0\}\times 
T^{3}$ of $\omega $. In this case, the kernel and the cokernel of $\dd_{c}$ 
are both trivial. (The operator has constant coefficients, so is 
straightforward to analyze.)
\end{proof}

\begin{proof}[Proof of Proposition \ref{thm:2.3}] 
The argument here is, modulo some 
notational changes, almost identical to that which proves the analogous 
assertion in the case where $X$ is compact; see, for example the books 
\cite{M} or 
\cite{KKM}. 
The largest modifications to the compact case argument are needed to 
address the orientation assertion, and in this regard, the reader can refer 
to the proof of Corollary 9.2 in \cite{MST}. 
\end{proof}

Lemmas \ref{thm:3.1}--\ref{thm:3.5} 
also have immediate applications to the subject of $\mm$'s 
compactness. The particular applications here are summarized below in 
Proposition \ref{thm:3.6}. Here is the background for this proposition: The 
proposition introduces the manifold with boundary $X_{0}$ as described in 
the beginning of Section \ref{sec:2}, 
except that here, $X_{0}$ is assumed to have 
non-empty boundary. Put a Riemannian metric on $X_{0}$ which is a flat, 
product metric on some neighborhood of $\partial X_{0}$ and extend this 
metric in 
the usual way to obtain a metric on $X=X_{0}\cup([0,\infty)\times\partial 
X_{0})$. Also, fix a self-dual $2$--form $\omega $ on $X$ 
which is non-zero and 
covariantly constant on $[0, \infty)  \times  \partial X_{0}$. Having 
made these 
selections, choose an element ${s}\in\sss_{0}(X)$ and then 
introduce, as in Section \ref{sec:2}d, 
the set $\varsigma({s})\subset H^{2}(X_{0},\partial X_{0};\Z)$ 
which consists of those elements $z$ which map 
to $c({s})$ in $H^{2}(X_{0};\Z)$.

What follows might be called a partial compactness assertion.

\begin{prop}\label{thm:3.6}
The restrictions of the chosen metric and 
form $\omega $ to $\partial X_{0}$ determine a constant $\delta  
> 0$ with the following significance: For each $z \in\varsigma({s})$, 
construct the moduli space $\mm$ and then for 
each $r\ge1$, introduce the subspace $\mm(r)\subset\mm$ of 
orbits of $(A,\psi)$ for which
\begin{equation*}
\int_{[r,\infty)\times\partial X_{0}}|F_{A}|^{2}\le\delta.
\end{equation*}
Then, $\mm(r)$ is compact in all cases, and actually empty for 
all but a finite set of $z\in\varsigma({s})$.
\end{prop}

\begin{proof}[Proof of Proposition \ref{thm:3.6}]\sloppy
As Witten pointed out \cite{W1}, the key to 
compactness theorems for the Seiberg--Witten moduli spaces is the 
Bochner--Weitzenb\"ock formula in \eqref{eq:3.6}. 
Of course, if $D_{A}\psi=0$, then 
the left hand side of \eqref{eq:3.6} 
is zero; thus contracting both sides with $\psi$ 
using the hermitian metric on $S_{+}$ yields a differential inequality for 
the function $|\psi|^{2}$. Moreover, when the first point in \eqref{eq:2.1} 
also 
holds, then, as noted in \cite{KM}, the maximum principle applies to this 
differential inequality and provides a uniform upper bound for $|\psi|^{2}$ 
in terms of $|s|$, $|\omega|$ and an asymptotic bound for $|\psi|^{2}$ 
on the ends of $X$. Lemma \ref{thm:3.5} provides such an asymptotic bound for 
$|\psi|$, so 
\begin{equation}\label{eq:3.11}
|\psi|^{2}\le\zeta\sup_{X}(|s| + |\omega|)
\end{equation}
on the whole of $X$. Here, $\zeta $ depends only on the Riemannian metric. 
Note, by the way, that this constant $\zeta $ is independent of both the 
$\Spin^{\C}$ structure ${s}$ and $z \in \varsigma({s})$. 

\fussy
This bound on $|\psi|$ together with Lemma \ref{thm:3.5} 
provide a uniform upper 
bound on the $L^{2}$ norm of $P_{+}F_{A}$. This upper bound is also 
independent of both ${s}$ and $z \in \varsigma({s})$. This 
last bound on $P_{+}F_{A}$ provides an $L^{2}$ bound on the 
anti-self dual part, $P_{-}F_{A}$, of $A$'s curvature $2$--form via the 
string of identities
\begin{equation}\label{eq:3.12}
c({s})\bullet c({s}) = -(4\pi^{2})^{-1}\int_{X}F_{A}\wedge F_{A}= 
(4\pi^{2})^{-1}\int_{X}
(|P_{+}F_{A}|^{2} - |P_{-}F_{A}|^{2}).
\end{equation}
In particular, notice that the resulting upper bound for the $L^{2}$ norm of 
$F_{A}$ is independent of the class $z\in\varsigma({s})$. 

Standard elliptic regularity techniques can be invoked to bootstrap these 
upper bounds on $F_{A}$ and $|\psi|$ into uniform and $z\in\varsigma 
({s})$ independent upper bounds for all $C^{k}$ norms for a suitable 
point on the gauge orbit of $(A,\psi)$. (See \cite{M} 
or \cite{KKM} to see how this 
is done.) In particular the latter imply that any sequence of gauge orbits 
in $\cup_{z\in\varsigma(s)}\mm(r)$ is defined by a corresponding 
sequence $\{(A_{k},\psi_{k})\}$ of solutions to \eqref{eq:3.1} which converges 
in the $C^{\infty}$ topology on compact subsets of $X$. Meanwhile, the uniform 
bounds that are provided by Lemma \ref{thm:3.5} 
imply that a sequence of gauge orbits 
in $\cup_{z\in\varsigma(s)}\mm(r)$ is actually defined by a sequence 
$\{(A_{k},\psi_{k})\}$ of solutions to \eqref{eq:2.1} which converges in the 
strong $C^{\infty}$ topology on the whole of $X$. This last fact implies 
convergence in each $\mm(r)$ 
and it implies that there can be only finitely many 
$z\in\varsigma({s})$ for which the corresponding $\mm(r)$ is not 
empty.
\end{proof}

\sh{d) The family version of Proposition \ref{thm:2.4}}

Though Proposition \ref{thm:3.6} 
asserts that each $\mm(r)$ is compact, there is no reason 
for the whole of $\mm$ to be compact. In fact, most probably, the statement in 
Proposition \ref{thm:2.4} 
is about as strong as can generally be made. As remarked in 
Section \ref{sec:2}d, 
the compactness asserted in Proposition \ref{thm:2.4} is strong enough to 
provide a definition of $\sw({s}, z)$ given a reasonable choice of 
$\omega $ and, in the $d > 0$ cases, a suitably generic choice of the set 
$\underline\Lambda$. 

However, the compactness asserted by Proposition \ref{thm:2.4} 
is not strong enough 
for use in the proof of Theorem \ref{thm:2.5}. 
Indeed, a comparison between the values 
of $\sw(\cdot)$ as defined by different, but still allowable choices for the 
triple of metric, $\omega$ and $\underline\Lambda $ involves an 
interpolating path in the space of 
such triples, and thus a corresponding 
$1$--parameter family of moduli spaces.     
Meanwhile, Proposition \ref{thm:2.4} says nothing 
about compactness for families of moduli spaces. This weakness in 
Proposition \ref{thm:2.4} 
is addressed below with the statement of Proposition \ref{thm:2.4}'s 
family version for use in the proof of Theorem \ref{thm:2.5}.

To set the stage for the family version of Proposition \ref{thm:2.4}, 
fix $s\in\sss_{0}(X_{0})$ and consider two sets of triples, $\Gamma_{0}  
\equiv(g_{0},\omega_{0}, \underline\Lambda_{0})$ and 
$\Gamma_{1}\equiv(g_{1},\omega_{1}, \underline{\Lambda}_{1})$, 
for use in defining $\sw({s}, z)$ for $z\in\varsigma({s})$. 
Thus, $g_{0}$ and $g_{1}$ are metrics on $X$ which restrict to 
flat, product metrics on $[0, \infty)\times\partial X_{0}\subset X$. 
Meanwhile, $\omega_{0}$ and $\omega_{1}$ are self-dual forms on $X$ for 
the respective metrics $g_{0}$ and $g_{1}$ which restrict to each component 
of $[0, \infty)  \times  \partial X_{0}$ 
as non-zero, covariantly constant $2$--forms. 
Finally, $\underline\Lambda_{0}$ 
and $\underline\Lambda_{1}$ 
are 
two sets of points
and the associated data needed to 
define $\sw$ when the number $d$ in Proposition \ref{thm:2.2} is positive. 

The triple $\Gamma_{0}$ has its associated moduli space $\mm({s}, z)$. 
There is, of course, an analogous moduli space that is defined by $\Gamma 
_{1}$ and these two spaces will be distinguished as $\mm^{0}$ and $\mm^{1}$, 
respectively. With this notation understood, require now of $\omega_{0}$ 
and $\omega_{1}$ that their moduli spaces $\mm^{0}$ and $\mm^{1}$ consist 
only of smooth points. Meanwhile, require of $\underline\Lambda_{0}$ 
and $\underline\Lambda_{1}$ that the conditions which define their 
respective versions of the set $\mm^{\underline\Lambda}$ in \eqref{eq:2.5} 
cut the latter out of 
the corresponding $\mm$ in a transversal fashion. These will be denoted 
respectively by $\mm^{0\underline\Lambda}$ and $\mm^{1\underline\Lambda}$.

Here is the final key assumption on $\Gamma_{0}$ and $\Gamma_{1}$. 
Assume that there exists a continuous, $1$--parameter path $\{\varpi_{t}  
\in H^{2}(X_{0};\mathbb{R}): t\in[0, 1]\}$ with each $\varpi_{t}$ 
having non-zero pull-back in the cohomology of 
each component of $\partial X_{0}$, 
and with $\varpi_{0}$ and $\varpi_{1}$ taming $\omega_{0}$ and 
$\omega_{1}$, respectively. 

With these assumptions in hand, consider now the following:

\begin{prop}\label{thm:3.7}
Under the assumptions just made, there 
exists a continuous, interpolating family, $\{\Gamma_{t}: t \in [0, 
1]\}$, of data triples for which the corresponding family of moduli 
spaces $\ww \equiv \cup_{t\in[0,1]}\mm^{t}$ has the 
following structure:
\begin{itemize}
\item  $\ww$ is a smooth, oriented,  $d+1$ dimensional 
manifold with boundary for which the tautological map to $[0, 1]$ is 
smooth and a product over a neighborhood of $\{0, 1\}$. Moreover, the 
induced boundary orientation on $\mm^{1}$ agrees with its orientation 
from Proposition \ref{thm:2.3}, while that on $\mm^{0}$ disagrees.
\item  Let $\ww^{\underline\Lambda} \equiv \cup_{t \in [0,1]}
\mm^{t\underline\Lambda} \subset\ww$. Then $\ww^{\underline\Lambda}$ has the 
structure of a finite, disjoint set of 
embedded, oriented intervals each mapping
to $[0, 1]$ as a product 
over a neighborhood of $\{0, 1\}$. 
Moreover, the induced boundary orientation on $\mm^{1\underline\Lambda}$ 
agrees with its orientation from the definition of $\sw({s},z)$, 
while that on $\mm^{0\underline\Lambda}$ disagrees with its
$\sw({s}, z)$ orientation. 
\end{itemize}
\end{prop}
The proof of this proposition is discussed in Section \ref{sec:6}b 
so accept its 
assertions for the time being.

With Proposition \ref{thm:3.7} in hand, consider the following:

\begin{proof}[Proof of Theorem \ref{thm:2.5}] 
The invariance assertions in the theorem are 
a standard consequence of the third point in Proposition \ref{thm:3.7}. 
Indeed, the 
intervals in $\ww^{\underline\Lambda}$ 
pair each point in $\mm^{0\underline\Lambda}$ either with 
another point in this space, but one with the opposite sign for the $\sw$ 
count, or else with a point in $\mm^{1\underline\Lambda}$ 
which has the same sign for 
the $\sw$ count. 
Meanwhile, each point in the latter space which is not paired 
by a component of $\ww^{\underline\Lambda}$ 
to one in $\mm^{0\underline\Lambda} $ is paired by a 
component of $\ww^{\underline\Lambda}$ 
with another such point, but one with the 
opposite $\sw$ count sign.
\end{proof}

\sh{e) Gluing moduli spaces}

In this subsection, $X_{0}$ denotes a compact, connected, oriented 
$4$--manifold 
with boundary consisting of a disjoint union of $3$--dimensional 
tori. Here, the boundary can be empty, but if so, require that $b^{2+}\ge1$. 
Let $M \subset X_{0}$ be an embedded, $3$--dimensional torus, and if $M$ is 
separating, write $X_{0} = X_{-} \cup X_{+}$, where $X_{+} \cap X_{-} = M$. 
Otherwise, let $X_{1} \subset X$ denote the complement of 
an open, tubular neighborhood of $M$. In the separating case, introduce the 
non-compact manifolds $\underline{X}_{\pm} \equiv X_{\pm} \cup 
([0, \infty) \times \partial X_{\pm})$, and in the non-separating case, 
$\underline{X}_{1} \equiv X_{1} \cup ([0, \infty) \times \partial X_{1})$. 

With the introduction of 
$\underline{X}_{\pm}$ and $\underline{X}_{1}$,
the remainder of this subsection 
describes how moduli spaces on $\underline{X}_{\pm}$, 
in the separating case, or on $\underline{X}_{1}$ otherwise, 
can be glued together over the ends that contain $M$ to produce portions of 
the moduli space for $X$.

To begin the presentation, choose a flat metric on $M$, and then choose a 
metric on $X_{0}$ which restricts as a flat product metric on an interval 
neighborhood of $M$. This is to say that the metric should allow for an 
isometric embedding of $(-\varepsilon,\varepsilon)\times M$ into $X$ 
which sends $\{0\} \times M$ to $M$. The chosen metric for $X_{0}$ should 
also restrict as a product, flat metric 
on a neighborhood of $\partial X_{0}$. 

With this metric fixed, select a self-dual $2$--form $\omega $ which is 
non-zero and constant on a neighborhood of $M$ and likewise on a neighborhood 
of each component of $\partial X_{0}$. The 
metric on $X_{0}$ and the form $\omega $ 
induce, in a presumably obvious way, a pair of metric and self-dual $2$--form 
on $\underline{X}$ and also on $\underline{X}_{\pm}$ or $\underline 
{X}_{1}$, as the case may be. Here, the metric on these spaces is flat and 
a product on the ends, and the self-dual $2$--form, 
still called $\omega $, is 
non-zero and constant on each end component. Moreover, in the separating 
case, $\underline{X}_{\pm} $ have one special end, that which contains 
$M\subset\partial X_{\pm}$. In $X_{-}$, 
this end has an orientation preserving 
isometry with $[0, \infty)\times M$, while in $X_{+}$, the orientation 
preserving isometry sends the end to $(-\infty, 0]  \times  M$. In the 
non-separating case, $\underline{X}_{1}$ has two special ends, one with an 
orientation preserving isometry to $[0, \infty) \times M$ and the other to 
$(-\infty, 0] \times  M$. In all of these cases, the metric and the 
form $\omega $ 
restrict to these special ends as the constant extension of the given metric 
and form on $(-\varepsilon,\varepsilon) \times M \subset X$.

With regard to the choice of $\omega $, Proposition \ref{thm:2.3} asserts the 
following: Fix neighborhoods of $M$ and 
$\partial X$  whose closure is not the whole of $X$ 
and there is a Baire subset of choices for $\omega $ which have the given 
restriction near $M$ 
and $\partial X$  and are such that all moduli spaces of solutions 
to \eqref{eq:2.1} on $X$, $\underline{X}_{+}$ and $\underline{X}_{-}$, or on 
$\underline{X}_{1}$, are smooth manifolds for which the operator $\dd_{c}$ 
has trivial cokernel at all points. In particular, when choosing the form 
$\omega $, be sure to take one for which this last conclusion holds.

The chosen metric on $X$, call it $g$, will now be used to construct a 
$1$--parameter 
family, $\{g_{R}\}_{R\ge0}$, of metrics on $X$. Here, $g_{0} = g$, 
while $X$ with the metric $g_{R}$ admits an isometric embedding of the 
cylinder $[-R, R] \times M$ whose complement with the metric $g_{R}$ is 
isometric to $X - M$ with the metric $g$. Alternately, the metric $g_{R}$ 
makes $X$ isometric to
\begin{equation}\label{eq:3.13}
(X - M) \cup ([-R, R] \times M),
\end{equation}
where $X - M$ has the metric $g$ and $[-R, R] \times M$ has the product, 
flat metric. With $g_{R}$ understood, introduce $X^{R}$ to denote $X$ as a 
Riemannian manifold with the metric $g_{R}$.

With regard to \eqref{eq:3.13}, 
note that in the respective cases where $M$ does and 
does not separate $X$, the Riemannian manifold $X^{R} - M$ admits a 
canonical, orientation preserving isometry
\begin{equation}\begin{split}\label{eq:3.14}
\Theta\colon X^R-M \to&\left(\underline{X}_{-}-([R,\infty)\times M)\right)\\
&\cup\left(\underline{X}_{+}-((-\infty,-R]\times M)\right)\subset
\underline{X}_{-}\cup\underline{X}_{+}\text{ or}\\
\Theta\colon
X^R-M \to&\underline{X}_1-\left(([R,\infty)\times M)
\cup((-\infty,-R]\times M)\right)\subset\underline{X}_{1},
\end{split}\end{equation}
respectively.

There is one last remark to make here about the metric $g_{R}$, which is 
that the form $\omega $ can be viewed as living on $X^{R}$ in as much as its 
restriction to $X - M \subset X$ defines it on the isometric $X - M  
\subset X^{R}$ and then there is an evident extension as a self-dual 
$2$--form on the whole of $X^{R}$ which is non-zero and constant on the 
cylinder $[-R, R] \times M$.

With the geometric preliminaries complete, now choose $({s}, z) \in 
\sss_{0M}(X_{0},\partial X_{0})$. In the case where $M$ 
is separating, each pair 
$(({s}_{-},z_{-}),({s}_{+},z_{+})) \in \wp^{-1}({s},z) \subset 
\sss_{0}(X_{-},\partial X_{-}) \times \sss_{0}(X_{+},\partial X_{+})$ 
determines moduli spaces $\mm_{-}$ and $\mm_{+}$ on 
$\underline{X}_{-}$ and $\underline{X}_{+}$, respectively. In the 
non-separating case, each $({s}_{1},z_{1}) \in \wp^{-1}({s},z)$ 
determines the moduli space $\mm_{1}$ on $\underline{X}_{1}$. 
With regard to these spaces, remember that the form $\omega $ has 
been chosen so that these moduli spaces are smooth manifolds with 
$\cokernel(\dd_{c}) = \{0\}$ at all points. Note for use below that the 
specification of $r \ge 1$ 
defines the subsets $\mm_{\pm}(r) \subset \mm_{\pm}$ 
and $\mm_{1}(r) \subset \mm_{1}$ of pairs $(A,\psi)$ which obey 
the curvature bound in Proposition \ref{thm:3.6}.

Meanwhile, the choice of $R \ge 0$ 
determines the Riemannian manifold $X^{R}$, 
and then, with the help of $\omega $, the pair $({s}, z)$ determine a 
corresponding moduli space, $\mm^{R}$. 

With all of this understood, what follows in Proposition \ref{thm:3.8} is a 
description of the fundamental fact that parts of $\mm^{R}$ 
can be constructed 
from $\mm_{-}(r) \times \mm_{+}(r)$ or from $\mm_{1}(r)$ as the case may be. 
Here, a lower bound for $R$ is determined by $r$. 

\begin{prop}\label{thm:3.8} 
Each pair $r'\ge0$ and 
$\delta > 0$ determines a lower bound for a choice of $r$ 
and then a choice of $r\gg r'$ which is greater than this 
lower bound determines a lower bound for a choice of $R\gg r$. Choose 
$r'\ge 0$, then $r\gg r'$ consistent with the lower bound, and 
finally $R\gg r$ consistent with its lower bound. Then the 
following are true:

\begin{itemize}
\item  In the case where $M$ separates, there exists an 
embedding
\begin{equation*}
\Phi_{r}\colon\bigcup_{((s_{-},z_{-}),(s_{+},z_{+}))\in\wp^{-1}((s,z))}
\mm_{-}(r)\times\mm_{+}(r)\to\mm^{R},
\end{equation*}
which maps the interior of its domain onto an open set of smooth 
points that contains the subspace of $(A,\psi)$ with
\begin{equation}\label{eq:3.16}
\int_{[-R+r',R-r')\times M} 
|F_{A}|^{2} \le \delta/2.
\end{equation}
In addition, if the lines $L_{-}$, $L_{+}$, and 
$L_{0}$ are oriented as in Theorem \ref{thm:2.7} and if the induced 
orientations on the respective moduli spaces are used, then $\Phi 
_{r}$ is orientation preserving. Moreover, if $(c_{-},c_{+})$ 
is in the domain of $\Phi_{r}$, there are 
points $(A_{\pm},\psi_{\pm}) \in c_{\pm}$ and $(A,\psi) \in 
\Phi_{r}(c_{-},c_{+})$ such that
\begin{equation*}
\sum_{0\le k\le2}\left(|\nabla^{k}(A - \Theta^*(A_{-},A_{+})| + |\nabla
_{A}^{k}(\psi - \Theta^*(\psi_{-},\psi_{+})|\right) \le e^{-R/\zeta},
\end{equation*}
where $\Theta$ is the map in the first line of \eqref{eq:3.14} 
and $\zeta\ge1$ depends only on the restriction 
to $M$ of the form $\omega $ and the metric $g$.

\item  In the case where $M$ does not separate, there 
exists an embedding
\begin{equation*}
\Phi_{r}\colon\bigcup_{((s_l,z_l))\in\wp^{-1}((s,z))}   
\mm_{1}(r)\to\mm^{R},
\end{equation*}
which maps the interior of its domain onto an open set of smooth 
points that contains the subspace of $(A,\psi)$ which obey 
\eqref{eq:3.16}. In addition, if the lines $L_{1}$ and 
$L_{0}$ are oriented as in Theorem \ref{thm:2.7} and if the induced 
orientations on the respective moduli spaces are used, then $\Phi 
_{r}$ is orientation preserving. Moreover, if $c_{1}$ is 
in the domain of $\Phi_{r}$, there are points $(A_{1},\psi_{1}) \in 
c_{1}$ and $(A,\psi)\in\Phi_{r}(c_{1})$ such that
\begin{equation*}
\sum_{0\le k\le2}\left(|\nabla^{k}(A-\Theta^*A_{1})| + |\nabla_{A}^{k}(\psi - 
\Theta^*\psi_{1})|\right) \le e^{-R/\zeta},
\end{equation*}
where $\Theta$ is the map in the second line of \eqref{eq:3.14} 
and $\zeta\ge 1$ depends only on the restrictions to 
$M$ of $\omega$ and $g$.
\end{itemize}
\end{prop}

As the gluing result in Proposition \ref{thm:3.8} 
is now standard fair in gauge 
theories, the proof will be omitted. However, note that a proof can be 
produced via a straightforward application of the techniques introduced in 
\cite{DK} for gluing self-dual moduli spaces in $SU(2)$ gauge theory. In this 
regard, the choice of $\omega $ to make either each pair $\mm_{\pm}$ or each 
$\mm_{1}$ smooth insures the absence of obstruction bundles. As a last remark 
on this subject, note that reference to Section 9.1 of \cite{MST} can also be 
made to deal with the assertions in the proposition that concern 
orientations.

\sh{f) Implications from gluing moduli spaces}

This subsection discusses the application of Proposition \ref{thm:3.8} 
to the proof of 
Theorem \ref{thm:2.7}. 
For this purpose, suppose that $X_{0}$, $M$, $X_{\pm} $ and 
$X_{1}$ are as described at the outset of the preceding subsection. 

Coupled with Proposition \ref{thm:2.4}, 
the gluing result in Proposition \ref{thm:3.8} is almost 
enough to prove Theorem \ref{thm:2.7}. 
Indeed, missing still is a guarantee that each 
point in $(\mm^{R})^{\underline\Lambda}$ satisfies \eqref{eq:3.16} 
for some fixed $r'$ and all $R$ 
sufficiently large. The next proposition provides such a guarantee:

\begin{prop}\label{thm:3.9}
Continuing the discussion and notation 
from Proposition \ref{thm:3.8} and the preceding subsection, introduce $d$ as 
defined in Proposition \ref{thm:2.2}. If $d > 0$, choose the data set 
$\underline\Lambda$ so that the set of base points $\Lambda$ is disjoint from
$M$ and use \eqref{eq:3.13} to define 
the analogous data sets for each $X^{R}$. Choose a class $\varpi\in 
H^{2}(X_{0};\mathbb{R})$ whose pull-back in the cohomology of 
$M$ and in that of every component of $\partial X_{0}$ is non-zero. 
Next, choose a self-dual $2$--form $\omega'$ on $X_{0}$ 
which is tamed by $\varpi $ and whose restriction to the fiducial 
tubular neighborhood of $M$ and to that of each component of 
$\partial X_{0}$ is constant and non-zero. Then there is a Baire set of 
choices for the self-dual forms $\omega $ which are tamed by 
$\varpi $, which agree with $\omega'$ on the complement 
of a fixed, tubular neighborhood of $M\cup\partial X_{0}$, and which 
have the following additional property: If $r$ is 
sufficiently large, and $R$ also, subject to its lower bound 
constraints, then the set $\Lambda $ can be chosen to insure that 
$(\mm^{R})^{\underline\Lambda} $ lies in the image of Proposition 
\ref{thm:3.8}'s map 
$\Phi_{r}$. In particular, there exists $r' \ge 1$ such 
that when $R$ is sufficiently large, then each point in 
$(\mm^{R})^{\underline\Lambda} $ obeys \eqref{eq:3.16}.
\end{prop}
This proposition is proved in Section \ref{sec:6}c. 

The preceding two propositions play the key roles in the following:

\begin{proof}[Proof of Theorem \ref{thm:2.7}] 
Except for two assertions, the theorem 
follows directly from Propositions \ref{thm:3.7} 
and \ref{thm:3.9} via arguments which are 
standard fair in gauge theories. (Arguments of this sort were first given by 
Donaldson in the context of $SU(2)$ gauge theories, see eg\ \cite{DK}.) The 
assertions which do not follow 
immediately from Propositions \ref{thm:3.7}
and \ref{thm:3.8} 
concern the vanishing of all but finitely many of the numbers 
$\sw({s}_{-},z_{-})$ and $\sw({s}_{+},z_{+})$ or 
$\sw({s}_{1},z_{1})$ as the case may be. However, these last 
assertions follow from Proposition \ref{thm:2.4}.
\end{proof}

\section{Energy and compactness}\label{sec:4}

The estimates provided here make the first steps towards the proof of 
Proposition \ref{thm:2.4}. 
Although they are not strong enough to give Proposition \ref{thm:2.4} 
and its brethren in their entirety, they do yield the previously mentioned 
result that $\mm$ is compact if the pull-back to each boundary component of 
$\varpi $ is not a multiple of an integral cohomology class. 

In the subsequent discussions of this section, $X_{0}$ is as described at 
the beginning of Section \ref{sec:2}, 
a compact, connected, oriented $4$--manifold with 
boundary such that each boundary component is a $3$--torus. As before, $X_{0}$ 
is endowed with a metric which is a product flat metric on a tubular 
neighborhood of the boundary. Also, $X_{0}$ is endowed with a class 
$\varpi\in H^{2}(X_{0};\mathbb{R})$ whose pull-back to the cohomology of each 
component of $\partial X_{0}$ is non-zero. As previously, let $X$ 
denote $X_{0}\cup([0,\infty)\times\partial X_{0})$ 
with the induced Riemannian metric, and 
let $\omega $ be a self-dual $2$--form on $X$ 
which is non-zero and covariantly 
constant on each component of $[0, \infty)\times\partial X_{0}$, 
and which is tamed 
by $\varpi $.

\sh{a) The first energy bound}

The first of the important energy bounds is presented below as Proposition 
\ref{thm:4.1}. Its proof then occupies the remainder of this subsection.

\begin{prop}\label{thm:4.1}
There exist universal, positive constants 
$\{\kappa_{j}\}_{j=1,\ldots,4}$, and a positive constant, 
$\zeta $, which depends on the Riemannian metric, $\varpi 
$ and $\omega $; and these constants have the following 
significance: As usual, let $\mm$ denote the moduli space for 
a given $({s}, z) \in\sss_{0}(X_{0},\partial X_{0})$, and let 
$(A,\psi) \in \mm$. Then,  
\begin{equation*}\begin{split}
&\int_{X}(|\nabla_{A}\psi|^{2} + 
|P_{+}F_{A}|^{2}) \le \zeta + 
\kappa_{1}z\bullet\varpi - \kappa 
_{2}c({s})\bullet c({s})\text{ and}\\            
&\int_{X}|F_{A}|^{2} \le \zeta + 2\kappa_{1} z 
\bullet \varpi - \kappa_{3} c({s}) \bullet c({s}).
\end{split}\end{equation*}
Moreover, if $\omega$ is a closed form, then
\[\int_{X}(|\nabla_{A}\psi|^{2} + 
|P_{+}F_{A}|^{2}) \le \zeta + \kappa_{4}z\bullet\varpi.\]
Finally, in the case where $X = \R \times T^{3}$, these 
results hold with $\zeta = 0$. 
\end{prop}

\begin{proof}[Proof of Proposition \ref{thm:4.1}]
To obtain the proposition's first assertion, 
start with the Bochner--Weitzenb\"ock formula in \eqref{eq:3.6}, 
take the inner product 
of both sides of the latter with $\psi $, and then integrate the result over 
$X$. Integration by parts (which \eqref{eq:3.3} justifies) and input from the 
Seiberg--Witten equations can then be used to derive the equality
\begin{equation*}
0 = \int_{X}\left(|\nabla_{A}\psi|^{2}+\frac{{s}}{{4}} 
|\psi|^{2}+\frac{{1}}{{2}}|P_{+}F_{A}|^{2} - 
\frac{{i}}{{2}} F_{A}\wedge\omega\right).
\end{equation*}
Next, fix a closed $2$--form $\mu $ which represents the class $\varpi $ and 
equals $\omega $ on $[0, \infty)  \times  \partial X_{0}$. 
Because both $s$ and $\omega-\mu$ 
are supported on $\partial X_{0}$, this last equation with \eqref{eq:3.11} 
implies 
the inequality
\begin{equation}\label{eq:4.2}
\int_{X}(|\nabla_{A}\psi|^{2} + |P_{+}F_{A}|^{2}) \le \zeta' 
+ \kappa_{0}\int_{X} iF_{A}\wedge\mu + 
\frac{{1}}{{16}}\int_{X_0}|F_{A}|^{2}.
\end{equation}
Here, $\zeta'$ depends only on the metric, $\omega $ and $\mu $, while 
$\kappa_{0}$ is a positive, universal constant. With regard to \eqref{eq:4.2}, 
note that when $\mu $ is self-dual, the integrand of the last term in 
\eqref{eq:4.2} 
can be replaced by $|P_{+}F_{A}|^{2}$. 

The first integral on the right side of \eqref{eq:4.2} 
is equal to $2\pi\cdot z\bullet\varpi$. 
Meanwhile, \eqref{eq:3.12} 
relates the integral in \eqref{eq:4.2} of $|F_{A}|^{2}$ to that 
of $|P_{+}F_{A}|^{2}$. With these relations understood, the inequality 
in the first point of the proposition follows directly from \eqref{eq:4.2}. 
The 
second point's inequality follows from the first and \eqref{eq:3.12}. 
Meanwhile, the 
proposition's third point follows directly from \eqref{eq:4.2} given that the 
integrand in the last term on the right is replaced by 
$|P_{+}F_{A}|^{2}$. 
\end{proof}

\textbf{\boldmath b) Uniform asymptotics of $(A,\psi)$}

The bounds on $P_{+}F_{A}$ and $\nabla_{A}\psi$ that Proposition \ref{thm:4.1} 
provides are crucial inputs to a key generalization of Lemma \ref{thm:3.5}. 
The latter 
is stated next as Proposition \ref{thm:4.2}. 
The proof of this proposition then 
occupies the remainder of this subsection.

\begin{prop}\label{thm:4.2}
There exist positive constants, $\zeta_{0}$ and $\zeta_{1}$, 
which depend on the Riemannian 
metric, $\varpi $ and $\omega $; and which have the 
following significance: Given $({s}, z) \in  
\sss_{0}(X_{0},\partial X_{0})$, set
\begin{equation}\label{eq:4.3}
f \equiv \kappa_{1} z\bullet\varpi - \kappa_{2} 
c({s})\bullet c({s}), 
\end{equation}
where $\kappa_{1}$ and $\kappa_{2}$ are the 
constants that are introduced in Proposition \ref{thm:4.1}. Let $\mm$ 
denote the 
moduli space for $({s}, z)$ and let $(A,\psi)\in\mm$. 
Then:
\begin{itemize}
\item  $|F_{A}| < \zeta_{0} + \zeta_{1}f$ 
everywhere on $X$. 
\item  Let $[0, \infty)  \times  T^{3}$ be an end of $X$. 
There exists some number $N \le \zeta_{0} + \zeta_{1}f$ 
of points $\{t_{i}\} \in [0, \infty)$ (not necessarily 
distinct) such that 
\begin{equation}\label{eq:4.4}
|F_{A}|(t,\cdot) \le \zeta_{0}\left(e^{-t/\zeta_{0}}+ 
\sum_{i}                                         
e^{-(t-t_{i})/\zeta_{0}}\right) 
\end{equation}
at all points $(t,\cdot)\in[0,\infty)\times T^{3}$.
\end{itemize}
\end{prop}

\sloppy
\begin{proof}[Proof of Proposition \ref{thm:4.2}] 
The first bound on $|F_{A}|$ follows from 
the $L^{2}$ bounds in Proposition \ref{thm:4.1} 
using standard elliptic regularity 
techniques. To obtain the bound in \eqref{eq:4.4}, let $\{t_{j}\} 
\subset [2, \infty) $
denote the distinct integer points where the bound in Lemma \ref{thm:3.3} 
is violated 
for the cylinder $Y = [t_{j} - 2, t_{j} + 2] \times T^{3}$. 
Proposition \ref{thm:4.1} 
guarantees that there are no more than $N = (\zeta+f)/\delta$ 
such points where $\zeta $ is given in Proposition \ref{thm:4.1}. The 
estimate in \eqref{eq:4.4} then follows from Lemma \ref{thm:3.4}.
\end{proof}\fussy

\sh{c) Refinements for the cylinder}

The following proposition describes a useful refinement of the inequality in 
\eqref{eq:4.4} which holds when $X$ is the cylinder $\R \times T^{3}$. 

\begin{prop}\label{thm:4.3}  
Let $R \ge 2$ and let $X = [-R, R] 
\times T^{3}$. Let $\omega $ be a non-zero, covariantly 
constant, self-dual $2$--form on $X$. There exist constants $\delta 
> 0$ {and} $\zeta_{0} \ge 1$ which are independent of $R$ 
and which have the following significance: Let $(A, \psi 
)$ be a solution to \eqref{eq:2.1} on $X$ which obeys $\int_{Y}|F_{A}|^{2} 
< \delta$ when $Y = [-R, -R + 4] \times 
T^{3}$ and $Y = [R - 4, R] \times  T^{3}$. {Then:}

\begin{itemize}
\item  $|F_{A}| < \zeta_{0}$ {everywhere}. 

\item  {There exists some number} $N$ {of distinct, integer 
valued points} $\{t_{i}\} \in [-R, R]$ {such that} 
\begin{equation}\label{eq:4.5}
|F_{A}| \le \zeta_{0} \left(\delta e^{-(R - |t|)/\zeta_{0}} + 
\sum_{i}e^{-(t-t_{i})/\zeta_{0}}\right) 
\end{equation}
at all points $(t,\cdot)$ where $t \in [-R+2, R-2]$.
\end{itemize}
\end{prop}

\begin{proof}[Proof of Proposition \ref{thm:4.3}]  
First, note that $|\psi|$ is bounded by some number $B$ on 
$X' = [-R+1, R-1] \times  T^{3}$ which depends only on the metric and 
$\omega $. To find the bound $B$, first use the Seiberg--Witten equations to 
rewrite $P_{+}F_{A}$ in \eqref{eq:3.7} in terms of $\psi $ and so derive the 
following differential inequality for $|\psi|^{2}$:
\begin{equation}\label{eq:4.6}
2^{-1} d^*d|\psi|^{2} + m|\psi|^{2}(|\psi|^{2} - 
|\omega|) \le 0.
\end{equation}
Here, $m > 0$ is a universal constant. The maximum principle applies to 
\eqref{eq:4.6} 
and bounds $|\psi|^{2}$ in terms of $|\omega|$ and its size near the 
ends of $X'$. Meanwhile, Lemma \ref{thm:3.3} 
bounds the size of $|\psi|^{2}$ near the 
ends of $X'$ in terms of $|\omega|$. 

For the next step in the proof, take $t \in [-R+2, R-2] \times 
T^{3}$, set $T =\linebreak
{[t-1, t+1]} \times \R$ and then multiply both sides of 
\eqref{eq:3.7} 
by a standard, smooth function on $[-R, R]$ which is $1$ on $T$ and $0$ on 
$[-R, t - 2]$ and on $[t + 2, R]$. Integrate the result over $X'$, and then 
integrate by parts to remove the operator $d^*d$ from $|\psi|^{2}$. With 
the bound on $|\psi|$ by $B$, a simple manipulation gives the inequality 
\begin{equation}\label{eq:4.7}
\int_{T}(|\nabla_{A}\psi|^{2} + |P_{+}F_{A}|^{2}) \le \zeta' 
(1 + B^{2}); 
\end{equation}
here, $\zeta'$ is another constant which depends only on $|\omega|$.

The next task is to control the $L^{2}$ norm of $|P_{-}F_{A}|$. Here, the 
vanishing of $dF_{A}$ is employed to conclude that $d(P_{-}F_{A}) = - 
d(P_{+}F_{A})$. Now, note that by differentiating \eqref{eq:2.1}, the form 
$d(P_{+}F_{A})$ can be expressed in terms of $\psi $ and $\nabla 
_{A}\psi $ and as a result, its norm is bounded by a uniform multiple of 
$\zeta' B |\nabla_{A}\psi|$, where $\zeta'$ is a universal 
constant. Thus, \eqref{eq:4.7} 
provides an $L^{2}$ bound on $d(P_{-}F_{A})$. With 
this bound in hand, it now proves useful to rewrite the equation 
$d(P_{-}F_{A}) = -d(P_{+}F_{A}) \equiv \sigma $ by separating 
out $t$--derivatives 
from derivatives along the tori $T^{3}$. There result two 
equations,
\begin{equation}\begin{split}\label{eq:4.8}
&\partial_{t}f-{*}\partial f = \sigma_{\bot}\text{ and}\\
&\partial{*}f = \sigma_{0}.
\end{split}
\end{equation}

Here, $\partial$ 
is the exterior derivative along the torus, and ${*}$ denotes the 
Hodge star along the torus. Also, $f$ is the $t$--dependent $1$--form 
on $T^{3}$ 
which is obtained from $P_{-}F_{A}$ by contracting with the unit vector in 
the $t$--direction. 

To analyze the preceding equations, write $f = g + \partial u$, 
where $\partial{*}g = 0$ and 
where $u$ is a time dependent function on $T^{3}$ obeying 
$\partial{*}\partial u = \sigma_{0}$. 
Letting $\|\cdot\|_{t}$ denote the $L^{2}$ norm over $\{t\} 
\times T^{3}$, it follows by standard arguments that there is a solution 
$u$ which obeys
\begin{equation*}
\|\partial u\|_{t} \le \zeta \cdot \|\sigma_{0}\|_{t} 
\end{equation*}
for all $t \in [-R+2, R-2]$. 

Next, consider the projection of the top line in \eqref{eq:4.8} 
onto the kernel of $\partial 
* $. The result is an equation for $g$ which reads
\begin{equation}\label{eq:4.10}
\partial_{t}g-{*}\partial g = \sigma',
\end{equation}
where $\sigma'$ is the $L^{2}$--orthogonal projection of $\sigma_{\bot}$ 
onto the kernel of $\partial{*}$. (Thus, $\|\sigma'\|_{t} \le \|\sigma\|_{t}$.)
To analyze \eqref{eq:4.10}, consider the projections $g_{+}$, $g_{-}$ and 
$g_{0}$ onto the respective positive, negative, and zero eigenspaces of ${*} 
\partial$ acting on the kernel of $\partial{*}$. 
Likewise, introduce the analogous 
projections of $\sigma'$, namely $\sigma'_{+,-,0}$. Then 
\eqref{eq:4.10} implies
\begin{equation}\begin{split}\label{eq:4.11}
&(\partial_{t} - \lambda)\cdot \|g_{+}\|_{t} \ge -\|\sigma'_{+}\|_{t};\\
&(\partial_{t} + \lambda)\cdot \| g_{-} \|_{t} \le \| \sigma'_{-} \|_{t};\\
&\partial_{t} g_{0} = \sigma'_{0}.
\end{split}\end{equation}
Here, $\lambda$ is the smallest non-zero absolute value of an eigenvalue of 
${*}\partial$ on $\kernel(\partial{*})$. 

The first and second lines in \eqref{eq:4.11} can be integrated to yield
\begin{equation*}\begin{split}
&\|g_{+}\|_{t} \le e^{-\lambda \cdot (R-1-t)} \| g_{+} 
\|_{R-1} + \sup_{s\in[-R+1,R-1]}\int_{[s-1,s+1]}\|\sigma'_{+}\|_{s}^{2} ds
\text{ and}\\
&\|g_{-}\|_{t} \le e^{\lambda \cdot (t+1-R)} \| g_{-} 
\|_{-R+1} + \sup_{ s \in [-R+1,R-1]} \int_{[s-1,s+1]} \| 
\sigma'_{-} \|_{s}^{2} ds.
\end{split}
\end{equation*}
The final line can be integrated to find that
\begin{equation}\label{eq:4.13}
\| g_{0} \|_{t} \le \sup_{s \in [-R+1,R-1]} \|P_{+}F_{A} 
\|_{s}.
\end{equation}
The explicit formula for $\sigma'_{0}$ must be used to derive 
\eqref{eq:4.13}. 
For this purpose, introduce the $t$--dependent $1$--form $h$ on $T^{3}$ 
by writing 
$P_{+}F_{A} = dt \wedge h +{*}h$. Then, up to a sign, $\sigma'_{0}$ 
is the time derivative of the projection of $h$ onto the space of 
harmonic $1$--forms on $T^{3}$.

These last two equations with \eqref{eq:4.7}, 
the equation in the first point of 
\eqref{eq:2.1}, the bound $|\psi| \le B$ and Lemma \ref{thm:3.3} 
have the following consequence: 
The $L^{2}$ norm of $|P_{-}F_{A}|$ is uniformly bounded on each $t = 
\text{constant}$ 
slice of the cylinder $[-R+1, R-1] \times  T^{3}$ in terms of 
$|\omega|$ and the assumed $L^{2}$ bound of $F_{A}$ on $[-R, -R+4] \times 
T^{3}$ and on $[R-4, R] \times  T^{3}$. With the latter understood, 
standard elliptic regularity techniques find a uniform pointwise bound for 
$|F_{A}|$ at points where $t \in [-R+2, R-2]$. Of course, the bound on the 
$L^{2}$ norm on $t = \text{constant}$ 
slices provides one on all length $4$ cylinders 
in $X'$, and with this understood, the bound in \eqref{eq:4.5} 
follows from Lemma \ref{thm:3.4}.
\end{proof}

\sh{d) Vortices on the cylinder}

This section serves as a digression of sorts to describe the space of 
solutions to the vortex equation on $\R \times S^{1}$. These solutions 
will be used to describe the Seiberg--Witten moduli space on 
$\R \times T^{3}$. 

To describe the vortices, endow $\R \times  S^{1}$ with its standard 
product metric and with the complex structure from the identification via 
the exponential map with $\C^* = \C-\{0\}$. Fix a constant $r > 0$. A vortex 
solution is a pair $(v,\tau)$ of imaginary valued $1$--form and complex 
function which obey the conditions
\begin{equation}\begin{split}\label{eq:4.14}
\bullet\text{ }&{*}dv = -i r (1 - |\tau|^{2});\\
\bullet\text{ }&\bar{\partial}\tau + v_{0,1}\tau = 0;\\
\bullet\text{ }&(1 - |\tau|^{2}) \text{ is integrable}.
\end{split}
\end{equation}
Here, $v_{0,1}$ is the $(0, 1)$ component of $v$ in $T(\R \times S^{1})_{\C}$. 

Let $\ccc$ 
denote the set of solutions to \eqref{eq:4.14} modulo the equivalence relation 
that identifies $(v,\tau)$ with $(v',\tau')$ when $v' = v + \varphi 
d\varphi^{-1}$ and $\tau' = \varphi\tau$ whenever $\varphi $ is 
a smooth map from $\R \times S^{1}$ to the unit circle $S^{1}  \subset 
\C$. Topologize $\ccc$ as with $\mm$. 
That is, first topologize the solution set to 
\eqref{eq:4.14} with the subspace topology by embedding the latter in $i\cdot 
\Omega^{1}\times\Omega^{0}_{\C}\times[0,\infty)$, where 
the first two coordinates are 
$v$ and $\tau $, respectively, and the last is 
$r\int(1 - |\tau|^{2}){*}1$.                  
Then, give $\ccc$ the quotient topology.

Here are some facts (without proofs) about $\ccc$ (see, eg\ Section 2 of 
\cite{T2} 
or \cite{JT} for the proofs of similar assertions about vortices on $\C$):

{\sl
\begin{itemize}
\item  $\ccc$ {is the disjoint union of components}, 
$\{\ccc_{n}\}_{n=0,1,\ldots}$. {The component} $\ccc_{n}$ {is a 
manifold of complex dimension} $n$ {and is} {diffeomorphic to}  
$\zz_{n} = \{(y_{1},\ldots,y_{n}) \in \C^{n}:y_{n} \ne 0\}$ 
{by the map} $\Upsilon\colon \zz_{n} \to \ccc$ {which sends} $y\in\zz_{n}$ 
{to a solution of the form} 
\begin{equation}\label{eq:4.15}
(v,\tau) = (\bar{\partial}u - \partial u, e^{-u} p[y]), 
\end{equation}
{where} $p[y]$ {sends} $\eta\in\C^* = \R \times S^{1}$ 
{to the polynomial} $\eta^{n} + y_{1}\eta^{n-1} + \cdots 
+ y_{n}$. {Meanwhile}, $u$ {is the unique, 
real valued function on} $\R \times S^{1}$ {which obeys}
\begin{equation*}\begin{split}
&2i{*}\partial\bar{\partial}u=-r(1-e^{-2u}|p[y]|^{2});\\
&u=n\ln t+o(1)\text{ at points }(t,\cdot)\in\mathbb{R}\times 
S^{1}\text{ with }t\gg1;\\
&u = \ln|y_{n}| + o(1)\text{ at points }(t,\cdot)\in\mathbb{R} 
\times S^{1}\text{ with }t\ll-1.
\end{split}\end{equation*}

\item If $(v,\tau)\in\ccc_{n}$, then 
\begin{equation}\label{eq:4.17}
r\int(1-|\tau|^{2}) = 2\pi n.
\end{equation}

\item  $|\tau| \le 1$ with equality 
if and only if $|\tau|=1$ everywhere and $n = 0$.   

\item  {There is a constant} $\xi ${ which depends only 
on the vortex number and is such that}
\begin{equation}\label{eq:4.18}
(1-|\tau|^{2}) + r^{-1/2}|\nabla_{v}\tau| \le \zeta \exp[- 
(2r)^{1/2}\dist(\cdot,\tau^{-1}(0))].
\end{equation}

\item $\zz_{n}$ {is diffeomorphic to} $\C^{n-1}\times \C^*$.

\item There is a \textit{gluing map}, $\gaga$, 
that sends any finite product
$\zz_{n_{1}}\times\cdots\times\zz_{n_{k}}$ 
to the corresponding $\zz_{n_1 + \cdots + n_k}$; it is {defined by the 
requirement that} $p[\gaga(y_{n_{1}},\ldots,y_{n_{k}})] 
= p[y_{n_{1}}]\cdots p[y_{n_{k}}]$, where $p[\cdot]$ is as described in the 
first point 
above.

\item {The group} $\C^* = 
\R \times S^{1}$ acts on each 
$\ccc_{n}$ {as pull-back         
via its natural action on itself. In terms of 
the parameterization of a vortex solution as a point} $y = (y_{1}, 
\ldots,y_{n}) \in \zz_{n}$, {this action has} $\lambda \in 
\C^*$ {send} $y$ to $(\lambda^{-1}y_{1},\ldots, 
\lambda^{-n}y_{n})$. The action of an element in 
$\C^*$ on a vortex will be called a \textit{translation}.

\item {A vortex parametrized by} $y \in \zz_{n}$ {with} $|y_{n}| = 1$ 
{will be called} {\textit{centered}}. In this regard, note that 
$\ln|y_{n}|$ 
{equals the average of the} $t$ {coordinates of the zeros in} $\R 
\times S^{1}$ of the corresponding $\tau$. {The 
value of} $\ln|y_{n}|$ will be called the \textit{center} of the vortex.
\end{itemize}}

Vortex solutions on $\R \times S^{1}$ give Seiberg--Witten solutions on $X 
= \R \times  S^{1}  \times  T^{2}$ as follows: Take $\omega $ on $X$ to 
equal $r\cdot P_{+}\omega_{1}$, where $\omega_{1}$ is the standard 
volume form on $\R \times S^{1}$. Now let $(v,\tau)$ be a vortex 
solution and set 
$A = A_{0} + 2v$, $\psi = (\sqrt r \tau, 0)$,   
where $A_{0}$ 
is the product connection. Here, the bundle $S_{+}$ has been split into 
eigenspaces for Clifford multiplication by $\omega $.

\sh{\boldmath e) The moduli space for $\R \times T^{3}$}

This section constitutes a second digression to consider the moduli space of 
solutions to \eqref{eq:2.1} 
for the case $X = \R \times  T^{3}$. For this purpose, 
take $\omega $ in \eqref{eq:2.1} 
to be a non-zero, covariantly constant form on the 
whole of $X$. When considering the possibilities for $\mm$ 
in this case, note that 
as $X_{0} = [-1, 1] \times  T^{3}$, there is just one element in 
$\sss_{0}(X_{0})$, the trivial $\Spin^{\C}$ structure. Meanwhile, 
$\sss_{0}(X_{0},\partial X_{0}) = H^{2}_{\text{comp}}(\R\times 
T^{3})$. Moreover, cup product with a generator of 
$H^{1}_{\text{comp}}(R;\Z) = 
\Z$ provides an isomorphism, $\Xi\colon H^{1}(T^{3};\Z)\approx 
H^{2}_{\text{comp}}(\R \times  T^{3}; \Z)$ 
and so $\sss_{0}(X_{0},\partial X_{0}) 
= H^{1}(T^{3}; \Z)$.

With the preceding understood, consider:

\begin{prop}\label{thm:4.4} 
{Fix} $z \in H^{1}(T^{3}; \Z)$ {and then
the corresponding Seiberg--Witten 
moduli space} $\mm$ on $\R \times 
T^{3}$ consists of a single point (the orbit of a solution with 
$F_{A}\equiv 0$) {unless} 
\begin{equation*}
\omega = P_{+}(dt \wedge  \theta ), 
\end{equation*}
{where} $\theta $ is a constant $1$--form on $T^{3}$ 
{whose cohomology class is a positive multiple of}  $z$ {in} 
$H^{1}(T^{3};\Z)$. {When the latter is true, then} $\mm$ {is 
diffeomorphic to the vortex moduli space} $\ccc_{n}$, {where} $n$ 
{is the divisibility of} $z$ in $H^{1}(T^{3}; \Z)$. 
{In particular, this diffeomorphism arises from the fact that each 
point in} $\mm$ is gauge equivalent to $(A = A_{0} + 2a, \psi = \sqrt r 
(\alpha,0))$, where the pair $(a,\alpha)$ {is the 
pull-back from} $\R \times  S^{1}$ {of a vortex in} $\ccc_{n}$ 
via the map that identifies the $\R$ {factors while it fibers} $T^{3}$ 
{over} $S^{1}$ {so that the constant $1$--form that gives 
} $S^{1}$ length $1$ {pulls back as the constant $1$--form which 
represents} $n^{-1}z$ {in} $H^{1}(T^{3};\Z)$. 
\end{prop}

The remainder of this subsection is occupied with the following:

\begin{proof}[Proof of Proposition \ref{thm:4.4}] 
Use $\omega $ to decompose the bundle 
$S_{+}$ on  $\mathbb{R} \times T^{3}$ as $E\oplus E^{-1}$, where $E \to\R 
\times T^{3}$ is a complex line bundle. Note that $E$ must be topologically 
trivial since the restriction of its first Chern class to $T^{3}$ must 
vanish.

Now consider $(A,\psi)$. Write $\psi = (\alpha,\beta)$ to 
correspond with the splitting of $S_{+}$. Then, Witten's arguments from 
\cite{W1} 
for compact K\"ahler manifolds (using Lemma \ref{thm:3.5} 
to justify integration by 
parts) can be employed here to prove that $\beta = 0$ and that $\alpha $ is 
holomorphic with respect to that complex structure on $X$ that is defined by 
the flat metric and the self-dual form $\omega $. In addition, 
$iP_{+}F_{A}=(1 - |\alpha|^{2})\omega $.

Note that the maximum principle insures that $|\alpha| < 1$ everywhere 
unless $\alpha\equiv1$. In this case, the solution is gauge 
equivalent to the constant solution on $T^{3}$. With this understood, assume 
below that $|\alpha| < 1$. 

To proceed, note that $dF_{A} = 0$ so $dP_{-}F_{A} = -dP_{+}F_{A}$. 
Differentiate this last identity to find an equation of the form $\nabla^
*\nabla(iP_{-}F_{A}) + \zeta_{1}|\alpha|^{2}i 
P_{-}F_{A} = i\zeta_{2} P_{-}(d_{A}\alpha^*\wedge d_{A}\alpha)$, 
where $\zeta_{1,2} > 0$ are universal constants. A 
similar equation holds for $P_{+}F_{A}$ and comparing these two equations 
with the help of the maximum principle gives the pointwise bound 
$|P_{-}F_{A}| \le |P_{+}F_{A}|$. This last inequality, with the 
condition $z\bullet z = 0$, 
implies that $|P_{-}F_{A}| = |P_{+}F_{A}|$ which is 
equivalent to $F_{A}  \wedge  F_{A} = 0$. In fact, it follows now that $i 
P_{-}F_{A} = (1 - |\alpha|^{2})\sigma $, where $\sigma $ is a 
constant, anti-self dual $2$--form with norm equal to that of $\omega $. 
Furthermore, $d_{A}\alpha^*\wedge d_{A}\alpha$ is 
proportional to $F_{A}$ everywhere. (See (4.28--30) in \cite{T2}.) 
It also follows 
that $\alpha $ is holomorphic with respect to the complex structure on $\R 
\times  T^{3}$ that is defined by the given flat metric and the self-dual 
form $\omega $. 

Let $\omega_{1} = \omega + \sigma $. This form has square zero, and 
its kernel defines a $2$--dimensional distribution on $\R \times  T^{3}$ on 
which $F_{A}$ and $d_{A}\alpha$ both vanish. This distribution is 
also invariant under the complex structure on $\R \times  T^{3}$ because 
$d_{A}\alpha$ has zero projection onto the corresponding $T^{0,1}X$. 
Furthermore, since $\omega_{1}$ is constant, the resulting foliation is 
the image in $\R \times  T^{3}$ of a linear foliation of the universal 
covering space $\R \times \R^{3} 
= \C^{2}$ by complex lines. Moreover, 
as $|F_{A}|^{2}$ has 
finite integral over $\R \times  T^{3}$, the time 
coordinate (the $\R$ factor in $\R \times  T^{3}$) is constant on each 
leaf of the foliation. Furthermore, as $\alpha $ is holomorphic, its zero 
set, a union of leaves of the foliation, is a smooth, compact, codimension 
$2$, complex submanifold of $X$. 
This implies that each leaf of the foliation is 
a closed, linear torus in $T^{3}$. 

The proposition follows directly from the preceding remarks.
\end{proof}

\sh{f) Compactness in some special cases}

This subsection returns now to consider the compactification of the 
Seiberg--Witten moduli spaces for those $X_{0}$ from Proposition 
\ref{thm:2.4}. In this 
case, Propositions \ref{thm:4.2}--\ref{thm:4.4} 
can be employed to compactify the space $\mm_{s,m}$ 
from \eqref{eq:2.3} 
as a stratified space. The basic tool for this is Proposition 
\ref{thm:4.5}, 
below. The statement of Proposition \ref{thm:4.5} 
reintroduces the function $f$ on 
$H^{2}_{\text{comp}}(X;\mathbb{Z})$ from \eqref{eq:4.3}.

\begin{prop}\label{thm:4.5} 
{Fix} $({s}, z) \in 
\sss_{0}(X_{0},\partial X_{0})$, {let} $\mm$ denote the corresponding 
moduli space {and let} $\{c_{j}\}_{j=1,2,\ldots} \in \mm$ be 
a sequence with no convergent subsequences. {Then, the following 
exist:}
\begin{enumerate}
\item $z' \in \varsigma({s})$ {with} $z'\bullet\varpi 
< z \bullet \varpi $.
\item {A point} $c_{\infty}$ {in the corresponding} $\mm'$.
\item {For each component} $T = [0, \infty)  \times  T^{3}  \subset [0, 
\infty) \times \partial X_{0}$ {an element} $z_{\ee} 
\in  H^{1}(T^{3};\Z)$ {and a} 
{sequence} $\{o_{j}\}_{j=1,2,\ldots}$ {in the corresponding 
moduli space} $\mm^{\ee}$ {on} $\R \times  T^{3}$.
\end{enumerate}
{This data has the following significance:} {There is a 
subsequence of} $\{c_{j}\}$, hence relabled consecutively, such that: 
\begin{itemize}
\item $\{c_{j}\}$ {converges to} $c_{\infty}$ {on compact subsets of} 
$X$ in the $C^{\infty}$ {topology. This is to say that there exists 
a pair} $(A_{\infty},\psi_{\infty}) \in c_{\infty}$, and, for 
each index $j$, {a pair} $(A_{j},\psi_{j})$ on $c_{j}$ 
{such that if} $K \subset X$ {is a compact set and} 
$k$ is a positive integer, then
\begin{equation}\label{eq:4.20a}
\sum_{0\le p\le k}\left(|\nabla^{p}(A_{j} - A')| + |(\nabla 
_{A'})^{p}(\psi_{j} - \psi'_{j})|\right) < \varepsilon 
\end{equation}
at all points in $K$.
\item {Let} $T = [0, \infty)  \times  T^{3}$ {be 
a component of} $[0, \infty)  \times  \partial X_{0}$. {Given} $\varepsilon  > 
0$ {and a positive integer} $k$, {there exists} $R \ge 2$ 
{and, for each index} $j$, {there is a pair} $(A'_{j},\psi'_{j})$ 
{on} $o_{j}$ {such that}
\begin{equation}\label{eq:4.20b}
\sum_{0\le p\le k}\left(|\nabla^{p}(A_{j} - A')| + |(\nabla 
_{A'})^{p}(\psi_{j} - \psi'_{j})|\right) < \varepsilon
\end{equation} 
{at all points} $(t,\cdot)$ with $t \ge R$.
\item The vortex number which corresponds to $o_{j}$ is no 
greater than $\zeta \cdot f$, where $\zeta$ is 
independent of both $z$ and the sequence $\{c_{j}\}$.
\end{itemize}
\end{prop}

Here are two immediate corollaries of Propositions 
\ref{thm:4.4} and \ref{thm:4.5}:

\begin{prop}\label{thm:4.6} 
{Let} $X_{0}$ {and} $X$ be as 
in Proposition \ref{thm:2.4}. Now, suppose that the restriction of $\omega$ 
to each component of $[0, \infty)  \times  \partial X_{0}$ has the form 
$P_{+}(dt \wedge\theta)$, {where} $\theta $ {is a 
non-zero, covariantly constant $1$--form whose cohomology class is not 
proportional to an integral class in} $H^{1}(T^{3};\mathbb{R})$. For each 
$({s}, z) \in \sss_{0}(X_{0},\partial X_{0})$, {let} $\mm$ 
denote the resulting Seiberg--Witten moduli space. Then $\mm$ is 
compact.
\end{prop}

\begin{prop}\label{thm:4.7}
Suppose that $X$ {is isometric to 
the product metric on}  $S^{1}  \times M$, {where} $M$ {is an 
oriented, Riemannian $3$--manifold 
whose ends are isometric to} $[0, \infty)  \times 
T^{2}$. {In addition, suppose that} $\omega$ {is 
invariant under the evident} $S^{1}$ action and that on each end, 
$\omega = P_{+}(dt \wedge  \theta)$, {where} $\theta $ 
{is a covariantly constant $1$--form on} $S^{1}  \times 
T^{2}$ which is {\textit {not}} pulled back via 
projection to the $T^{2}$ {factor}. {For each} $({s}, 
z) \in \sss_{0}(X_{0}, \partial X_{0})$, {let} $\mm$ denote the 
corresponding Seiberg--Witten moduli space and let $\mm_{S} \subset \mm$ 
{denote the subset of} $S^{1}$--invariant orbits. Then 
$\mm_{S}$ {is compact}. 
\end{prop}

The remainder of this subsection contains the proofs of these 
propositions.

\begin{proof}[Proof of Proposition \ref{thm:4.6}] 
If the assertion were false, then 
Proposition \ref{thm:4.5} 
would find a solution on $\R \times T^{3}$ with $F_{A}\ne 
0$ identically. The latter is outlawed by Proposition \ref{thm:4.4}.
\end{proof}

\begin{proof}[Proof of Proposition \ref{thm:4.7}] 
If this assertion were false, Proposition 
\ref{thm:4.5} 
would find a non-trivial, $S^{1}$--invariant solution on $\R \times 
S^{1} \times  T^{2}$ which is not obtained from a vortex solution via 
a map which factors through the projection to $\R \times  T^{2}$. This is 
impossible, for if 
$(A, (\alpha, 0))$ has an $S^{1}$--invariant orbit under\linebreak
$C^{\infty}(\R \times  T^{3};S^{1})$, then $\alpha^{-1}(0)$ must be a 
union of $S^{1}$ orbits in $T^{3}$. 
\end{proof}

\begin{proof}[Proof of Proposition \ref{thm:4.5}]  
Proposition \ref{thm:4.2} describes each $c = (A, 
\psi) \in \{c_{j}\}$ on the ends of $X$. In particular, for each 
component $T \subset  [0, \infty)  \times  \partial X_{0}$ 
and each such $c$, there is 
the corresponding set $\{t_{i}\equiv t[T,c]_{i}\} \subset [0,\infty)$
that appears in \eqref{eq:4.5}. If there is no convergent subsequence, then 
Proposition \ref{thm:3.6} 
requires at least one end $T$ for which the sequence of sets 
$\left\{\{t_{i}[T,c]\}: c\in\{c_{j}\}_{j=1,2,\ldots}\right\}$ is not uniformly 
bounded. 
Even so, a subsequence of $\{c_{j}\}$ can be found for which the 
sets $\{t_{i}[T,c]\}$ for each fixed end $T$ all have the same number of 
elements as $c$ ranges over the subsequence. (Henceforth, all subsequences 
will be implicitly relabled by consecutive integers starting from $1$.) By 
passing again to a subsequence, these sets can be assumed to converge on 
compact domains in each end $T$. With this last point understood, then a 
limit, $c_{\infty}$, 
of a subsequence of $\{c_{j}\}$ is obtained using relatively 
standard compactness arguments. 

The sequence $\{o_{j}\}$ for a component $T \subset  [0, \infty)  \times 
\partial X_{0}$ 
is obtained by translating along $[0, \infty)$  to follow elements in 
$\{t[T,c]_{i}\}$ which do not stay bounded as $c$ ranges through the sequence 
$\{c_{j}\}$. In this regard, Propositions \ref{thm:4.2}, 
\ref{thm:4.3} and \ref{thm:4.4} together with the 
translation invariance of the equations in \eqref{eq:2.1} on 
$\R \times  T^{3}$ 
play the key role. Indeed, the construction of $\{o_{j}\}$ begins by 
obtaining a finite set of centered, limit vortices for each end by 
translating each $c\in\{c_{j}\}$ on the end and then, after passing to 
a subsequence, one follows each `clump' of energy. Here, the centers of 
these `clumps' on the end $T$ for the element $c  \in \{c_{j}\}$ are, by 
definition, obtained by first partitioning the set $\{t[T,c]_{i}\}$ which 
appears in \eqref{eq:4.4} 
into subsets whose elements are much closer to each other 
than to the other subsets of the partition. Then, the center of the clump 
subset is declared to be the average of the $t$ coordinates of the elements in 
the subset. 

After passing again to a subsequence, these `clump' partitions of the 
sets\linebreak
$\{t[T, c]_{i}\}$ can be assumed to produce the same number of subsets as $c$ 
ranges through $\{c_{j}\}$ and to be labeled for each such $c$ so that the 
following is true: First, the distance between subsets with different labels 
diverges as the index $j$ on $c_{j}$ tends to infinity. Second, the subsets 
with the same label have a fixed number of elements as $c$ ranges through 
$\{c_{j}\}$, and these elements can be themselves labeled so that the 
resulting ordered sets converge as the index $j$ tends to infinity.

One then translates the restriction of each $c\in\{c_{j}\}$ on an end 
$T$ so that the average $t$--coordinate of a given, labeled clump subset of 
$\{t[T, c]_{i}\}$ is $0$. After passing to a subsequence, the resulting 
translated sequence of Seiberg--Witten solutions will then converge strongly 
in the $C^{\infty}$ topology on compact domains in $\R \times  T^{3}$ to a 
solution to \eqref{eq:2.1} 
on $\R \times  T^{3}$. This limit is equivalent to a 
vortex solution as described by Proposition \ref{thm:4.4}. 
No generality is lost by 
translating this vortex in $\R$ so it is centered. 

By the way, the bound provided in Proposition \ref{thm:4.2} 
on the size of $\{t[T, 
c]_{i}\}$ explains the bound in Proposition \ref{thm:4.5} 
on the number of vortex 
solutions that arise.

In any event, a component, 
$T \subset  [0, \infty)  \times  \partial X_{0}$ provides 
a well defined, finite set of centered vortex solutions, each corresponding 
to one of the labels of the clump partition just described. This labeled set 
of limit vortex solutions can be characterized by the corresponding data 
$\{y^{(\alpha )}\}$, with each $y^{(\alpha )}\in\zz_{m}$ for $m = 
m_{\alpha} $. This characterization of the vortices is then used to 
construct, for each $c  \in \{c_{j}\}$, a vortex solution on $\R \times  
S^{1}$ which is obtained by gluing with the map 
$\gaga$ the translated (via the 
$\R$ 
factor in $\C^* = \R \times  S^{1}$) versions of the vortices from the limit 
set. Here, the particular translation for a vortex depends on the particular 
$c\in\{c_{j}\}$ and the particular clump label. To be precise, the 
translation is chosen so that the center in $\R$ of the translated vortex 
agrees with the average $t$--coordinate of the corresponding clump subset in 
$\{t[T, c]_{i}\}$. 

After the application of the gluing map $\gaga$, the result is a vortex which 
gives a Seiberg--Witten solution on $\R \times  T^{3}$ that is close to $c$, 
where the latter differs substantially from the limit $c_{\infty}$ 
and which is 
close to the trivial solution elsewhere.

With the preceding understood, the convergence assertion in \eqref{eq:4.20a} 
and \eqref{eq:4.20b} follows 
with standard elliptic regularity arguments for compact domains together 
with \eqref{eq:4.5}, \eqref{eq:4.18} and Lemma \ref{thm:3.4} 
to predict the form of each $c  \in
\{c_{j}\}$ on those parts of $[0, \infty)  \times  \partial X_{0}$ where the 
connection component of $c$ has small curvature.

The proof of Proposition \ref{thm:4.5} 
ends with an explanation for the fact that the 
solution $c_{\infty}$ 
sits in a moduli space defined by the original $\Spin^{\C}$ 
structure ${s}$ but with a class $z'\in\varsigma({s})$ 
with $z'\bullet \varpi < z \bullet \varpi $. The explanation starts 
with the observation that the convergence behavior described by 
\eqref{eq:4.20a} and \eqref{eq:4.20b}
insures that the original ${s}$ is also the $\Spin^{\C}$ structure for 
$c_{\infty}$. 
Keep in mind here that each $\{o_{j}\}$ gives solutions to \eqref{eq:2.1} on 
$\R \times  T^{3}$ whose $\Spin^{\C}$ structure is sent to zero by the map 
$c(\cdot)$ in \eqref{eq:1.3}. 
The explanation ends with the observation that the 
difference between the cup products of $z$ and $z'$ with $\varpi$ is a 
consequence of the final line in Proposition \ref{thm:4.1} 
when the convergence 
behavior in \eqref{eq:4.20b} is noted. 
\end{proof}

\section{Refinements for the cylinder}\label{sec:5}

The step from the compactness which is asserted in Proposition 
\ref{thm:4.6} to the 
assertions in Proposition \ref{thm:2.4} 
requires a substantial refinement of the 
asymptotic estimates from the preceding section. This section derives the 
required estimates for the behavior of solutions to \eqref{eq:2.1} 
on subcylinders in 
$\R \times  T^{3}$. 

\sh{\boldmath a) The operator $\dd_{c}$ when $X = \R \times T^{3}$}

The compactness study requires a more in depth study of the operator 
$\dd_{c}$ 
in the case where $X = \R \times  T^{3}$ and where $\omega $ in 
\eqref{eq:2.1} is a 
non-zero, covariantly constant form. The discussion here is broken into five 
steps.

\textbf{Step 1}\qua Fix $z \in  H^{1}(T^{3};\Z)$ and use 
the identification of the latter with\linebreak
$\sss_{0}(X_{0},\partial X_{0})$ to specify 
a Seiberg--Witten moduli space $\mm$. 
However, if $\mm$ is to contain more than the 
$F_{A} = 0$ solutions, it is necessary to further assume that $\omega= 
P_{+}(dt\wedge\theta)$, where $\theta $ is a covariantly constant 
form whose cohomology class in $H^{1}(T^{3};\Z)$ is proportional to $z$. Use 
the metric on $X$ and $\omega $ to define a complex structure, $J$. Note that 
$J\cdot dt = r\theta$, where $r \ne 0$, and note 
that ${*}\theta$ (a $2$--form on 
$T^{3}$) is invariant under the induced action of $J$ on $\Lambda^{2}T^*X$. 
Write $T^{0,1}X = \varepsilon_{0} \oplus \varepsilon_{1}$, where the 
first factor, $\varepsilon_{0}$, signifies the span of $\nu_{0} 
= dt - i r\theta$, and where the second factor is the orthogonal 
complement. Thus, the factor $\varepsilon_{1}$ has a covariantly 
constant, unit length frame $\nu_{1}$ with the property that the wedge of 
$\nu_{1}$ with its conjugate is proportional to ${*}\theta $. 
Meanwhile, the anti-canonical bundle is spanned by 
$\nu_{0}\wedge\nu_{1}$ and so is isomorphic to $\varepsilon_{01}  \equiv
\varepsilon_{0}\varepsilon_{1}$.

Given the preceding, the bundle $S_{+}$ can be written as $S_{+}= 
\varepsilon_{\C} \oplus \varepsilon_{01}$, where $\varepsilon_{\C}  
\to  T^{3}$ is a topologically trivial complex line bundle. Of course, 
$\varepsilon_{01}$ is also topologically trivial, but these lines in 
$S_{+}$ are distinguished by being the eigenbundles for Clifford 
multiplication by $\omega $. The point is that when $c\in\mm$, then the 
domain of $\dd_{c}$ consists of the $L^{2}_{1}$ sections over $\R\times 
T^{3}$ of
\begin{equation}\label{eq:5.1}
(\varepsilon_{0} \oplus \varepsilon_{\C}) \oplus (\varepsilon_{1} \oplus 
\varepsilon_{01}).
\end{equation}
Meanwhile, Clifford multiplication on the $\varepsilon_{\C}$ summand of 
$S_{+}$ identifies $S_{-}$ with $\varepsilon_{0} \oplus \varepsilon_{1}$ 
and so identifies the range of $\dd_{c}$ with the space of $L^{2}$ sections 
over $\R \times  T^{3}$ of
\begin{equation}\label{eq:5.2}
(\varepsilon_{\C} \oplus \varepsilon_{0}) \oplus (\varepsilon_{01} \oplus 
\varepsilon_{1}).
\end{equation}
Here, the factor $\varepsilon_{\C} \oplus \varepsilon_{01}$ corresponds to 
the $i \mathbb{R} \oplus \Lambda_{+}$ portion of $\dd_{c}(\cdot)$. 
In this regard, the 
real part of $\varepsilon_{\C}$ corresponds to the $i \mathbb{R}$ 
summand and the 
imaginary part to $i\mathbb{R} \omega   \subset  \Lambda_{+}$.

\textbf {Step 2}\qua The 
ordering of the factors and the 
placing of the parentheses              
in \eqref{eq:5.1} and \eqref{eq:5.2} 
have been chosen so as to 
decompose the domain and range of $\dd_{c}$ into direct sums, and this 
decomposition induces the following $2\times2$ block decomposition of 
$\dd_{c}$: 
\begin{equation}\label{eq:5.3}
\dd_{c} = \left(\begin{array}{cc} 
\Theta & \delta_F^{\dagger} \\ 
\delta_F & -\Theta^{\dagger}
\end{array} \right).
\end{equation}
Here, $\delta_{F}$ is a differential operator along the foliation of $M$ 
given by the kernel of $\theta $ which annihilates the defining constant 
sections of the first two summands in 
\eqref{eq:5.1} and whose adjoint annihilates 
those of the last two in \eqref{eq:5.2}. 
In this regard, remember that $c$ is defined 
from a vortex solution on $\R \times S^{1}$ via a projection of the 
form identity $ \times  \varphi $ from $\R \times  T^{3}$ to $\R \times 
S^{1}$. Here, $\varphi $ pulls back the standard $1$--form on $S^{1}$ to a 
multiple of $\theta $. The leaves of the foliation are the fibers of the map 
$\varphi $. Meanwhile, the operator $\Theta $ corresponds to a version of 
the operator which gives the linearized vortex equations. To be explicit 
here, 
\begin{equation}\label{eq:5.4}
\Theta(a,\lambda) = (\partial a + 2^{-1} r \bar{\tau}\lambda, \bar 
{\partial}_{v}\lambda + \tau a).
\end{equation}
Here, $r \equiv |\theta|$, $\partial_{v} \equiv \partial + v_{0,1}$ and 
$\partial = 
2^{-1}(\partial/\partial t - i\theta^*)$, where 
$\partial/\partial t$ denotes the tangent vector field to 
the lines $\R\times\text{point}$ while $\theta^*$ 
denotes the vector field which 
is metrically dual to the $1$--form $\theta $. The adjoint, 
$\Theta^{\dagger}$, of 
$\Theta $ is given by
\begin{equation}\label{eq:5.5}
\Theta^{\dagger}(b,\eta) = (-\bar{\partial}b + \bar{\tau}\eta,
-\partial_{v}\eta + 2^{-1} r \tau b).
\end{equation}

\textbf {Step 3}\qua The following lemma describes key 
properties of the operator $\Theta $.

\begin{lem}\label{thm:5.1} 
{Given a non-negative integer} $n$, let 
$(v,\tau)$ {be a vortex solution on} $\R \times  S^{1}$ 
{which has vortex number} $n$. {Then:}

\begin{itemize}
\item  {The operator} $\Theta $ is Fredholm from 
$L^{2}_{1}$ {to} $L^{2}$.
\item  {The} $L^{2}$ {kernel of} $\Theta $ is a 
$2n$--dimensional vector space of smooth forms.
\item  {The} $L^{2}$ {kernel of} $\Theta^{\dagger}$ 
{is empty}. {In fact, there exists a constant} $E$ {which 
depends only on} $|\theta|$ {and} $m$ {and is such that} $\| 
\Theta^{\dagger}w\|_{2}^{2} \ge E\| w \|_{2}^{2}$ {for all} 
$L^{2}_{1}$ {sections} $w$. 
\item  {By the same token,} $\|\Theta w\|_{2}^{2} \ge E \| 
w \|_{2}^{2}$ {for all} $L^{2}_{1}$ {sections} $w$ 
{which are}  $L^{2}$ {orthogonal to the kernel of} $\Theta $.
\item  {If} $w \in \kernel(\Theta)$, {then} $|(\nabla 
_{v})^{p}w| \le \zeta_{p} \| w \|_{2} \exp\left(-(2r)^{1/2} \dist( 
\cdot,\tau^{-1}(0))\right)$. {Here}, $\zeta_{p}$ {is 
independent of} $w$ {and} $(v,\tau)$.
\end{itemize}
\end{lem}

\begin{proof}[Proof of Lemma \ref{thm:5.1}] 
These assertions are all derived using the 
Weitzenb\"ock 
formula
\begin{equation}\label{eq:5.6}
\Theta\Theta^{\dagger}(b,\lambda)=(-\partial\bar{\partial}b
+2^{-1}r|\tau|^{2}b,
-\bar{\partial}_{v}\partial_{v}\eta+2^{-1}r|\tau|^{2}\eta),
\end{equation}
or the corresponding formula for $\Theta^{\dagger}\Theta$. The latter 
switches $\partial$ 
with $\bar{\partial}$ everywhere and adds an extra term which 
is proportional to $(1 - |\tau|^{2})$. As this term is small at large 
distances from $\tau^{-1}(0)$ (see \eqref{eq:4.18}), the formula for $\Theta
^{\dagger}\Theta $ implies that $\Theta $ is Fredholm on the $L^{2}_{1}$ 
completion of its domain as a map into the $L^{2}$ completion of its range. 
With this understood, \eqref{eq:5.6} 
implies that the cokernel of $\Theta $ (which is 
the kernel of $\Theta^{\dagger}$) is empty. Furthermore, 
\eqref{eq:5.6} plus the last 
point in \eqref{eq:4.18} 
implies the estimate in the third point of the lemma. The 
latter estimate plus the Fredholm alternative implies the estimate in the 
fourth point. In this regard, remark that when $w$ is $L^{2}$--orthogonal to 
the kernel of $\Theta $, then $w = \Theta^{\dagger}g$ 
for some element $g$. Thus, 
\begin{equation*}
\|\Theta w\|_{2} = \|\Theta\Theta^{\dagger}g\|_{2} \ge 
\|\Theta^{t}g\|_{2}^{2}/\|g\|_{2} \ge E \|\Theta^{t}g\|_{2} = E 
\| w \|_{2}.
\end{equation*}

The fact that the kernel dimension is $2n$ can be proved as follows: 
Differentiate the map which associates to $y \in \C^{m}$ a vortex $(v, 
\tau)$. (See \eqref{eq:4.15}.) Each such directional derivative gives an 
independent element in $\kernel(\Theta )$. Conversely, every element in 
$\kernel(\Theta)$ can be integrated to give an element in the tangent space 
to $\zz_{n}$. This follows from the vanishing of the kernel of 
$\Theta^{\dagger}$.

For the exponential decay estimate, consider first the $C^{0}$ case. If $w 
\in \kernel(\Theta)$, then \eqref{eq:4.18} and the Weitzenb\"ock formula for 
$\Theta^{\dagger}\Theta$ imply the following: Where the distance from $\tau 
^{-1}(0)$ is greater than $1$, the norm of $w$ obeys
\begin{equation}\label{eq:5.8}
d^*d|w| + 2r |w| \le \upsilon \cdot |w|,
\end{equation}
where $\upsilon $ is a function on $\R \times S^{1}$ which obeys 
\begin{equation}\label{eq:5.9}
\upsilon \le \zeta \exp\left(-\sqrt r \dist(\cdot,\tau^{-1}(0))/\zeta\right), 
\end{equation}
with $\zeta $ a universal constant. The control of $|w|$ where the distance to 
$\tau ^{-1}(0)$ is less than $1$ comes via standard elliptic regularity 
which finds a constant $\xi_{1}$ such that
\begin{equation*}
|w| \le \xi_{1} \|w\|_{2} 
\end{equation*}
at all points. Here, $\xi_{1}$ depends only on $r$. With the preceding 
understood, let $\{t_{i}\}$ denote the time coordinates of the $n$ points in 
$\tau^{-1}(0)$. An application of the comparison principle to \eqref{eq:5.8} 
now 
yields the following: Fix $\rho < (2r)^{1/2}$ and there exists $\zeta 
_{\rho} $ which is independent of $(v,\tau)$ and such that
\begin{equation}\label{eq:5.11}
|w| \le \zeta_{x} \| w \|_{2} \sum_{i}\exp(-\rho |t - t_{i}|).
\end{equation}

To obtain the analog of \eqref{eq:5.11} 
for the case where $\rho = (2r)^{1/2}$, 
introduce the Green's function $G(\cdot, z')$ for the operator $d^*d + 2r$ 
with a pole at $z'\in\R\times S^{1}$. An application of the 
comparison principle to the operator $d^*d + 2 r$ finds a constant $\kappa 
_{0}$ which makes the following assertion true: If the $\R$ coordinate of 
$z'$ 
is $s$, that of $z$ is $t$, and also $|t - s|\ge\zeta\ge1$, then 
\begin{equation}\label{eq:5.12}
0 < G(z, z') \le \kappa_{0} \exp(-(2r)^{1/2}|t - s|). 
\end{equation}
Furthermore, the absolute values of the derivatives of $G(\cdot,z')$ at 
such $z$ enjoy similar upper bounds. 

With the Green's function in hand, multiply both sides of \eqref{eq:5.8} 
by $G(z,\cdot)$ 
and integrate both sides of the result over the region, $U$, where 
the distance to the set $\tau^{-1}(0)$ is at least one. This operation 
produces an integral inequality which can be further manipulated to yield 
the following bound:
\begin{multline}\label{eq:5.13}
|w(z)| \le \xi_{2}\| w \|_{2} \sum_{i} \exp(-(2r)^{1/2}|t - t_{i}|)\\
+ \xi_{2}\| w \|_{2}\int_{U} ds \exp(-(2r)^{1/2} |t - 
s|) \sum_{i} \exp(-(2r)^{1/2}(1+\delta)|s - t_{i}|).
\end{multline}
Here, $\delta  > 0$ is a universal constant, while $\xi_{2}$ depends 
only on $r$. To explain, 
the first term in \eqref{eq:5.13}   
is due to the integration by 
parts boundary term that arises when the operator $d^*d$ in 
\eqref{eq:5.8} is moved from 
$w$ to $G(\cdot,z)$. Meanwhile, the second term in \eqref{eq:5.13} 
comes from the 
integral of $G(z,\cdot)\upsilon\cdot|w|$. In this regard, \eqref{eq:5.9} is 
used to bound $\upsilon $, \eqref{eq:5.11} 
with $\rho $ very close to $(2r)^{1/2}$ is 
used to bound $|w|$, and \eqref{eq:5.12} is used to bound $G(z,\cdot)$. 

The $\rho = (2r)^{1/2}$ version of \eqref{eq:5.11} 
follows directly from \eqref{eq:5.13} 
since the second term in this equation is not greater than $\xi_{3}\| w 
\|_{2}\sum_{i} \exp(-(2r)^{1/2}|t - t_{i}|)$, with $\xi_{3}$ 
depending only on the parameter $r$.

The proof for the asserted bounds on the higher derivatives of $w$ is obtained 
by first differentiating the equation $\Theta w = 0$ say, $p$ times, and then 
writing the latter as an equation of the form $\Theta(\nabla^{p}w) = 
\text{lower order derivatives of }w$. 
The preceding argument for the $C^{0}$ bound 
can be copied to obtain the desired estimates.
\end{proof}

\textbf {Step 4}\qua It is important to note that there is one 
particularly canonical element in the kernel of $\Theta $, this being
\begin{equation}\label{eq:5.14}
\pi_{c}=(2^{-1}r(1-|\tau|^{2}),\partial_{v}\tau).
\end{equation}
Viewed as a section over $\ccc_{n}$ of the complex tangent space, the vector 
$\pi_{c}$ generates the translation induced $\C^*$ action. 

Note that the pointwise bound in the last assertion of Lemma \ref{thm:5.1} 
is sharp 
for $\pi_{c}$. The following lemma makes a precise statement:

\begin{lem}\label{thm:5.2} 
{Given the form} $r \equiv |\theta|$ and a positive integer $n$, 
there exists a constant 
$\zeta\ge 1$ such that when $(v,\tau)\in\ccc_{n}$ 
then
\begin{equation}\label{eq:5.15}
(1 - |\tau|^{2}) \ge \zeta^{ - 1} \sum_{j}\exp( - (2r)^{1/2}|t - 
t_{j}|),
\end{equation}
{where} $\{t_{j}\}$ {denote the} $t$--coordinates of the 
zeros of $\tau $.
\end{lem}

\begin{proof}[Proof of Lemma \ref{thm:5.2}]  
Let $x = (1 - |\tau|^{2})$. Then, by 
virtue of the vortex equations in \eqref{eq:4.14}, this function obeys
\[
d^*dx + 2 r x = 4 |\partial_{v}\tau |^{2} + 2 r |x|^{2};\]
thus
\begin{equation}
d^*dx + 2 r x \ge 0\label{eq:5.17}
\end{equation}
everywhere. The first consequence of 
\eqref{eq:5.17} comes via the maximum principle, 
this being the previously mentioned fact that $(1 - |\tau |^{2}) > 0$ as 
long as $n \ge 0$. This last point can be parlayed to give a universal lower 
bound for $(1 - |\tau|^{2})$ at points at fixed distance from $\tau ^{ 
- 1}(0)$. In particular, the following is true:
{\sl\begin{equation}\begin{split}\label{eq:5.18}
&\text{Given }r = |\theta |\text{, the vortex number }n\text{ and 
also }\rho > 0\text{, there exists}\\
&\xi > 0\text{ such that when }
(v,\tau)\in\ccc_{n}\text{ then }(1 - |\tau|^{2}) \ge \xi 
\text{ at points}\\
&\text{with distance }\rho\text{ or less from }\tau^{-1}(0).
\end{split}\end{equation}}

Accept \eqref{eq:5.18} 
for the moment to see its application first. In this regard, 
\eqref{eq:5.18} 
is applied here to supply a uniform, positive lower bound, $\xi $, 
for $(1 - |\tau|^{2})$ on the 
constant $t$ circles where $|t - t_{j}| \le 1$ 
for at least one $t_{j}$. With 
such a bound in place, reapply the maximum
principle to \eqref{eq:5.17} 
but use $x = (1 - |\tau|^{2}) - n^{-1}\xi 
\sum_{j} \exp( - (2r)^{1/2}|t - t_{j}|)$ 
and restrict to points in $\R 
\times  S^{1}$ where $|t - t_{j}| \ge 1$ 
for all $t_{j}$. The 
resulting conclusion (that $x \ge 0$) and \eqref{eq:5.18} 
together give \eqref{eq:5.15}.

To justify \eqref{eq:5.18}, 
consider the ramifications were the claim false. In 
particular, there would exist a sequence $\{(v_{i},\tau_{i})\} \in 
\ccc_{n}$, and a corresponding sequence of pairs of points $\{(z_{j}, 
z'_{j})\}\subset\R \times  S^{1}$, where $\tau_{j}(z_{j}) = 
0$, $\dist(z_{j}, z'_{j}) \le \rho $ and $\lim_{j \to \infty} |\tau 
_{j}(z'_{j})| = 1$. After translating each $(v_{j},\tau_{j})$ 
appropriately, all of the points $\{z_{j}\}$ can be taken to be a fixed 
point $z\in\R \times S^{1}$. Meanwhile, the sequence $\{z'_{j}\}$, 
has a convergent subsequence with limit $z'$ whose distance is $\rho $ or less 
from $z$. Also, a subsequence $\{(v_{j},\tau_{j})\}$ converges strongly 
in the $C^{\infty}$ topology on compact domains in $\R\times S^{1}$ to a 
vortex solution, $(v,\tau)$, although the latter may lie in some $\ccc_{n'}$ 
for $n' \le n$. Indeed, the $C^{1}$ convergence of $\tau $ and $C^{0}$ 
convergence of $v$ follows directly from \eqref{eq:4.17} 
and \eqref{eq:4.18} by appeal to the 
Arzela--Ascoli theorem; convergence in $C^{k}$ can then be deduced by 
differentiating the vortex equations. 
In particular, the $C^{1}$ convergence 
here implies that $\tau(z) = 0$ and 
$|\tau(z')| = 1$. However, these 
two conclusions are not compatible; as argued previously, $|\tau| < 1$ 
everywhere if $|\tau|$ is less than $1$ at any point.
\end{proof}

\textbf {Step 5}\qua This last step uses the results of Lemma 
\ref{thm:5.1} to draw conclusions about $\dd_{c}$. 
The latter are summarized by:

\begin{lem}\label{thm:5.3} 
{Suppose} $c = (A, \psi) \in \mm$ 
{comes from a vortex solution}, $(v, \tau) \in \ccc_{n}$ 
on $\C^*$ {via a fibration map} $\varphi $, {from} $\R \times 
(S^{1}\times S^{2})$. {Then:}

\begin{itemize}
\item  {The kernel of} $\dd_{c}$ {is the} 
$2n$--dimensional vector space which consists of configurations $((a, 
\lambda),(0,0))$ {which are annihilated by both} $\Theta $ 
{and differentiation along the fibers} $F$ {of} $\varphi $. 
{In particular,} $(a,\lambda)$ comes from the kernel of 
$\Theta$ {on} $\C^* = \R \times S^{1}$ {via pull-back by 
the map} $\varphi $. 

\item  {The kernel of} $\dd_{c}^{\dagger}$ is the 
$2n$--dimensional vector space which consists of configurations $((0, 
0),(a,\lambda))$ {which are in the kernel 
of} $\Theta $ 
{and differentiation along the fibers of} $\varphi$.  
{These 
also come from the kernel of} $\Theta $ {on} $\R\times S^{1}$ 
{via pull-back by the map} $\varphi$.

\item  {If} $w$ {is an} $L^{2}_{1}$ {section of 
\eqref{eq:5.1} which is} $L^{2}$--orthogonal to the kernel of $\dd_{c}$, 
{then}
\begin{equation}\label{eq:5.17'}
\|\dd_{c}w\|_{2} \ge E\| w \|_{2}, 
\end{equation}
{where} $E$ {depends only on} $\theta $ {and} $n$ 
{and is, in particular, independent of} $(v,\tau)$.  

\item {Likewise}, {if} $w$ {is an} $L^{2}_{1}$ 
{section of} \eqref{eq:5.2} which is $L^{2}$--orthogonal to 
the kernel of $\dd_{c}^{\dagger}$, {then} $\|\dd_{c}^{\dagger}w\|_{2} \ge E 
\| w \|_{2}$.
\end{itemize}
\end{lem}

\begin{proof}[Proof of Lemma \ref{thm:5.3}] 
To consider the assertion about the kernel of 
$\dd_{c}$, write $w$ in $2$--component 
form with respect to the splitting in \eqref{eq:5.1} 
as $(\alpha,\beta)$, where each of $\alpha $ and $\beta $ also have two 
components. Then, 
\begin{equation*}
\|\dd_{c}w\|_{2}^{2} = \|\Theta\alpha\|_{2}^{2} + \|\bar 
{\partial}_{F}\alpha\|_{2}^{2} + \|\Theta^{\dagger}\beta 
\|_{2}^{2} + \|\bar{\partial}_{F}^{\dagger}\beta\|_{2}^{2}. 
\end{equation*}
Integration by parts along the fibers of $\varphi $ equates 
$\|\bar{\partial}_{F}\alpha\|_{2}^{2} = 4^{-1}\|d_{F}\alpha\|_{2}^{2}$, and 
$\|\bar{\partial}_{F}^{\dagger}\beta\|_{2}^{2} = 4^{ - 1} \| 
d_{F}\beta\|_{2}^{2}$. Now, write $\alpha = \alpha_{0} + 
\alpha_{1}$, where $\alpha_{0}$ is constant along each fiber of 
$\varphi $, and where $\alpha_{1}$ is $L^{2}$--orthogonal to the constants 
along each fiber of $\varphi $. It then follows that
\begin{equation}\label{eq:5.19}
\|\dd_{c}w\|_{2}^{2} \ge E'(\|\alpha_{1}\|_{2}^{2} + \| 
\beta_{1}\|_{2}^{2}) + \|\Theta\alpha_{0}\|_{2}^{2} 
+ \|\Theta^{\dagger}\beta_{0}\|_{2}^{2}. 
\end{equation}
The assertions that concern the kernel of $\dd_{c}$ follow from 
\eqref{eq:5.19} using 
Lemma \ref{thm:5.1}. 

An analogous argument proves the assertions in the lemma that concern the 
kernel of $\dd_{c}^{\dagger}$.
\end{proof}

\sh{\boldmath b) Decay bounds for $\kernel(\dd_{c})$ when $c\in 
\mm_{P}$}

The proof of Proposition \ref{thm:2.4} 
requires a refinement of Proposition \ref{thm:4.3}'s 
decay estimates along a finite cylinder $[-R, R] \times T^{3}$. In 
particular, \eqref{eq:4.5} 
establishes that a solution $(A,\psi)$ to \eqref{eq:2.1} on $[-R, 
R] \times T^{3}$ with $|F_{A}|$ small at all points is exponentially 
close in the middle of the cylinder to a $T^{3}$ solution, $(A_{0},\psi 
_{0}) \in\mm_{P}$. However, the bound in \eqref{eq:4.5} was not concerned with 
the size of the decay constant in the exponential. This subsection supplies 
a precise estimate for the constant $\zeta_{0}$ which appears in 
\eqref{eq:4.5}. 

For this purpose, consider the operator $\dd_{c}$ from \eqref{eq:5.3} 
when $c= 
(A_{0},\psi_{0})$ is the pull-back to $[-R, R] \times T^{3}$ of a 
solution on $T^{3}$ which defines $\mm_{P}$. In this case, write the 
operator $\dd_{c}$ as
\begin{equation}\label{eq:5.20}
\dd_{c} = 2^{-1}\cdot(\partial/\partial t + \oo).
\end{equation}
Here, $\oo$ is a $t$--independent, symmetric, first order operator which 
differentiates along the $T^{3}$ directions. Moreover, there are natural 
trivializations of the summands of $S_{+}=\varepsilon_{\C}\oplus
\varepsilon_{01}$ and $T^{0,1}X$ which make $\oo$ a constant 
coefficient operator. In this regard, the trivialization of $\varepsilon 
_{\C}$ makes $\psi_{0}$ the constant section with vanishing imaginary 
part and positive real part, while that of $T^{0,1}X$ is as described prior 
to \eqref{eq:5.1}. 

With the afore-mentioned trivialization understood, then Fourier transforms 
can be employed to investigate the spectrum of $\oo$. In particular, this 
spectrum is a nowhere accumulating subset of $(-\infty, -(2r)^{1/2}]\cup 
[(2r)^{1/2},\infty)$ 
which is unbounded in both directions, and invariant under 
multiplication by $-1$. Furthermore, the minimal, positive eigenvalue $E_{0}
= (2r)^{-1/2}$ is degenerate, with two eigenvectors over $\C$, these being the 
sections
\begin{equation}\label{eq:5.21}
s^{+}_{+} = (\sqrt r, \sqrt 2, 0, 0)\text{ and }s^{+}_{-} = (0, 0, \sqrt 2, 
-\sqrt r) 
\end{equation}
of \eqref{eq:5.1}. 
Meanwhile, the eigenvalue $-E_{0}$ also has two eigenvectors over 
$\C$, the sections 
\begin{equation}\label{eq:5.22}
s^{-}_{+} = (\sqrt r, -\sqrt 2, 0, 0)\text{ and }s^{-}_{-} = 
(0, 0, \sqrt 2, \sqrt r) 
\end{equation}
of \eqref{eq:5.1}. 

Any section $w$ of \eqref{eq:5.1} 
over $[-R, R] \times T^{3}$ which is annihilated 
by $\dd_{c}$ has the form
\begin{equation}\label{eq:5.23}
w = \sum_{E>0} e^{-E\cdot(t+R)}\cdot s^{+}_{E} + \sum_{E>0} 
e^{-E\cdot(R-t)}\cdot s^{-}_{E},
\end{equation}
where $s^{+}_{E}$ is an eigenvector of $\oo$ with eigenvalue $E$, and where 
$s^{-}_{E}$ is likewise an eigenvector, but with eigenvalue $-E$. Note that 
$\sum_{E>0}\|s^{+}_{E}\|_{2,T}^{2} + \sum_{E>0}\|s^{-}_{E}\|_{2,T}^{2} 
\le \zeta\cdot\| w \|_{2}^{2}$. Here, $\|\cdot\|_{2,T}$ 
denotes the $L^{2}$ norm on $T^{3}$ and $\|\cdot\|_{2}$ 
denotes the $L^{2}$ norm over $[-R, R] \times T^{3}$.

Equation \eqref{eq:5.23} is a linear version of the following:

\begin{lem}\label{thm:5.4} 
{Suppose that} $\upsilon $ is a 
homomorphism over $[-R, R] \times  T^{3}$ {from the bundle in 
\eqref{eq:5.1} 
to that in \eqref{eq:5.2} which obeys the bound} $|\upsilon| \le \zeta 
\cdot e^{-(R-|t|)/\zeta}$, {where} $\zeta  > 0$ {is a 
constant. Then, there are constants} $z > (2r)^{1 /2}$ {and} $\zeta'$ 
which depend only on $\zeta $ {and which have the 
following significance: Let} $w$ {be a section over} $[-R, R] \times  
T^{3}$ {of \eqref{eq:5.1} which obeys} 
\begin{equation}\label{eq:5.24}
\dd_{c}w + \upsilon w = 0.
\end{equation}
{Then} $w =w_{0} + w_{1}$ with
\begin{equation}\begin{split}\label{eq:5.25}
&\| w_{1} \|_{2,T}\big|_{t} \le \zeta' K e^{-z\cdot(R-|t|)}, 
\text{ where }K = \sup_{t \in [-R,-R+1]\cup [R-1,R]}|w|\text{ and}\\
&\begin{split}
w_{0} =&\exp(-(2r)^{1/2}(R + t))(u^{+}_{+} s^{+}_{+} + 
u^{+}_{-} s^{+}_{-}) +\\
&\exp(-(2r)^{1/2}(R - t))(u^{-}_{+}s^{-}_{+} + u^{-}_{-} s^{-}_{-}),
\text{ where }u^{\pm}_{\pm}\text{ are constants.}\end{split}
\end{split}\end{equation}
\end{lem}

\begin{proof}[Proof of Lemma \ref{thm:5.4}] 
Let $\Pi_{+}$ denote the $L^{2}$--orthogonal 
projection (on $T^{3}$) onto the span of the eigenvectors of $\oo$ with 
eigenvalue $E > (2r)^{1/2}$. Meanwhile, let $\Pi_{-}$ denote the 
corresponding 
projection onto those eigenvectors with eigenvalue $E < -(2r)^{1/2}$. Let 
$f_{+}(t)$ denote the $L^{2}$ norm of the time $t$ version of $\Pi_{+}w$, and 
likewise define $f_{-}(t)$. Then $f_{\pm}$ obey a pair of coupled 
differential inequalities of the form
\begin{equation}\begin{split}\label{eq:5.26}
&(\partial/\partial t + E_{2})f_{+} \le \zeta e^{-(R-|t|)/\zeta} 
\left(f_{+} + f_{-} + K \exp(-(2r)^{1/2} (R - |t|))\right);\\
&(\partial/\partial t - E_{2})f_{-} \ge -\zeta e^{-(R-|t|)/\zeta} 
\left(f_{+} + f_{-} + K \exp(-(2r)^{1/2} (R - |t|))\right).
\end{split}
\end{equation}
Here, $E_{2} > (2r)^{1/2}$ is the second smallest positive eigenvalue of 
$\oo$. 
This last equation can be integrated (after some algebraic manipulations) 
to obtain the bound 
\begin{equation*}
(f_{+} + f_{-})\big|_{t} \le \zeta' (f_{+}\big|_{-R} + f_{-}\big|_{R} + 
K) e^{-z'\cdot(R-|t|)}, 
\end{equation*}
where $z' > (2r)^{1/2}$ and $\zeta' > 0$ depends only on $\zeta $.

Meanwhile, write $w = b^{+}_{+} s^{+}_{+} + b^{+}_{-} 
s^{+}_{-} + b^{-}_{+} s^{-}_{+} + b^{-}_{-} s^{-}_{-} + 
\Pi_{+}w + \Pi_{-}w$ and consider the equations for $\mathbf{b}^{+} = 
(b^{+}_{+},b^{+}_{-})$ and for the corresponding $\mathbf{b}^{-}$. 
In particular, these equations have the form
\begin{equation*}\begin{split}
&(\partial/\partial t + 2r)\mathbf{b}^{+} = g^{+}\cdot\mathbf{b}^{+} + 
g^{-}\cdot\mathbf{b}^{-} + \upsilon^{+}\text{ and}\\
&(\partial/\partial t - 2r)\mathbf{b}^{-} = h^{+}\cdot\mathbf{b}^{+} + 
h^{-}\cdot\mathbf{b}^{-} + \upsilon^{-},
\end{split}
\end{equation*}
where $|g^{\pm}| + |h^{\pm}| \le \zeta e^{-(R-|t|)/\zeta}$ and 
$|\upsilon^{\pm}| \le \zeta e^{-z\cdot(R-|t|)}$ with $z > 
(2r)^{1/2}$. Integrating these last two equations gives
\begin{equation*}\begin{split}
&\mathbf{b}^{+}\big|_{t} = \mathbf{b}^{+}\big|_{-R}\cdot 
\exp(-(2r)^{1/2}(t + R)) + \mathbf{c}^{+}\text{ and}\\
&\mathbf{b}^{-}\big|_{t} = \mathbf{b}^{-}|_{R}\cdot 
\exp(-(2r)^{1/2} (R - t)) + \mathbf{c}^{-},
\end{split}
\end{equation*}
where $|\mathbf{c}^{\pm}| \le \zeta' K e^{-z'\cdot(R - |t|)}$ with $z' > 
(2r)^{1/2}$ and $\zeta'$ depending only on $\zeta $.
\end{proof}

By way of an application, choose a vortex $(v,\tau)$ on $\R\times 
S^{1}$ and use the latter to define the gauge orbit, $c$, of a solution to 
\eqref{eq:2.1} 
on $\R \times  T^{3}$. The element $\pi_{c}$ in \eqref{eq:5.14} can be viewed 
as either an element, $\pi_{c+}$, in the kernel of $\dd_{c}$ or as an element, 
$\pi_{c-}$, in $\cokernel(\dd_{c})$. 
Here, $\pi_{c+}$ is the section of \eqref{eq:5.1} whose 
first two components are those of $\pi_{c}$ and whose second two are zero. 
Meanwhile, $\pi_{c-}$ is the section of \eqref{eq:5.2} 
whose first two components are 
zero and whose second two are the components of $\pi_{c}$. (The applications 
below only use $\pi_{c-}$.) 
As there are points on the gauge orbit $c$ which are 
asymptotic as $t \to \pm \infty$ 
to a solution $(A_{0}, \psi_{0})$ which 
defines $\mm_{P}$, 
Lemma \ref{thm:5.4} can be applied to $\pi_{c\pm}$ with the 
following effect:

\begin{lem}\label{thm:5.5} 
{Let} $(v,\tau) \in \ccc_{n}$ {and 
use the latter, as instructed 
in Proposition \ref{thm:4.4}, 
to define the gauge orbit,} 
$c$, of a solution to 
\eqref{eq:2.1} on $\R \times T^{3}$. 
{Define} $\pi_{c\pm}$ as in the preceding paragraph.
\begin{itemize}
\item {Let} $(A_{+},\psi_{+})$ {denote a point on the 
orbit} $c$ {such that} $|A_{+} - A_{0}| +\linebreak |\psi_{+} - \psi 
_{0}| \le \zeta e^{-t/\zeta}$ {at all points} $(t,\cdot)$
with $t > 0$. {Here,} $\zeta  > 0$. Use $\psi 
_{0}$ to define the constant real section of the first summand of 
$S_{+} = \varepsilon_{\C} \oplus \varepsilon_{0}$. {Then, as} $t\to\infty$,
\begin{equation}\begin{split}\label{eq:5.30}
&\pi_{c+} = \exp(-(2r)^{1/2} t) u^{+} s^{+}_{+} + \oo(e^{-zt}) 
\text{ and}\\
&\pi_{c-} = \exp(-(2r)^{1/2} t) u^{+} s^{-}_{-} + 
\oo(e^{-zt}).
\end{split}\end{equation}
\item {Let} $(A_{-}, \psi_{-})$ denote a point on the 
orbit $c$ such that $|A_{-} - A_{0}| +\linebreak |\psi_{-} - \psi 
_{0}| \le \zeta  e^{t/\zeta}$ at all points $(t, \cdot 
)$  with $t < 0$. Here, $\zeta  > 0$ also. 
Use $\psi_{0}$ to define the constant real section of 
the first summand of $S_{+} = \varepsilon_{\C} \oplus\varepsilon_{0}$.
Then, as $t \to  -\infty$,
\begin{equation}\begin{split}\label{eq:5.31}
&\pi_{c+} = \exp((2r)^{1/2}t) u^{-} s^{-}_{+} + \oo(e^{zt})\text{ and}\\ 
&\pi_{c-} = \exp((2r)^{1/2}t) u^{-} s^{+}_{-} + \oo(e^{zt}).
\end{split}\end{equation}
\end{itemize}
{In} {these equations,} $s^{\pm}_{\pm}$ {are 
given in} \eqref{eq:5.21} {and} \eqref{eq:5.22}, and $z > (2r)^{1/2}$ 
{is a constant which depends only on} $r$ and the constant 
$\zeta $, {but otherwise is independent of the chosen vortex.} 
{Meanwhile,} $u^{\pm}  > 0$ {are constants which depend on the 
chosen vortex} $c = (v, \tau)$ {and which obey}
\begin{equation}\begin{split}\label{eq:5.32}
&u^{+} \ge \xi^{-1}\cdot \sum_{j} \exp((2r)^{1/2}t_{j})\text{ and}\\
&u^{-} \ge \xi^{-1}\cdot \sum_{j} \exp(-(2r)^{1/2}t_{j}),
\end{split}
\end{equation}
{where} $\{t_{j}\}$ {are the} $t$--coordinates of the 
zeros of $\tau $ {and where} $\xi $  depends only on the 
vortex number.
\end{lem}

\begin{proof}[Proof of Lemma \ref{thm:5.5}] 
The assertions of \eqref{eq:5.30} and \eqref{eq:5.31} about 
$\pi_{c+}$ follow directly from Lemma \ref{thm:5.4}. 
Meanwhile, the assertions about 
$\pi_{c-}$ follow from Lemma \ref{thm:5.4} 
after changing $t$ to $-t$ in \eqref{eq:5.20}. Then, given 
these assertions, the bounds in \eqref{eq:5.32} follow from Lemma 
\ref{thm:5.2}.
\end{proof}

\sh{c) More asymptotics for solutions on a cylinder}

Proposition \ref{thm:5.6}, 
below, constitutes a second application of Lemma \ref{thm:5.4}, here 
to the asymptotics of a solution to \eqref{eq:2.1} 
on the cylinder $X = [-R - 2, R + 
2] \times  T^{3}$ with the form $\omega  = P_{+}(dt \wedge 
d\theta)$, where $\theta $ is non-zero and covariantly constant. As in 
previous subsections, the lemma below takes $(A_{0},\psi_{0})$ to be 
the pull-back to $X$ of a solution to \eqref{eq:3.1} 
which defines $\mm_{P}$ and it takes 
$r \equiv |\theta |$. 

\begin{prop}\label{thm:5.6} 
{The metric on} $T^{3}$ {and the 
form} $\omega $ {determine a constant} $\zeta\ge1$ that has 
the following significance: Let $R \ge 4$ and suppose that 
$(A,\psi)$ {obeys} \eqref{eq:2.1} and the assumptions of Lemma 
\ref{thm:3.4} on $[-R - 2, R + 2] \times  T^{3}$. {Then, there is 
a gauge equivalent pair} $(A_{0} + b, \psi_{0} + \eta)$, 
{where} $(b,\eta)$ {defines a section of \eqref{eq:5.1} that} 
{obeys} 
\begin{equation}\begin{split}\label{eq:5.33}
(b,\eta) &= (u^{+}_{+} s^{+}_{+} + u^{+}_{-} s^{+}_{-}) 
\exp(-(2r)^{1/2} (t + R))\\
&+ (u^{-}_{+} s^{-}_{+} + u^{-}_{-} 
s^{-}_{-}) \exp(-(2r)^{1/2} (R - t)) + w_{1} 
\end{split}\end{equation}
{at all} $(t,\cdot) \in [-R+\zeta, R-\zeta] \times 
T^{3}$. {Here,} \[|w_{1}| \le \zeta \cdot \exp\left(-((2r)^{1/2} + 
\zeta^{-1})\cdot(R - |t|)\right),\] $\{u^{\pm}_{\pm}\}$ {are 
constant, and}
\begin{equation}\label{eq:5.34}
\sum_{0\le k\le3}|\nabla^{k}(b,\eta)| \le \zeta \cdot 
\exp(-(2r)^{1/2} (R - |t|)) 
\end{equation}
{at all points} $(t,\cdot)\in[-R, R] \times  T^{3}$.
\end{prop}

\begin{proof}[Proof of Proposition \ref{thm:5.6}]
According to Lemma \ref{thm:3.4}, there exists a pair\linebreak
$(A_{0} + b, \psi_{0} + \eta)$ on the gauge orbit of $(A,\psi)$ 
which obeys 
\begin{equation}\label{eq:5.35}
\sum_{0\le k\le3} |\nabla^{k}(b,\eta)| \le \zeta \cdot 
e^{-(R-|t|)/\zeta} 
\end{equation}
at all points $(t,\cdot)$ where $t \in [-R, R]$. Thus, on some smaller 
length cylinder, this pair differs very little from $(A_{0},\psi_{0})$ 
and so can be analyzed by treating \eqref{eq:2.1} 
as a perturbation of a linear 
equation. Moreover, as explained below, there is a fiducial choice of such a 
gauge equivalent pair $(A_{0} + b, \psi_{0} + \eta)$, where $(b, \eta)$ 
obeys \eqref{eq:5.34} and also
\begin{equation}\label{eq:5.36}
d^*b - 2 i \Ima(\psi_{0}^{\dagger}\eta) = 0 
\end{equation}
on a subcylinder of the form $[-R + \zeta_{1}, R - \zeta_{1}] 
\times  T^{3}$, where $\zeta_{1}$ depends only on the metric and $\omega 
$. 

With the preceding understood, view $w = (b, \eta)$ as a section of 
\eqref{eq:5.1}. 
By virtue of \eqref{eq:2.1} and 
\eqref{eq:5.26}, this section obeys an equation of the form in 
\eqref{eq:5.24},
where $\upsilon $ is a linear function of the components of $w$. In 
particular, Lemma \ref{thm:5.4} 
is applicable here with $R$ replaced by $R' \equiv R 
- \zeta_{1}$ because the condition in \eqref{eq:5.35} insures that $\upsilon $ 
obeys the requisite bounds on $[-R + \zeta_{1}, R - \zeta_{1}] 
\times  T^{3}$. Thus, \eqref{eq:5.33} can be seen to follow from 
\eqref{eq:5.25}. 

The refined derivative bounds for $(b,\eta)$ in \eqref{eq:5.34} 
are obtained from 
the $C^{0}$ bounds via standard elliptic techniques via \eqref{eq:5.24}. In 
particular, to obtain the $C^{1}$ bounds, simply differentiate \eqref{eq:5.24} 
and, 
remembering that $\upsilon $ is a linear functional of the components of $w$, 
observe that the result has the same schematic form as \eqref{eq:5.24}. 
Thus, a 
second appeal to Lemma \ref{thm:5.4} 
provides the $C^{1}$ bounds on $(b,\eta)$. The 
derivation of the $C^{3}$ bound is only slightly more complicated.

It remains now to justify the asserted existence of the point $(A_{0} + b, 
\psi_{0} + \eta)$ on the gauge orbit of $(A,\psi)$ for which both 
\eqref{eq:5.35} and \eqref{eq:5.36} 
hold. For this purpose, use Lemma \ref{thm:3.4} to conclude that $(A, 
\psi)$ is gauge equivalent to $(A_{0} + b', \psi_{0} + \eta')$, 
where $(b', \eta')$ obeys
\begin{equation}\label{eq:5.37}
\sum_{0\le k\le4} |\nabla^{k}(b',\eta')| \le \zeta'e^{-(R-|t|)/\zeta'}
\end{equation}
at all points $(t,\cdot)$ where $t \in [-R, R]$.
Here, $\zeta' > 0$ is a constant which is independent of $R$. The pair 
$(A_{0} + b, \psi_{0} + \eta)$ is guaranteed to come from the same 
gauge orbit as $(A_{0} + b', \psi_{0} + \eta')$ if $(b, \eta )$ and 
$(b', \eta')$ are related via the identity 
\begin{equation}\label{eq:5.38}
(b, \eta) = (b' - 2 i\cdot du, e^{i\cdot u}\cdot \eta' + 
(e^{i\cdot u} - 1)\cdot \psi_{0}), 
\end{equation}
where $u$ is a smooth function on $[-R - 2, R + 2] \times  T^{3}$. Thus, 
the goal is to find $u$ in \eqref{eq:5.38} so that \eqref{eq:5.36} 
holds on an appropriate 
subcylinder. In this regard, \eqref{eq:5.36}      
should be viewed as an equation for the 
function $u$. In particular, if $u$ 
has a suitably small $C^{2}$ norm, then this 
equation has the schematic form
\begin{equation}\label{eq:5.39}
- 2 i (d^*du + 2 r u) + d^*b' - 2 i 
\Ima(\psi_{0}^{\dagger}\eta') + 
\Re(u) = 0,                          
\end{equation}
where $|\Re| \le \zeta_{1} (|u|^{2} + |u| (|b'| + |\eta'|))$ with $\zeta 
_{1}$ a constant which is independent of $(b', \eta')$ and $R$. 

The analysis of \eqref{eq:5.39} 
is facilitated by the following observation: Because 
of \eqref{eq:5.37}, 
the pair $(b',\eta')$ is uniformly small on uniformly large 
subcylinders of $[-R, R] \times  T^{3}$. That is, given $\varepsilon  > 
0$, there exists $\xi  > 2$ which is independent of $(b', \eta')$ and $R$ 
such that at all points $(t,\cdot) \in [-R + \xi, R - \xi] 
\times  T^{3}$, all derivatives from orders $0$ through $4$ of $(b',\eta')$ 
are bounded in size by $\varepsilon e^{-(R-|t|)/\xi} $. Meanwhile, the 
Green's function, $G$, for $d^*d + 2 r$ 
defines a bounded operator from $L^{2}(\R 
\times  T^{3})$ to $L^{2}_{2}(\R \times T^{3})$ and satisfies 
the pointwise bound in \eqref{eq:5.12}. 

The preceding observations suggest a contraction mapping construction of a 
solution $u$ to \eqref{eq:5.39} 
on a uniformly large subcylinder of $[-R - 2, R + 2] 
\times  T^{3}$. For this purpose, fix $\xi  > 2$ and introduce a smooth, 
non-negative function $\beta $ on $\R$ 
which equals $1$ on $[-R + \xi, R - \xi]$, 
vanishes where $|t| > R - \xi  + 1$ and has first and second derivatives 
bounded by $10$. Then, consider the map from $C^{0}(\R \times  T^{3})$ to 
itself which sends $u$ to 
\begin{equation*}
T(u) \equiv - i 2^{-1} G(\beta[d^*b'- 2 i \Ima(\psi_{0}^{\dagger}\eta') 
+ \Re(u)]).
\end{equation*}
Here, the fact that $T$ defines a self map on $C^{0}(\R \times  T^{3})$ is 
insured by the right hand inequality in \eqref{eq:5.12}. 
Moreover, \eqref{eq:5.12} and \eqref{eq:5.37} 
imply the following: There exist constants $\xi \ge 2$ and $\xi' > 0$ which 
are independent of $R$ and $(b',\eta')$ and such that $T$ is a contraction 
mapping on the radius $\xi'$ ball in $C^{0}(\R \times  T^{3})$. For 
such $\xi $, the map $T$ has a unique fixed point, $u$, in this ball. 

\sloppy
Of course, \eqref{eq:5.12} 
insures that this $u$ decays to zero exponentially fast as 
$|t|\to\infty$ on  $\R \times  T^{3}$. Moreover, \eqref{eq:5.12} 
in conjunction with 
\eqref{eq:5.37} 
can be used to prove that $u$ and its derivatives obey the required 
norm bounds throughout in\linebreak$[-R -2, R + 2] \times  T^{3}$. These last 
derivations are straightforward and omitted.
\end{proof}\fussy

\sh{d) The distance to a non-trivial vortex}

Now, consider a half infinite tube of the form $ Y = [-2R, \infty)  \times  
T^{3}$, where $R \ge 4$. Fix $B_{0} \ge 0$ and $B_{1} \ge 1$ 
and suppose that $(A, 
\psi )$ is a solution to \eqref{eq:2.1} on $Y$ with the following properties: 
\begin{equation}\begin{split}\label{eq:5.41}
&\int_{Y}|F_{A}|^{2} \le B_{0};\\
&|F_{A}| \le B_{1}(\exp(t/B_{1}) + \exp(-(2R + 
t)/B_{1}))\text{ at points with }t \in  [-2R, 0].
\end{split}
\end{equation}

With Propositions \ref{thm:4.5} and \ref{thm:5.6} 
in mind, it can be expected that when $R$ is 
large, then $(A, \psi )$ is close on $[-R, \infty)  \times T^{3}$ to the 
restriction of a solution of \eqref{eq:2.1} 
on the whole of  $ \mathbb{R} \times  T^{3}$. 
Indeed, such is the case, as the subsequent proposition attests.

\begin{prop}\label{thm:5.7} 
{Given} $B_{0} \ge 0$ {and} $B_{1}\ge 
1$ {as above, there exists} $R_{0}\ge 4$ {and} $\zeta\ge 1$ 
{with the following significance: Suppose that} $R \ge R_{0}$ 
{and that} $(A, \psi)$ {obeys} \eqref{eq:2.1} {and} \eqref{eq:5.41} 
{on the half infinite tube} $Y = [-2R, \infty) \times  T^{3}$. 
{Then, there is a solution} $(A_{1},\psi_{1})$ {on} $\R 
\times  T^{3}$ to \eqref{eq:2.1} {and a gauge transform of} $(A,\psi)$ 
{on} $Y$ {such that} $(A,\psi) = (A_{1}, \psi_{1}) + w$ 
{on} $Y$, {where}
\begin{equation}\begin{split}\label{eq:5.42}
&|w| \le \zeta \exp(-(2r)^{1/2} R) \exp(-(2r)^{1/2} (R + t)) 
\text{ if }t \in  [-R, -R/2];\\
&\int_{t \ge -R/2} |w|^{2} \le \zeta \cdot \exp(-3 
(2r)^{1/2} R).
\end{split}
\end{equation}
\end{prop}

The remainder of this subsection contains the following:

\begin{proof}[Proof of Proposition \ref{thm:5.7}] 
Use Lemma \ref{thm:5.6} to find $\zeta_{0}\ge 1$ 
and a gauge for $(A,\psi)$ on the cylinder $[-2R + \zeta_{0}, -\zeta 
_{0}] \times  T^{3}$ of the form $(A,\psi) = (A_{0} + b, \psi 
_{0} + \eta)$, where $(b,\eta)$ come from Proposition \ref{thm:5.6} and obey 
\eqref{eq:5.33}, \eqref{eq:5.34} and \eqref{eq:5.36}, while 
$(A_{0},\psi_{0})$ is the pull-back 
from $T^{3}$ of a pair that defines $\mm_{P}$. Now, choose a smooth, 
non-negative function $\beta $ on $\R$ which equals $1$ where $t \ge 1$ 
and $0$ where $t\le0$. 
Use $(b,\eta)$ and the function $\beta $ to define a configuration 
$(A',\psi')$ on $\R\times  T^{3}$ as follows:
\begin{equation*}\begin{split}
&(A',\psi') = (A_{0},\psi_{0})\text{ where }t < -R - 2;\\
&(A',\psi') = (A_{0},\psi_{0}) + \beta(t+R) (b, 
\eta)\text{ where }t \in [-R-2, -R+2];\\
&(A',\psi') = (A,\psi)\text{ where }t \ge -R + 2. 
\end{split}
\end{equation*}

Note that $\hh = (P_{+}F_{A'} - \tau(\psi'\otimes(\psi')^{\dagger}) + i\cdot 
\omega, D_{A'}\psi')$ vanishes except when the coordinate $t \in 
[-R + 1, -R + 2]$. Moreover, when $t \in [-R + 1, -R +2]$, then 
\eqref{eq:5.34} 
guarantees that 
\begin{equation}\label{eq:5.44}
\sum_{0\le k\le2} |\nabla^{k}\hh| \le \zeta_{1} \exp(-(2r)^{1/2} R), 
\end{equation}
where $\zeta_{1}$ is independent of $R$; it depends only on the constants 
$B_{0}$ and $B_{1}$. Thus, for large $R$, the pair $(A',\psi')$ is very 
close to solving \eqref{eq:2.1} on $\R \times  T^{3}$. 
The following lemma makes 
this notion precise:

\begin{lem}\label{thm:5.8} 
{Under the assumptions of Proposition \ref{thm:5.7}, there 
exists} $m \ge 0$, $\varepsilon_{0} > 0$ {and, given} $\varepsilon   
\in (0, \varepsilon_{0})$, {there exists} $R_{\varepsilon} $ 
and these have the following properties: 
First, an upper bound for 
$m$, the numbers $\varepsilon_{0}$, and a lower bound 
for $R_{\varepsilon}$ depend only on $B_{0}$ and $B_{1}$.
{Second, when} $R \ge R_{\varepsilon}$, {then the pair} $(A', 
\psi')$ {has} $C^{2}$ {distance less than} $\varepsilon $ 
{from gauge orbits of solutions to} \eqref{eq:2.1} {on} $\R \times 
T^{3}$ that come from vortex solutions with vortex number 
$m$. Third, such an orbit contains a unique pair $(A_{1}, \psi 
_{1})$ {for which} $(b_{1},\eta_{1}) \equiv (A' - 
A_{1}, \psi' - \psi_{1})$ {obeys}
\begin{equation}\begin{split}\label{eq:5.45}
\bullet\text{ }&d^*b_{1} - 2\cdot i\cdot \Ima(\psi_{1}^{\dagger}\eta 
_{1}) = 0;\\
\bullet\text{ }&\int(|b_{1}|^{2} + |\eta_{1}|^{2}) \le \xi  
\varepsilon^{2};\\
\bullet\text{ }&\sum_{0\le k\le2} |\nabla^{k}(b_{1},\eta_{1})| \le \xi 
\varepsilon\text{ everywhere};\\
\bullet\text{ }&|(b_{1},\eta_{1})| \le \xi\varepsilon 
\left(\exp(-(2r)^{1/2}(2R - \xi + t)) + \exp((2r)^{1/2}(t + \xi))\right),\\
&\text{where }t \in [-2 R + \xi, -\xi].
\end{split}
\end{equation}
{Here,} $\xi ${ depends only on the vortex number} $m$; 
{in particular, it is independent of} $\varepsilon $, $R$, {and 
the original pair} $(A,\psi)$. 
\end{lem}
The proof of Lemma \ref{thm:5.8} is given below. 

Lemma \ref{thm:5.8} 
enters the proof of Proposition \ref{thm:5.7} in the following way: Fix some 
positive $\varepsilon\ll\varepsilon_{0}$ and then $R_{\varepsilon} $ 
as in Lemma \ref{thm:5.8}. An upper bound for $\varepsilon $ is derived in the 
subsequent arguments. In any event, suppose that $R > R_{\varepsilon} $. Let 
$\uu \subset \ccc_{m}$ denote the open set of elements which provide gauge 
orbits of solutions to \eqref{eq:2.1} on $\R\times T^{3}$ that have $C^{2}$ 
distance less than $\varepsilon $ from the pair $(A',\psi')$. Given a 
vortex solution from $\uu$, let $(A_{1},\psi_{1})$ denote the 
corresponding solution to \eqref{eq:2.1} 
that is provided by Lemma \ref{thm:5.8}. Introduce the 
resulting $(b_{1},\eta_{1})$ and observe that $w_{1} = (b_{1}, 
\eta_{1})$ obeys an equation of the form
\begin{equation}\label{eq:5.46}
\dd_{1}w_{1} + \qq(w_{1}) + \hh = 0.
\end{equation}
Here, $\dd_{1}$ denotes the operator $\dd_{c}$ from 
\eqref{eq:5.3} as defined using $c = 
(A_{1},\psi_{1})$, $\qq(\cdot)$ is a universal, quadratic, fiber 
preserving map from \eqref{eq:5.1} to \eqref{eq:5.2} 
and $\hh$ is interpreted as a section of 
\eqref{eq:5.2}. 

With $w_{1}$ and $\dd_{1}$ understood, consider:

\begin{lem}\label{thm:5.9} 
{Under the assumptions of Proposition \ref{thm:5.7}, there 
exists} $\varepsilon_{1}  > 0$ {which depends only on} $B_{0}$ 
{and} $B_{1}$ and which has the following significance: If 
$\varepsilon < \varepsilon_{1}$, {then there exists} $(A_{1}, 
\psi_{1})$ {as described in the preceding paragraph for which the 
corresponding} $w_{1}$ {is} $L^{2}$--orthogonal to the kernel 
of the operator $\dd_{1}$.
\end{lem}
The proof of this lemma is also given below.

Given the statement of Lemma \ref{thm:5.9}, 
the proof of Proposition \ref{thm:5.7} continues 
with the following observation: There exists $\varepsilon_{2} > 0$ and 
$\zeta_{2}$ which depend only on $B_{0}$ and $B_{1}$ and are such that 
when $\varepsilon < \varepsilon_{2}$, then the $L^{2}$ norm of Lemma 
\ref{thm:5.9}'s section $w_{1}$ satisfies
\begin{equation}\label{eq:5.47}
\| w_{1} \|_{2} \le \zeta_{2} \exp(-(2r)^{1/2} R).
\end{equation}
Indeed, the existence of such a pair $(\varepsilon_{2}, \zeta_{2})$ 
is guaranteed by \eqref{eq:5.17}, \eqref{eq:5.44} 
and the third point in \eqref{eq:5.45}. As is 
explained below, this upper bound on the $L^{2}$ norm of $w_{1}$ implies the 
pointwise bounds:
\begin{equation}\begin{split}\label{eq:5.48}
\bullet\text{ }&|w_{1}| \le \zeta_{3} \exp(-(2r)^{1/2} R) 
\text{ everywhere};\\
\bullet\text{ }&|w_{1}| \le \zeta_{3} \exp(-(2r)^{1/2} R) 
\left[\exp(-(2r)^{1/2} (R + t)) + \exp((2r)^{1/2}(t + \zeta_{3}))\right]\\ 
&\text{when }t \in [-R, -\zeta_{3}].
\end{split}
\end{equation}
Here, $\zeta_{3}$ depends only on $B_{0}$ and $B_{1}$; in particular, it 
is independent of $R$. Indeed, the first line in \eqref{eq:5.48} 
follows from \eqref{eq:5.44}, 
\eqref{eq:5.46} and \eqref{eq:5.47} 
using standard elliptic estimates, while the second 
follows from the first after an appeal to Lemma \ref{thm:5.4}.

Notice that the first point of \eqref{eq:5.48} 
implies the first point in \eqref{eq:5.42}. The 
second point in \eqref{eq:5.42} 
is obtained as follows: Start with the second line of 
\eqref{eq:5.48} and so conclude that
\begin{equation}\label{eq:5.49}
|w_{1}| \le \zeta_{4} \exp(-3 (2r)^{1/2}R/2)
\text{ where }t \in 
\left[-\frac{R}2 -2, -\frac{R}2 + 2\right]. 
\end{equation}
Here, $\zeta_{4}$ is independent of $R$ and is determined solely by $B_{0}$ 
and $B_{1}$. This last bound is used to derive an upper bound for the 
$L^{2}$ norm of the section $w_{2}$ of \eqref{eq:5.1} 
that is defined by the rule
\begin{itemize}
\item  $w_{2} = 0$ where $t \le -R/2$,
\item  $w_{2} = \beta(-\frac{R}2+t)w_{1}$ where $t \ge -R/2$.
\end{itemize}
Note that a suitable upper bound for the $L^{2}$ norm of $w_{2}$ will yield, 
with \eqref{eq:5.49}, the second point in \eqref{eq:5.42}. 

To obtain such a bound, use \eqref{eq:5.46} to conclude that $w_{2}$ obeys an 
equation of the form
\begin{equation}\label{eq:5.51}
\dd_{1}w_{2} + \qq(w_{2}) + \hh_{2} = 0,
\end{equation}
where $\hh_{2} = 0$ except when $t \in [-\frac{R}2, -\frac{R}2 + 2]$. 
Moreover, where 
$\hh_{2}$ is not zero, it obeys the bound $|\hh_{2}| \le \zeta_{4} \exp(-3 
(2r)^{1/2} R/2)$ by virtue of \eqref{eq:5.49}; here $\zeta_{4}$ 
is independent of 
$R$ as it is determined solely by $B_{0}$ and $B_{1}$. 

With \eqref{eq:5.51} 
understood, write $w_{2} = w_{20} + w_{21}$, where $w_{20}$ 
is the $L^{2}$--orthogonal projection of $w_{2}$ onto the kernel of $\dd_{1}$. 
In this regard, $w_{20}$ enjoys the following upper bound: 
\begin{equation}\label{eq:5.52}
|w_{20}| + \| w_{20} \|_{2} \le \zeta_{5} \exp(- 3 (2r)^{1/2} R/2),
\end{equation}
where $\zeta_{5}$ is independent of $R$, being determined solely by $B_{0}$ 
and $B_{1}$. Indeed, remember that $w_{1}$ is orthogonal to the kernel of 
$\dd_{1}$ and so the projection of $w_{2}$ onto this kernel is the same as 
that of $[1 - \beta(-\frac{R}2+(\cdot))] w_{1}$. 
In particular, the size of 
the latter projection obeys \eqref{eq:5.52}        
for the following reasons: First, $w_{1}$ 
enjoys the bound in \eqref{eq:5.49}. 
Second, Lemma \ref{thm:5.4} guarantees that any $\varpi 
\in \kernel(\dd_{1})$ is bounded by $\zeta_{6}\|\varpi\|_{2} 
\exp((2r)^{1/2}t)$ at the points where $t \le -\zeta_{6}$. Finally, any 
$\varpi \in \kernel(\dd_{1})$ obeys $|\varpi| \le \zeta_{6} \| 
\varpi \|_{2}$ at all points. Here, again, $\zeta_{6}$ is independent 
of $R$; it depends only on the vortex number $m$ and hence only on $B_{0}$ and 
$B_{1}$.

Given \eqref{eq:5.52}, 
the previously mentioned bound on $\hh_{2}$, and the fact that 
\eqref{eq:5.51} can be written as $\dd_{1}w_{21} + \qq(w_{21}) + \qq'(w_{20}, 
w_{21}) + \qq(w_{20}) + \hh_{2} = 0$, another appeal to \eqref{eq:5.17} finds 
$\varepsilon_{3} > 0$ and $\zeta_{7}$ which are independent of $R$, 
depend only on $B_{0}$ and $B_{1}$ and are such that when $\varepsilon < 
\varepsilon_{3}$, then $\| w_{21} \|_{2} \le \zeta_{7} \exp(-3 
(2r)^{1/2} R/2)$. This last bound, \eqref{eq:5.49} and \eqref{eq:5.52} 
directly yield the 
final point in \eqref{eq:5.42}. 
\end{proof}

\begin{proof}[Proof of Lemma \ref{thm:5.8}] 
All of the arguments for Lemma \ref{thm:5.8} closely 
follow arguments previously given, and so, except for the outline that 
follows, the details are left to the reader. The outline starts with the 
observation that a slightly modified version of the argument for Proposition 
\ref{thm:4.5} 
establishes the existence of an upper bound for $m$ and $R_{\varepsilon} $ 
such that when $R > R_{\varepsilon} $, then $(A',\psi')$ has $C^{2}$ 
distance $\varepsilon $ or less from a solution to \eqref{eq:2.1} 
on $\R \times 
T^{3}$. The bound on the vortex number comes from the $L^{2}$ constraint in 
the statement of Proposition \ref{thm:5.7}. 
Given that $(A',\psi')$ is close to a 
solution to \eqref{eq:2.1}, 
then the points in \eqref{eq:5.45} are proved by arguments which 
are essentially the same as those used above to prove Proposition 
\ref{thm:5.6}. Note 
that these may require an increase of the lower bound for $R$. The arguments 
for the lemma's uniqueness assertion are basically those used to prove the 
slice theorem for the action of the gauge group on the space of solutions to 
\eqref{eq:2.1}.
\end{proof}

\begin{proof}[Proof of Lemma \ref{thm:5.9}]
As explained previously, Lemma \ref{thm:5.8} provides the non\-empty,
open set $\uu \subset \ccc_{m}$ of elements that determine solutions
to \eqref{eq:2.1} on $\R \times T^{3}$ whose $C^{2}$ distance is less
than $\varepsilon $ from $(A',\psi')$. The assignment to a point in $\uu$ of 
the corresponding $w_{1}$ defines a map from $\uu$ into the space of $L^{2}$ 
sections of \eqref{eq:5.1}. 
This map is smooth; the proof is straightforward so its 
details are left to the reader. The assignment to a point in $\uu$ 
of the square 
of the $L^{2}_{1}$ norm of the corresponding $w_{1}$ then defines a smooth 
function, $f$, on $\uu$. As is explained momentarily, the differential of $f$ 
vanishes at precisely the points where $w_{1}$ is $L^{2}$--orthogonal to the 
kernel of $\dd_{1}$. Indeed, let $b$ denote a tangent vector to $\uu$ 
at the point 
that corresponds to $(A_{1},\psi_{1})$. Then, the differential of 
$w_{1}$ in the direction defined by $b$ has the form $w_{1b} + \delta 
_{b}$, where $w_{1b}\in \kernel(\dd_{1})$ and $\delta_{b}$ is 
tangent to the orbit through $(A_{1},\psi_{1})$ of the gauge group 
$C^{\infty}(\R\times T^{3}, S^{1})$. This is a consequence of the fact 
that $(A',\psi')$ is fixed and only $(A_{1},\psi_{1})$ is moved by 
$b$. Meanwhile, $w_{1}$ is $L^{2}$--orthogonal to $\delta_{b}$ by virtue of 
the first point in \eqref{eq:5.45} 
and so the differential of $f$ vanishes in the 
direction of $b$ if and only if $w_{1}$ is $L^{2}$--orthogonal to $w_{1b}$. As 
Proposition \ref{thm:4.4} 
guarantees that the kernel of $\dd_{1}$ is the span of the 
$\{w_{1b}\}$, so $w_{1}$ is orthogonal to $\kernel(\dd_{1})$ if and only if 
$(A_{1},\psi_{1})$ comes from a critical point of $f$.

With the preceding understood, it remains only to demonstrate that $f$ has 
critical points. For this purpose, a return to \eqref{eq:5.46} is in order. In 
particular, with \eqref{eq:5.17}, \eqref{eq:5.34} 
and in conjunction with standard elliptic 
regularity arguments, \eqref{eq:5.46} 
leads to the following conclusion: There exist 
constants $\xi\ge 1$ and $\delta  > 0$ which are independent of $\varepsilon 
$ and $R$ and such that if $|w_{1}| 
< \delta $, then the $C^{1}$ norm of 
$w_{1}$ is bounded by 
$\xi (\| w_{1} \|_{2} + 
\exp(-R/\xi))$. Since the 
$C^{0}$ norm of $w_{1}$ is bounded by its $C^{1}$ norm, this last point 
and the second point in \eqref{eq:5.45} 
guarantee that $f$ is proper when $\varepsilon 
$ is small and $R$ is large. Said precisely, $\varepsilon_{1} > 0$ and 
$R_{1} \ge 1$ must exist such that when $\varepsilon < \varepsilon_{1}$ 
and $R > R_{1}$, then the map $f$ is a proper map on the set $\uu$. 
As a proper 
function has at least one critical point so there exists at least one 
$(A_{1},\psi_{1})$ where the corresponding $w_{1}$ and $\kernel(\dd_{1})$ 
are orthogonal. 
\end{proof}

\section{Compactness}\label{sec:6}

This last section contains the final arguments for Propositions 
\ref{thm:2.4}, \ref{thm:3.7} and 
\ref{thm:3.9}. 

\sh{a) Proof of Proposition \ref{thm:2.4}}

The assertion that $\mm_{s,m}$ is both compact and consists of a finite number 
of strata follows from Proposition \ref{thm:4.5}. 
Meanwhile, the assertion in the 
first point of the proposition is a restatement of the conclusions of 
Proposition \ref{thm:2.3}. 
Thus, the only remaining issue is that of the second point. 
The latter's assertion is proved in the subsequent seven steps. 

\textbf {Step 1}\qua To begin, consider some $z \in 
\varsigma({s})$ and a sequence $\{c_{j}\} \subset \mm({s}, 
z)$ with no convergent subsequences. After passage to a subsequence, as 
always renumbered consecutively from $1$, the sequence $\{c_{j}\}$ can be 
assumed to determine data $c_{\infty}$ 
and, for each component of $\partial X_{0}$, a 
sequence $\{o_{j}\}$ of solutions on $\R \times  T^{3}$ to \eqref{eq:2.1}, 
all as 
described in Proposition \ref{thm:4.5}. 
In this regard, remember that each $o_{j}$ is 
determined by a solution, $(\tau_{j}, v_{j})$, of the vortex equations 
in \eqref{eq:4.14}. 

Fix a component $[0, \infty)  \times  T^{3}$ of 
$[0, \infty)  \times  \partial X_{0}$, let 
$\{o_{j}\}$ denote the corresponding sequence of solutions on $\R \times 
T^{3}$, and for each $j$, let $t_{j}$ denote the smallest of the time 
coordinates of the zeros of the corresponding $\tau_{j}$. These 
$\{t_{j}\}$ can be assumed to define an increasing and unbounded sequence.

Propositions \ref{thm:5.6} and \ref{thm:5.7} 
determine certain data with certain special 
properties. Here is the data: A constant $\zeta \ge 1$; a pair 
$(A_{\infty}, 
\psi_{\infty})$ on the gauge orbit $c_{\infty}$ over $[0, \infty)  \times 
T^{3}$; for each sufficiently large index $j$, a pair $(A_{j}, \psi 
_{j})$ on the gauge orbit $c_{j}$ over $[0, \infty)  \times  T^{3}$; and, for 
each such $j$, a pair $(A_{1j},\psi_{1j})$ on $o_{j}$ over $\R \times 
T^{3}$. This data enjoys the special properties that are listed 
\eqref{eq:6.1}--\eqref{eq:6.3} 
below. Note that in these equations, $(A_{0},\psi_{0})$ denotes the 
solution to \eqref{eq:2.1} on $\R\times T^{3}$ which is given by the vortex 
equation solution with vortex number zero. Moreover, the bundle $S_{+}$ is 
implicitly written as in \eqref{eq:5.1}, $S_{+} = \varepsilon_{\C} \oplus 
\varepsilon_{01}$ and $\psi_{0}$ is used to trivialize the 
$\varepsilon_{\C}$ summand and thus defines the section with vanishing 
imaginary part and positive real part. Finally, the sections $s^{+}_{\pm} $ 
of \eqref{eq:5.1} are defined in \eqref{eq:5.21}.

Here is the promised list of properties: 
{\sl
\begin{equation}\begin{split}\label{eq:6.1}
\bullet\text{ }&(A_{\infty},\psi_{\infty}) = (A_{0},\psi_{0}) + 
(u^{+}_{+} s^{+}_{+} + u^{+}_{-} s^{+}_{-}) 
\exp(-(2r)^{1/2} (t - \zeta)) + w_{\infty},\\
&\text{where }t \ge \zeta.
\text{ Here, }|w_{\infty}| \le \zeta \cdot \exp(-z t) 
\text{ and }z \ge (2r)^{1/2} + \zeta^{ - 1}.
\end{split}
\end{equation}} 
This is by virtue of Proposition \ref{thm:5.6}.
{\sl
\begin{equation}\begin{split}\label{eq:6.2}
\bullet\text{ }&(A_{j},\psi_{j}) = (A_{0},\psi_{0}) + 
(u_{j+}^{+} s^{+}_{+} + u_{j-}^{+} s^{+}_{-}) 
\exp(-(2r)^{1/2} (t - \zeta)) + w_{j},\\
&\text{where }t \in 
[\zeta, t_{j} - \zeta].\\ 
&\text{Here, }|w_{j}| \le \zeta [\exp(-z t) + \exp(-z (t_{j} 
- t))]\text{ and }z \ge (2r)^{1/2} + \zeta^{-1}.\\
\bullet\text{ }&|u_{j+}^{+} - u^{+}_{+}| + |u_{j-}^{+} - 
u^{+}_{-}| \to 0\text{ as }j \to \infty.\\
\bullet\text{ }&t_{j}\to\infty\text{ as }j \to \infty.
\end{split}
\end{equation}}
The points in \eqref{eq:6.2} also follow from Proposition \ref{thm:5.6}.
{\sl
\begin{equation}\begin{split}\label{eq:6.3}
\bullet\text{ }&(A_{j},\psi_{j}) = (A_{1j},\psi_{1j}) + w_{1j} 
\text{ where }t \ge t_{j}/2.\\
&\text{Here, }|w_{1j}| \le \zeta \exp(-(2r)^{1/2} 
t_{j}/2) \exp(-(2r)^{1/2}(t - t_{j}/2)),\\
&\text{where }t \in 
[t_{j}/2, 3t_{j}/4].\\
\bullet\text{ }&\int_{t \ge 3t_{j} /4} |w_{j}|^{2}  \le \zeta  \exp(-3 
(2r)^{1/2} t_{j}/2).\\
\bullet\text{ }&w_{j}\text{ satisfies an equation of the form }
\dd_{j} w_{1j} + \qq(w_{1j}) = 0,
\text{ where}\\
&t \ge t_{j}/2.\text{ Here, }w_{1j}\text{ is viewed as 
a section of \eqref{eq:5.1}, }\dd_{j}\equiv\dd_{c}\text{ with}\\
&c=o_{j}\text{ and }\qq(\cdot)\text{ is a universal, quadratic, fiber 
preserving map}\\
&\text{from \eqref{eq:5.1} to \eqref{eq:5.2}.}
\end{split}\end{equation}}

The points in \eqref{eq:6.3} follow from Proposition 
\ref{thm:5.7} after translating the origin 
so that the new origin corresponds to $t_{j}$ and hence $-t_{j}$ corresponds 
to the old $0$. Then, take $R = t_{j}/2$. Note that $t_{j}$ is the most 
negative of the $t$--coordinates 
of the zeros of the $o_{j}$ version $\tau $ in 
\eqref{eq:4.15}.

\textbf{Step 2}\qua The solution $(A_{1j},\psi_{1j})$ on 
$\R\times T^{3}$ to \eqref{eq:2.1} is defined by a solution $(\tau_{j}, 
v_{j})$ to the vortex equation in \eqref{eq:4.14} 
on $\R \times  S^{1}$. As long 
as the latter is not the trivial vortex (gauge equivalent to the pair $(0, 
1)$), then the operator $\Theta $ in \eqref{eq:5.4} 
has complex multiples of the 
element in \eqref{eq:5.14} 
in its kernel. To underscore the index dependence, use 
$\pi_{j}$ to denote this element. Then $\dd_{j}$ 
has the element $\pi_{j-} = (0, 
\pi_{j})$ in its cokernel. Here, $\pi_{j-}$ 
is defined as in Lemma \ref{thm:5.4} using $c 
= o_{j}$.

Note that \eqref{eq:5.31} and Lemma \ref{thm:5.4} 
assert that at $t \in [\frac{t_{j}}2 - 2, 
\frac{t_{j}}2 + 2]$, this $\pi_{j-}$ has the form
\begin{equation}\label{eq:6.4}
\pi_{j-} = \exp(-(2r)^{1/2} t_{j}/2) h_{j} s^{+}_{-} + \sigma_{j}'.
\end{equation}
Here, $h_{j}$ is a positive constant, $h_{j} \ge \xi^{-1}$, while 
$|\sigma_{j}'| \le \xi \cdot \exp(-z\cdot t_{j}/2)$ with $z > 
(2r)^{1/2} + \xi^{-1}$. In this regard, $\xi\ge1$ is a constant which 
is independent of the index $j$. (The vector $s^{+}_{-}$ which appears in 
\eqref{eq:6.4} is given in \eqref{eq:5.21}.)

\textbf {Step 3}\qua Take the inner product of both sides of 
the equality $\dd_{j} w_{1j} + \qq(w_{1j}) = 0$ with $\pi_{j-}$ and then 
integrate over the region where $t \ge t_{j}/2$. With \eqref{eq:5.20} 
in mind, apply 
integration by parts to the resulting expression and so find that
\begin{equation}\label{eq:6.5}
\int_{t \ge t_{j} /2}              
\langle \pi_{j-}, 
w_{1j}\rangle = 2 \int_{t 
\ge t_{j}/2} \langle \pi_{j-}, \qq(w_{1j}) \rangle.
\end{equation}
It follows from \eqref{eq:6.2} that the left hand side of \eqref{eq:6.4} 
is equal to
\begin{equation}\label{eq:6.6}
h_{j} u_{j-}^{+} \exp(-(2r)^{1/2} t_{j}) + \xi \exp(- z t_{j}),
\end{equation}
where $z > (2r)^{1/2}$ and $\xi $ are independent of $j$.  Furthermore, for 
large $j$, the number $u_{j-}^{+}$ is determined also by $c_{\infty}$ up to 
a small error because of the last point in \eqref{eq:6.2}. 
Indeed, \eqref{eq:6.6} can be 
rewritten as
\begin{equation}\label{eq:6.7}
h_{j} (u^{+}_{-} + \varepsilon_{j}) \exp(-(2r)^{1/2} t_{j}), 
\end{equation}
where $|\varepsilon_{j}| \to0$ as 
$j \to \infty$ and $u^{+}_{-}$        
is given in \eqref{eq:6.1}.

Meanwhile, the second and third points of \eqref{eq:6.3} and \eqref{eq:5.23} 
imply that the 
absolute value of the right hand side of \eqref{eq:6.5} is no greater than
\begin{equation}\label{eq:6.8}
\xi \exp(-3 (2r)^{1/2} t_{j}/2),
\end{equation}
where $\xi $ is, once again, independent of the index $j$.

\textbf {Step 4}\qua It follows from \eqref{eq:6.7}, \eqref{eq:6.8} 
and the fact 
that $h_{j} > \zeta ^{-1}$ that the sequence $\{c_{j}\}$ can exist as 
described only if the configuration $c_{\infty}$ is such that the number 
$u^{+}_{-}$, which arises in \eqref{eq:6.1}, is zero.

Associate to each element $c  \in \mm_{s}$ 
and each component of $\partial X_{0}$ 
the complex number $u^{+}_{-}$ using \eqref{eq:6.1}. 
The proof of Proposition \ref{thm:2.4} 
is completed with a demonstration of the fact that for a suitably generic 
choice of $\omega $ in \eqref{eq:2.1}, 
the set of those $c$ in $\mm_{s}$ where any given 
$u^{+}_{-}$ is zero has codimension at least $2$. 

To see that such is the case, first fix a fiducial choice $\omega_{0}$ of 
self-dual $2$--form which is non-zero and covariantly constant on each end of 
$X$. Next, fix an open set $K \subset  X_{0}$ with compact closure in the 
interior of $X_{0}$, and consider the space $\underline{\mm}_{s}$ that is 
defined as follows: As a set, this space consists of pairs $(c, \omega)$ 
where $\omega $ is a self-dual $2$--form on $X_{0}$ that agrees with $\omega 
_{0}$ on the complement of $K$, and where $c$ is the orbit under the action of 
$C^{\infty}(X;S^{1})$ of a solution to the version of \eqref{eq:2.1} 
which is defined 
by the given form $\omega $. The topology on $\underline{\mm}_{s}$ is the 
minimal topology so that the assignment to a pair $(c,\omega)\in 
\underline{\mm}_{s}$ of the form $\omega - \omega_{0}$ provides a 
continuous map to the space of smooth, compactly supported sections over $K$ 
of $\Lambda_{+}$ whose fibers are topologized as before. The strata of 
$\underline{\mm}_{s}$ are labeled by the elements of $\varsigma({s})$ 
and are defined in the obvious way. They are all smooth manifolds for which 
the assignment to $(c,\omega)$ of $\omega - \omega_{0}$ is a smooth 
map into the Fr\'echet space of compactly supported, smooth sections over $K$ of 
$\Lambda_{+}$. (Given Proposition \ref{thm:2.2}, 
these last assertions about the 
strata of $\underline{\mm}_{s}$ are proved with virtually the same 
arguments that establish the analogous assertion for compact $4$--manifolds.) 

Fix attention on a stratum, $\underline{\mm}({s},z) \subset 
\underline{\mm}_{s}$. 
Given a component of $\partial X_{0}$, define a map, $\varphi 
\colon \underline{\mm}({s}, z) \to \C$, by associating the complex 
number $u^{+}_{-}$ from \eqref{eq:6.1} to each $(c,\omega)$. An argument is 
provided below for the following assertion:

\begin{lem}\label{thm:6.1} 
{The map} $\varphi $ {is a smooth map 
and it has no critical points}.
\end{lem}

Given Lemma \ref{thm:6.1}, 
an appeal to the Sard--Smale theorem finds the Baire subset 
which makes the second point of Proposition \ref{thm:2.4} true.

With regard to this Sard--Smale appeal, remember that the latter considers 
maps between Banach manifolds and so an extra step is required for its use 
here. This extra step requires the introduction of a sequence of Banach 
space versions of $\underline{\mm}_{s}$ whose intersection is the smooth 
version given above. However, the use of Banach spaces modeled on $C^{p}$ 
for $p \ge 2$ or $L^{2}_{k}$ for $k \gg 1$ 
makes no essential difference in any 
of the previous or subsequent arguments. 

\textbf {Step 5}\qua This step, and the subsequent steps 
contain the following:

\begin{proof}[Proof of Lemma \ref{thm:6.1}] 
The proof that $\varphi $ is smooth is 
straightforward and will be omitted. The subsequent discussion addresses the 
question of whether $\varphi$ has critical points. For this purpose, it 
proves convenient to identify the tangent space to a given stratum 
$\underline{\mm}$ of $\underline{\mm}_{s}$ 
at a point of interest, $(c = (A, \psi),\omega)$,    
with a Fr\'echet space of smooth, square integrable 
pairs $(w = (a, \eta), \kappa)$                     
that obey the equation $\dd_{c} w + 
\kappa = 0$. Here, $a$ is an imaginary valued section of $T^*X$, $\eta $ is a 
section of $S_{+}$, and $\kappa $ is an imaginary valued section of $\Lambda 
_{+}$ with compact support on $K$. It is a consequence of Lemma \ref{thm:5.4} 
and 
Proposition \ref{thm:5.6} 
that the restriction of $w$ to each component of $[0, \infty)  
\times  \partial X_{0}$ has the form
\begin{equation}\label{eq:6.10}
w = \exp(-(2r)^{1/2} t) (w^{+}_{+} s^{+}_{+} + w^{+}_{-} 
s^{+}_{-}) + w',
\end{equation}
where $w^{+}_{\pm} $ are complex numbers and $|w'| \le \zeta e^{-zt}$ 
with $z$ a universal constant which is greater than $(2r)^{1/2}$. 

With this last equation understood, then the statement of Lemma 
\ref{thm:6.1} follows 
with the verification that any value for the complex number $w^{+}_{-}$ 
can be obtained by a suitable choice of pairs $(w,\kappa)$ which obey 
$\dd_{c}w + \kappa = 0$. This verification is the next order of business.

\textbf {Step 6}\qua The constant $w^{+}_{-}$ in \eqref{eq:6.10} is 
given by
\begin{equation}\label{eq:6.11}
w^{+}_{-} = (2r)^{1/2} \lim_{R \to \infty} \exp(2(2r)^{1/2} R) \int 
_{t>R} (s^{+}_{-})^{\dagger}w \exp(-(2r)^{1/2} t).
\end{equation}
With \eqref{eq:6.11} 
understood, let $\alpha $ be a non-zero complex number and let 
$x_{R}$ denote the section of $i\cdot T^*X \oplus S_{+}$ which is zero except 
where $t \ge R$ 
on the given component of $[0, \infty)  \times  \partial X_{0}$ in which 
case 
\begin{equation}\label{eq:6.12}
x_{R} \equiv (2r)^{1/2} \alpha 
s^{+}_{-} \exp((2r)^{1/2} (2R - t)).       
\end{equation}
Thus,
\begin{equation}\label{eq:6.13}
\Rea(\bar{\alpha} w^{+}_{-}) =              
\lim_{R \to \infty} \langle x_{R}, w\rangle,
\end{equation}
where $\langle\;,\;\rangle$ denotes the $L^{2}$ inner product over $X$ and 
$\Rea$ denotes the real part. This 
section $x_{R}$ is introduced for reasons which should be clear momentarily. 

Now, as $\dd_{c} w + \kappa = 0$, the section $w$ can also be written as 
\begin{equation*}
w = w_{0} - \dd_{c}^{-1} \kappa, 
\end{equation*}
where $\dd_{c}w_{0} = 0$ holds and $\dd_{c}^{-1}$ maps the $L^{2}$--orthogonal 
complement in\linebreak 
$L^{2}(i\cdot (\mathbb{R} \oplus \Lambda_{+}) \oplus S_{-})$ of 
$\cokernel(\dd_{c})$ 
to the $L^{2}$--orthogonal complement of $\kernel(\dd_{c})$ in 
$L^{2}_{1}(i T^*X \oplus S_{+})$. With the help of this decomposition, 
\eqref{eq:6.13} 
becomes
\begin{equation}\label{eq:6.15}
\Rea(\bar{\alpha}w^{+}_{-}) = - \lim_{R \to \infty} (\langle x_{R}, 
\dd_{c}^{-1}\kappa \rangle - \langle x_{R}, w_{0}\rangle).
\end{equation}

To proceed, introduce $x^{0}_{R}$ to denote the $L^{2}$--orthogonal 
projection of $x_{R}$ onto the kernel of $\dd_{c}$ and then introduce $y_{R} 
\equiv (\dd_{c}^{\dagger})^{-1} (x_{R} - x^{0}_{R})$. Equation 
\eqref{eq:6.15} can be rewritten with the help of $y_{R}$ as
\begin{equation}\label{eq:6.16}
\Rea(\bar{\alpha} w^{+}_{-}) = - \lim_{R \to \infty} (\langle y_{R}, 
\kappa \rangle - \langle x^{0}_{R}, w_{0}\rangle).
\end{equation}
In the meantime, it follows from Lemma \ref{thm:5.4} 
(after changing $t$ to $-t$) that 
there exists $\zeta  > 1$ which is independent of $R$ and such that $|y_{R}| 
< \zeta $ on $X_{0}$. Moreover, there exists $L \ge 0$ which is independent of 
$R$ and such that
\begin{equation*}
|y_{R}| \ge \zeta^{-1} \exp((2r)^{1/2} t), 
\end{equation*}
where $t\in[L, R]$ on the given component of $[0,\infty)\times\partial X_{0}$. 
In this regard, note that Lemma \ref{thm:5.4} 
insures that $|x^{0}_{R}|$ enjoys an $R$ 
and $\alpha $ independent bound on $X$ when $|\alpha| \le 1$. 

Now, there are two possibilities to consider. The first is that there exists 
a complex number $\alpha\ne0$ such that 
\begin{equation}\label{eq:6.18}
\lim_{R \to 0} \sup_{X} |x^{0}_{R}| = 0.
\end{equation}
The second possibility is that there is an unbounded subsequence of values 
for $R$ such that the corresponding sequence $\{\sup|x^{0}_{R}|\}$ has a 
non-zero limit for each unit length $\alpha \in \C$. Now, in this second case, 
it follows from \eqref{eq:6.16} that all complex numbers can be realized by 
$w^{+}_{-}$ in \eqref{eq:6.10} 
by tangent vectors at $(c,\omega)$ of the form $(w, 
0)$ with $w\in \kernel(\dd_{c})$. 

With the preceding understood, suppose that \eqref{eq:6.18} 
holds for some unit 
length complex number $\alpha $. Then, as $\dd_{c}^{\dagger} y_{R} = - 
x^{0}_{R}$ except where $t \ge R$ on the given component of 
$[0, \infty)  \times 
\partial X_{0}$, the sequence $\{y_{R}\}_{R\gg1}$ 
converges as  $R \to \infty$ to a 
non-zero section, $y$, of the bundle 
$i\cdot (\mathbb{R} \oplus \Lambda_{+}) \oplus S_{-}$ which obeys
\begin{equation}\begin{split}\label{eq:6.19}
&\dd_{c}^{\dagger} y = 0;\\
&\Rea(\bar{\alpha} w^{+}_{-}) = - \langle y, \kappa \rangle.
\end{split}\end{equation}
(As $\kappa $ has compact support, the integral in \eqref{eq:6.19} 
is well defined.)

\textbf {Step 7}\qua According to the preceding discussion, if 
$(c,\omega)$ is a critical point of $\varphi $, then there exists such 
unit length $\alpha\in\C$ such that 
\begin{equation}\label{eq:6.20}
0 = \langle y, \kappa \rangle 
\end{equation}
for all sections $\kappa $ of $i\Lambda_{+}$ which have compact support 
on $U$.

To make use of \eqref{eq:6.20}, it is also important to realize that $y$ 
is also 
$L^{2}$--orthogonal to all compactly supported sections of the $i \mathbb{R}$
summand in 
$i\cdot (\mathbb{R} \oplus \Lambda_{+}) \oplus S_{-}$. 
To see that such is the case, let 
$q$ denote a compactly supported section of the summand in question. Given 
$q$, 
there exists a unique, $L^{2}$ function $u$ on $X$ which obeys 
$d^*du + 2|\psi 
|^{2} u = q$. With $u$ in hand, then $\dd_{c}\underline {u} = q$ where 
$\underline{u} \equiv (du, - 2^{-1} u \psi )$. Thus,
\begin{equation}\label{eq:6.21}
\langle y_{R}, q\rangle = \langle x_{R}, \underline{u}\rangle. 
\end{equation}
To see that the left hand side of \eqref{eq:6.21} 
vanishes in the limit as $R \to\infty$, 
note that $x_{R}$ lies in the last two summands of \eqref{eq:5.1}. 
Meanwhile, Lemma 
\ref{thm:5.4} 
implies that the projection of $\underline {u}$ into these summands is 
$\ooo(e^{-zs})$ where $z > (2r)^{1/2}$ is independent of $s$. Given this last 
bound, the vanishing as $R$ tends to infinity of the right hand side of 
\eqref{eq:6.21} 
follows directly from \eqref{eq:6.12}.

To summarize: If $(c, \omega )$ is a critical point of the map $\varphi $, 
then there exists a non-zero section, $y$, 
over $X$ of $i\cdot(\mathbb{R}\oplus\Lambda_{+})\oplus S_{-}$ 
with the following properties:
{\sl
\begin{equation}\begin{split}\label{eq:6.22}
\bullet\text{ }&y\text{ is }L^{2}\text{--orthogonal to sections with compact}\\
&\text{support on }U\text{ of the }
i\cdot(\mathbb{R}\oplus\Lambda_{+})\text{ summand}.\\ 
\bullet\text{ }&\dd_{c}^{\dagger} y = 0.
\end{split}\end{equation}}
Now, the first point here implies that $y$ restricts to $U$ as a section of 
$S_{-}$ only. With this understood, then the projection of the equation in 
the second point of \eqref{eq:6.22} onto the $i T^*X$ 
summand of $i T^*X \oplus S_{+}$ asserts 
that 
\begin{equation*}
\Ima(\psi^{\dagger} \clif(e) y) = 0
\end{equation*}
for all sections $e$ of $T^*X$ 
with compact support on $U$. This last condition can 
hold only if $y$ is identically zero on $U$. 

Having established that $y$ vanishes identically on $U$, 
it then follows that $y 
\equiv 0$ on the whole of $X$ since there is a version of Aronszajn's 
unique continuation principle \cite{Ar} 
which holds for elements in the kernel of 
$\dd_{c}^{\dagger}$. 
The preceding conclusion establishes that $\varphi $ has no 
critical points as claimed. 
\end{proof}

\sh{b) Proof of Proposition \ref{thm:3.7}}

Given the details of the preceding proof of Proposition \ref{thm:2.4}, 
the assertions 
here follow via standard applications of the Sard--Smale theorem. The details 
for the application to this particular case are left to the reader.

\sh{c) Proof of Proposition \ref{thm:3.9}}

The argument given below considers only the case where $M$ separates $X$. 
As the 
discussion in the case where $X - M$ is connected is identical at all 
essential points to that given below, the latter discussion is omitted.

To begin, suppose, for the sake of argument, that there exists an increasing 
and unbounded sequence, $\{R_{j}\}_{j=1,2,\ldots}$ 
and a corresponding sequence 
$\{c_{j}\} \subset \mm^{R_{j}}$ with the property that for each fixed $r'$, 
the inequality in \eqref{eq:3.16} 
holds when $(A,\psi) = c_{j}$ for only finitely 
many $j$. Arguing as in Step 1 of the proof of Proposition \ref{thm:2.4} 
finds a 
subsequence of $\{c_{j}\}$ (hence renumbered consecutively), and data 
$c_{\infty-}$, $c_{\infty+}$ 
and $\{o_{j}\}_{j=1,2,\ldots}$, where now $c_{\infty\pm}$ are 
orbits of solutions to \eqref{eq:2.1} 
on the respective $X_{\pm} $. In this regard, 
each $X_{\pm} $ may have other ends besides the end where $M$ lived, and each 
such end has an associated map $\varphi $ as described in the preceding 
subsection. In this regard, assume 
that for each such end, $\varphi(c)\ne0$ 
for $c = c_{\infty\pm}$.

To return to the data $\{c_{j}\}$, 
$c_{\infty\pm}$, $\{o_{j}\}$, remark that an 
evident analog of \eqref{eq:6.1}--\eqref{eq:6.3} 
exists here. In particular, the orbit          
$c_{\infty-}$ supplies the complex 
number $u^{+}_{-}$ as defined in \eqref{eq:6.1} for 
the end $[0, \infty)  \times M$ of $X_{-}$, while $c_{\infty+}$ supplies the 
analogous $u^{-}_{-}$ which comes from the `time reversed' version of 
\eqref{eq:6.1} on the end 
$(-\infty,0]\times M$ of $X_{+}$. There is also an analog of 
\eqref{eq:6.7} and \eqref{eq:6.8} here:
\begin{equation}\begin{split}\label{eq:6.24}
&h_{j-} (u^{+}_{-} + \varepsilon_{j-}) \exp(-(2r)^{1/2} 
(R_{j} - t_{j-}))\\
+&h_{j+} (u^{-}_{-} + \varepsilon_{j+}) 
\exp(-(2r)^{1/2}(R_{j} - t_{j+})) = 0,
\end{split}\end{equation}
where
\begin{equation*}\begin{split}
&|\varepsilon_{j\pm}| \to 0\text{ as }j \to \infty;\\
&t_{j\pm}\to\infty\text{ and }|t_{j\pm}|/R_{j}\to0\text{ as }j\to\infty.
\end{split}\end{equation*}
In \eqref{eq:6.24}, 
the data $h_{j\pm}$ and $t_{j\pm}$ are supplied by the vortex 
$o_{j}$. In particular, $h_{j\pm} $ are real numbers, both greater than a 
$j$--independent, positive $\zeta^{-1}$, while $-t_{j-}$ and $t_{j+}$ are, 
respectively, the most negative and most positive of the $t$--coordinates of 
the zeros of the $o_{j}$ version $\tau$ in \eqref{eq:4.15}. With the preceding 
understood, it now proves useful to package \eqref{eq:6.24} as 
\begin{equation*}
(u^{+}_{-} + \varepsilon_{j-}) \Delta_{j} + (u^{-}_{-} + 
\varepsilon_{j+})(1 - \Delta_{j}) = 0,
\end{equation*}
where $\Delta_{j} \in [0, 1]$ is supplied by the vortex $o_{j}$. 

As $|\varepsilon_{j\pm}| \to 0$ as $j \to \infty$, the sequence 
$\{\Delta_{j}\}$ has a unique limit, $\Delta\in[0, 1]$ with
\begin{equation}\label{eq:6.27}
u^{+}_{-} \Delta + u^{-}_{-} (1 - \Delta ) = 0.
\end{equation}

Now, the set of non-zero pairs $(u^{+}_{-}, u^{-}_{-})$ which obey a 
relation as in \eqref{eq:6.27} for some $\Delta $ determines a codimension $1$ 
subvariety in $\C \times \C$. In particular, coupled with Lemma \ref{thm:6.1} 
and the 
Sard--Smale theorem, this last observation implies the next lemma. 

\begin{lem}\label{thm:6.2} 
{Given a form} $\omega'$ {as described in 
Proposition \ref{thm:3.9} and open sets} $U_{\pm}$ {with respective compact 
closures in the $\pm$ components of} $X_{0} - M$, {there exists a 
Baire subset of choices for} $\omega $ {in} \eqref{eq:2.1} {which agree 
with} $\omega'$ {on} $X - (U_{-}\cup U_{+})$, {and 
which have the following additional property:} {Let} $({s}, z) 
\in \sss_{0}(X_{0},\partial X_{0})$, let $(({s}_{-}, 
z_{-}),({s}_{+},z_{+})) \in \wp^{-1}(({s}, 
z))$ {and use} $\mm_{-}$ and $\mm_{+}$ to denote the 
corresponding moduli space of solutions to \eqref{eq:2.1} on $X_{-}$ 
{and}  $X_{+}$, {respectively.} {Then, there is a 
codimension one subvariety in the product,} $\mm_{-} \times \mm_{+}$, 
{which contains the only pairs} $(c_{-},c_{+})\in\mm_{-}  
\times \mm_{+}$ for which the corresponding pair of complex numbers 
$(u^{+}_{-}, u^{-}_{-})$ satisfies \eqref{eq:6.27} {for some 
choice of} $\Delta\in[0, 1]$. 
\end{lem}

Lemma \ref{thm:6.2} 
directly implies the assertions of Proposition \ref{thm:3.9} in the $d = 0$ 
case. Indeed, it is a consequence of Lemma \ref{thm:6.2} 
that the relevant subvariety 
in $\mm_{-} \times \mm_{+}$ must be empty as each $\mm_{\pm}$ has 
dimension zero. Thus, no solution to \eqref{eq:6.27} 
will exist and so for $R$ large, 
each point in $\mm^{R}$ is in the image of the map $\Phi $ from Proposition 
\ref{thm:3.8}.

To argue the $d > 0$ assertion of the proposition, choose all of the points in 
$\Lambda $ to lie on the $X_{-}$ side of $X$. According to Lemma 
\ref{thm:6.2}, the 
conclusion that $(\mm^{R})^{\underline\Lambda}$ 
lies in $\Phi$'s image for large $R$ 
fails only if the following is true: There exist $\mm_{\pm}$ with $\mm_{-}$ 
having dimension $2d$, $\mm_{+}$ having dimension $0$, and 
$\mm_{-}^{\underline\Lambda}   
\times \mm_{+} \subset \mm_{-} \times \mm_{+}$ intersecting the 
subvariety from Lemma \ref{thm:6.2}. However, as Lemma \ref{thm:6.2}'s 
subvariety has 
codimension $1$ and $\mm_{-}^{\underline\Lambda}\times\mm_{+}$ 
is a finite set, 
such an intersection is precluded by a suitably generic choice for the 
points in $\Lambda $.

\section{\boldmath $3$--dimensional implications}\label{sec:7}

The purpose of this final section is to discuss the implications of the 
theorems and propositions of the preceding sections in the special case 
where $X_{0} = S^{1}  \times  Y_{0}$ and $Y_{0}$ is either a compact 
$3$--manifold with positive first Betti number or else a compact $3$--manifold 
with boundary whose boundary components are all tori. 

Here are the key points to note with regard to such $X_{0}$: An argument 
from \cite{KM} readily adapts to prove that the first Chern class, $c({s})$, 
for any pair $({s}, z) \in \sss_{0}(X_{0},\partial X_{0})$ with 
non-zero Seiberg--Witten invariant is pulled up from $Y_{0}$. Moreover, as is 
explained below, the Seiberg--Witten invariants of $X_{0}$ are identical to 
invariants that are defined for $Y_{0}$ by counting solutions on $Y  
\equiv  Y_{0}  \cup ([0, \infty)\times\partial Y_{0})$ of a version of 
\eqref{eq:3.1}. 

To describe this version of \eqref{eq:3.1}, a Riemannian metric, 
a $\Spin^{\C}$ 
structure and a closed $2$--form $\omega_{0}$ must first be chosen. Having 
made these choices, the two equations in \eqref{eq:3.1} 
make perfectly good sense on 
any oriented $3$--manifold. However, for $Y$ as described in the preceding 
paragraph, the metric should be a product, flat 
metric on $[0,\infty)\times
\partial Y_{0}$, 
and the form $\omega_{0}$ must be non-zero and constant on each 
component of $[0, \infty)  \times\partial Y_{0}$. 
With this understood, the equations 
in \eqref{eq:3.1} should be augmented with the extra `boundary condition'
\begin{equation}\label{eq:7.1}
\int_{Y} |F_{A}|^{2} < \infty. 
\end{equation}

The invariants for $Y_{0}$ that the preceding paragraph mentions are then 
obtained via a count with appropriate algebraic weights of the orbits under 
the action of $C^{\infty}(Y;S^{1})$ of the solutions to \eqref{eq:3.1} 
which satisfy 
\eqref{eq:7.1}.

These $3$--manifold equations on $Y$ can be viewed as versions of 
\eqref{eq:2.1} on $X = 
S^{1}  \times Y$ which is how the equivalence between the Seiberg--Witten 
invariants for $Y_{0}$ and $X_{0}$ arise. For this purpose, consider the 
version of \eqref{eq:2.1} 
on $X = Y \times S^{1}$ when $X$ has its product metric 
and when $\omega = d\tau\wedge\theta + \omega_{0}$. Here, 
$\theta $ is the metric dual on $Y$ to $\omega_{0}$ and $d\tau$ is an 
oriented, unit length, constant $1$--form on $S^{1}$. In this case, note that 
solutions to \eqref{eq:3.1} 
on $Y$ which obey \eqref{eq:7.1} provide solutions to \eqref{eq:2.1} on $X$. 
Moreover, in the case where $c({s})$ is pulled up from $Y_{0}$, an 
integration by parts argument, much like that used to prove Proposition 
\ref{thm:4.1}, 
proves that all solutions to \eqref{eq:2.1} on 
$X$ are pull-backs of solutions to \eqref{eq:3.1} 
and \eqref{eq:7.1} on $Y$. 
Thus, Theorems \ref{thm:1.1} and \ref{thm:2.7} directly imply Mayer--Vietoris 
like theorems for the $3$--dimensional Seiberg--Witten equations. (In 
particular, Theorem \ref{thm:1.1} implies Theorem 5.2 in \cite{MT}.)   

The question arises as to whether the implications of Theorems 
\ref{thm:1.1} and \ref{thm:2.7} 
for the $3$--dimensional Seiberg--Witten equations 
can be proved with a strictly 
$3$--dimensional version of the arguments of the previous sections. The short 
answer is no as there are additional complications that arise and make for a 
somewhat more involved story. And, as the story here is already long enough, 
the additional discussion will not be provided, save for the brief comments 
of the next paragraph. 

Because the solutions to \eqref{eq:3.1} on $Y$ 
provide solutions to \eqref{eq:2.1} on $X$, the 
analysis in Sections \ref{sec:3}--\ref{sec:5} 
can be directly employed to characterize the 
behavior of the solutions to \eqref{eq:3.1}. 
However, there is one caveat: Properties 
that hold for solutions to \eqref{eq:2.1} 
on $X$ when the form $\omega $ is chosen from 
a Baire set may not hold for the solutions on $Y$ because a Baire set may be 
devoid of $2$--forms given by $d\tau\wedge\theta + \omega_{0}$ 
where $\omega_{0}$ is a closed form on $Y$. In particular, the following 
analog of Proposition \ref{thm:2.4} holds in the purely 
$3$--dimensional context:

\begin{prop}\label{thm:7.1} 
{Let} $\theta_{0}$ denote a 
closed $1$--form on $Y$ which is non-zero and constant on each component 
of $[0,\infty)\times\partial Y_{0}$. {With} $\theta_{0}$ {given, 
use a form} $\omega_{0}$ {in} \eqref{eq:3.1} whose metric dual 
agrees with $\theta_{0}$ {on} $[0, \infty)  \times\partial Y_{0}$. 
{Then each} $\mm_{s,m}\subset\mm_{s}$ is compact and 
contains {only a finite number of strata.} {Moreover, fix a 
closed $2$--form} $\omega'$ whose metric dual agrees with $\theta 
_{0}$ on $[0, \infty)  \times\partial Y_{0}$; and fix a non-empty, 
open set $U \subset Y_{0}$. {Then, there is a Baire set of 
smooth, closed $2$--forms} $\omega$ that agree with $\omega'$ 
{on} $X - U$ and have the following properties:
\begin{itemize}
\item  {Each stratum of} $\mm_{s}$ is a smooth manifold of 
dimension $0$. Moreover, the cokernel of the operator $\dd_{c}$ 
{vanishes for each} $c \in  \mm_{s}$.
\item  {The boundary of the closure in} $\mm_{s}$ {of any 
stratum intersects the remaining strata as a codimension $1$ submanifold.}
\end{itemize}
\end{prop}
Note that this last proposition implies the third assertion in Equation (4) 
of \cite{MT}.

The appearance in Proposition \ref{thm:7.1} 
of codimension $1$ submanifolds as opposed 
to codimension $2$ causes the added complications in the proof of the purely 
$3$--dimensional version of Theorem \ref{thm:2.7}. In particular, the purely 
$3$--dimensional versions of Propositions 
\ref{thm:3.7} and \ref{thm:3.9} may not hold. Even so, 
somewhat more complicated analogs of these propositions can be established 
that are sufficient to provide the purely $3$--dimensional proof of Theorem 
\ref{thm:2.7}.

\end{document}

%% file: gtoutput.tex
\def\ifplaintex{\expandafter\ifx\csname documentclass\endcsname\relax}

\ifplaintex 
\else
\fi

\expandafter\ifx\csname beginpicture\endcsname\relax
\expandafter\ifx\csname documentclass\endcsname\relax
\input pictex \else
\input prepictex \input pictex \input postpictex \fi\fi

\def\gt{{\mathsurround=0pt\it $\cal G\mskip-2mu$eometry \&\ 
$\cal T\!\!$opology}}        

\def\gtp{{\mathsurround=0pt\it $\cal G\mskip-2mu$eometry \&\ 
$\cal T\!\!$opology $\cal P\!$ublications}}  

\def\lognumber#1{\def\thelognumber{#1}}
\def\volumenumber#1{\def\thevolumenumber{#1}}
\def\papernumber#1{\def\thepapernumber{#1}}
\def\volumeyear#1{\def\thevolumeyear{#1}}

\def\pagenumbers#1#2{\def\startpage{#1}\def\finishpage{#2}}
\def\published#1{\def\publishdate{#1}}
\def\proposed#1{\def\theproposer{#1}}
\def\seconded#1{\def\theseconders{#1}}
\def\received#1{\def\receiveddate{#1}}

\def\accepted#1{\def\accepteddate{#1}}
\def\asciititle#1{\def\theasciititle{#1}}

\long\def\asciiabstract#1{\long\def\theasciiabstract{#1}}
\def\asciikeywords#1{\def\theasciikeywords{#1}}

\let\\\par\let\thelognumber\relax
\let\thevolumenumber\relax\let\thepapernumber\relax
\let\thevolumeyear\relax\let\thesamplenumber\relax\let\startpage\relax
\let\finishpage\relax\let\publishdate\relax\let\receiveddate\relax
\let\reviseddate\relax\let\accepteddate\relax\let\theasciititle\relax
\let\theasciiauthors\relax
\let\theasciiabstract\relax\let\theasciikeywords\relax
\let\theasciiemail\relax\let\theshortauthors\relax\let\theshorttitle\relax

\long\def\maketitlep{   

\count0=\startpage

\gt\hfill      
\beginpicture
\setcoordinatesystem units <0.33truein, 0.33truein> point at 2.2 0.9
\setplotsymbol ({$\cal G$})
\plotsymbolspacing=9truept
\circulararc 315 degrees from 0 1 center at 0 0
\setplotsymbol ({$\cal T$})
\circulararc 315 degrees from 1 -1 center at 1 0
\endpicture
\break
{\small\ifx\thesamplenumber\relax 
Volume \else Sample
\fi\thevolumenumber\ (\thevolumeyear)
\startpage--\finishpage\nl
Published: \publishdate}
\vglue 0.5truein plus 0.4fil minus 0.1truein

{\parskip=0pt\leftskip 0pt plus 1fil\def\\{\par\smallskip}{\ifplaintex\large
\else\Large\fi\bf\thetitle}\par\medskip}   

\vglue 0pt plus 0.1fil 

{\parskip=0pt\leftskip 0pt plus 1fil\def\\{\par}{\sc\theauthors}
\par\medskip}

\vglue 0pt plus 0.1fil 

{\small\parskip=0pt\let\newline\\
{\leftskip 0pt plus 1fil\def\\{\par}{\sl\theaddress}\par}
\expandafter\ifx\theemail\relax    
\relax\else\vglue 5pt plus 0.02fil minus 2pt\def\\{\stdspace{\rm 
and}\stdspace} 
\cl{Email:\stdspace\tt\theemail}\fi
\ifx\theurl\relax                  
\relax\else\vglue 5pt plus 0.02fil minus 2pt\def\\{\stdspace{\rm 
and}\stdspace}
\cl{URL:\stdspace\tt\theurl}\fi\par}

\vglue 7pt plus 0.3fil minus 3pt

{\bf Abstract}
\vglue 5pt plus 0.1fil minus 2pt

\theabstract

\vglue 7pt plus 0.3fil minus 3pt

{\bf AMS Classification numbers}\quad Primary:\quad \theprimaryclass

Secondary:\quad \thesecondaryclass

\vglue 5pt plus 0.3fil minus 2pt

{\bf Keywords}\quad \thekeywords

\vglue 10pt plus 0.5fil minus 5pt

{\small  Proposed: \theproposer\hfill Received: \receiveddate\nl
Seconded: \theseconders\hfill 
\ifx\reviseddate\relax                         
Accepted: \accepteddate                        
\else
Revised: \reviseddate                          
\fi}
\eject
}       

\let\maketitlepage\maketitlep
\let\maketitle\maketitlepage

\font\phead=cmsl9 scaled 950
\font=cmsl9 scaled 1050
\font\pnum=cmbx10 scaled 913
\font\lnum=cmbx10 
\font\pfoot=cmsl9 scaled 950
\font=cmsl9 scaled 1050
\ifplaintex
\headline{\vbox to 0pt{\vskip -4.5mm\line{\small\phead\ifnum
\count0=\startpage ISSN 1364-0380 (on line)
1465-3060 (printed) \hfill {\pnum\folio}\else\ifodd\count0\def\\{ }%
\ifx\theshorttitle\relax\thetitle\else\theshorttitle\fi\hfill{\pnum\folio}
\else\def\\{ and }{\pnum\folio}\hfill\ifx\theshortauthors\relax\theauthors
\else\theshortauthors\fi\fi\fi}\vss}}
\footline{\vbox to 0pt{\vglue 0mm\line{\small\pfoot\ifnum\count0=\startpage
\copyright\ \gtp\hfill\else
\gt, Volume \thevolumenumber\ (\thevolumeyear)\hfill\fi}\vss
}}
\else
\makeatletter
\def\@oddhead{{\small\lhead\ifnum\count0=\startpage ISSN 1364-0380 (on line)
1465-3060 (printed) \hfill {\lnum\number\count0}\else\ifodd\count0
\def\\{ }\ifx\theshorttitle\relax \thetitle \else\theshorttitle\fi\hfill
{\lnum\number\count0}\else\def\\{ and }{\lnum\number\count0}
\hfill\ifx\theshortauthors\relax 
\theauthors\else\theshortauthors\fi\fi\fi}}\def\@evenhead{\@oddhead}
\def\@oddfoot{\small\lfoot\ifnum\count0=\startpage\copyright\ \gtp\hfill\else
\gt, Volume \thevolumenumber\ (\thevolumeyear)\hfill\fi}
\def\@evenfoot{\@oddfoot}
\makeatother
\fi

\newwrite\gtoutfile
\long\gdef\makeheadfile{  
{\def\\{, }\def\s{ }
\immediate\openout\gtoutfile head.xxx
\immediate\write\gtoutfile{To: math@arxiv.org}
\immediate\write\gtoutfile{Subject: put or rep NNNNN:pppp}
\immediate\write\gtoutfile{--text follows this line--}
\immediate\write\gtoutfile{Proxy-for: \ifx\theasciiauthors\relax
\theauthors\else\theasciiauthors\fi\s<\ifx\theasciiemail\relax\theemail\else\theasciiemail\fi>}
\immediate\write\gtoutfile{\noexpand\\}
\immediate\write\gtoutfile{Authors: \ifx\theasciiauthors\relax
\theauthors\else\theasciiauthors\fi}
\immediate\write\gtoutfile{Title: \ifx\theasciititle\relax
\thetitle\else\theasciititle\fi}
\immediate\write\gtoutfile{Subj-class: GT or SG or MG etc}
\immediate\write\gtoutfile{MSC-class: \theprimaryclass\ifx\thesecondaryclass\relax\else, \thesecondaryclass\fi}
\immediate\write\gtoutfile{Journal-ref: Geom. Topol. \thevolumenumber
(\thevolumeyear) \startpage-\finishpage}
\immediate\write\gtoutfile{Comments: Published by Geometry and Topology at}
\immediate\write\gtoutfile{\s\s http://www.maths.warwick.ac.uk/gt/GTVol\thevolumenumber/paper\thepapernumber.abs.html}
\immediate\write\gtoutfile{\noexpand\\}
\immediate\write\gtoutfile{}
\ifx\theasciiabstract\relax
\immediate\write\gtoutfile{\theabstract}\else
\immediate\write\gtoutfile{\theasciiabstract}\fi
\immediate\write\gtoutfile{}
\immediate\write\gtoutfile{\noexpand\\}
\immediate\write\gtoutfile{}
\immediate\closeout\gtoutfile}}  

\def\maketitlepage{\maketitlep\makeheadfile}
\let\maketitle\maketitlepage

\def\ifplaintex{\expandafter\ifx\csname documentclass\endcsname\relax}

\ifplaintex 
\else
\fi

\expandafter\ifx\csname beginpicture\endcsname\relax
\expandafter\ifx\csname documentclass\endcsname\relax
\input pictex \else
\input prepictex \input pictex \input postpictex \fi\fi

\def\gt{{\mathsurround=0pt\it $\cal G\mskip-2mu$eometry \&\ 
$\cal T\!\!$opology}}        

\def\gtp{{\mathsurround=0pt\it $\cal G\mskip-2mu$eometry \&\ 
$\cal T\!\!$opology $\cal P\!$ublications}}  

\def\lognumber#1{\def\thelognumber{#1}}
\def\volumenumber#1{\def\thevolumenumber{#1}}
\def\papernumber#1{\def\thepapernumber{#1}}
\def\volumeyear#1{\def\thevolumeyear{#1}}

\def\pagenumbers#1#2{\def\startpage{#1}\def\finishpage{#2}}
\def\published#1{\def\publishdate{#1}}
\def\proposed#1{\def\theproposer{#1}}
\def\seconded#1{\def\theseconders{#1}}
\def\received#1{\def\receiveddate{#1}}

\def\accepted#1{\def\accepteddate{#1}}
\def\asciititle#1{\def\theasciititle{#1}}

\long\def\asciiabstract#1{\long\def\theasciiabstract{#1}}
\def\asciikeywords#1{\def\theasciikeywords{#1}}

\let\\\par\let\thelognumber\relax
\let\thevolumenumber\relax\let\thepapernumber\relax
\let\thevolumeyear\relax\let\thesamplenumber\relax\let\startpage\relax
\let\finishpage\relax\let\publishdate\relax\let\receiveddate\relax
\let\reviseddate\relax\let\accepteddate\relax\let\theasciititle\relax
\let\theasciiauthors\relax
\let\theasciiabstract\relax\let\theasciikeywords\relax
\let\theasciiemail\relax\let\theshortauthors\relax\let\theshorttitle\relax

\long\def\maketitlep{   

\count0=\startpage

\gt\hfill      
\beginpicture
\setcoordinatesystem units <0.33truein, 0.33truein> point at 2.2 0.9
\setplotsymbol ({$\cal G$})
\plotsymbolspacing=9truept
\circulararc 315 degrees from 0 1 center at 0 0
\setplotsymbol ({$\cal T$})
\circulararc 315 degrees from 1 -1 center at 1 0
\endpicture
\break
{\small\ifx\thesamplenumber\relax 
Volume \else Sample
\fi\thevolumenumber\ (\thevolumeyear)
\startpage--\finishpage\nl
Published: \publishdate}
\vglue 0.5truein plus 0.4fil minus 0.1truein

{\parskip=0pt\leftskip 0pt plus 1fil\def\\{\par\smallskip}{\ifplaintex\large
\else\Large\fi\bf\thetitle}\par\medskip}   

\vglue 0pt plus 0.1fil 

{\parskip=0pt\leftskip 0pt plus 1fil\def\\{\par}{\sc\theauthors}
\par\medskip}

\vglue 0pt plus 0.1fil 

{\small\parskip=0pt\let\newline\\
{\leftskip 0pt plus 1fil\def\\{\par}{\sl\theaddress}\par}
\expandafter\ifx\theemail\relax    
\relax\else\vglue 5pt plus 0.02fil minus 2pt\def\\{\stdspace{\rm 
and}\stdspace} 
\cl{Email:\stdspace\tt\theemail}\fi
\ifx\theurl\relax                  
\relax\else\vglue 5pt plus 0.02fil minus 2pt\def\\{\stdspace{\rm 
and}\stdspace}
\cl{URL:\stdspace\tt\theurl}\fi\par}

\vglue 7pt plus 0.3fil minus 3pt

{\bf Abstract}
\vglue 5pt plus 0.1fil minus 2pt

\theabstract

\vglue 7pt plus 0.3fil minus 3pt

{\bf AMS Classification numbers}\quad Primary:\quad \theprimaryclass

Secondary:\quad \thesecondaryclass

\vglue 5pt plus 0.3fil minus 2pt

{\bf Keywords}\quad \thekeywords

\vglue 10pt plus 0.5fil minus 5pt

{\small  Proposed: \theproposer\hfill Received: \receiveddate\nl
Seconded: \theseconders\hfill 
\ifx\reviseddate\relax                         
Accepted: \accepteddate                        
\else
Revised: \reviseddate                          
\fi}
\eject
}       

\let\maketitlepage\maketitlep
\let\maketitle\maketitlepage

\font\phead=cmsl9 scaled 950
\font=cmsl9 scaled 1050
\font\pnum=cmbx10 scaled 913
\font\lnum=cmbx10 
\font\pfoot=cmsl9 scaled 950
\font=cmsl9 scaled 1050
\ifplaintex
\headline{\vbox to 0pt{\vskip -4.5mm\line{\small\phead\ifnum
\count0=\startpage ISSN 1364-0380 (on line)
1465-3060 (printed) \hfill {\pnum\folio}\else\ifodd\count0\def\\{ }%
\ifx\theshorttitle\relax\thetitle\else\theshorttitle\fi\hfill{\pnum\folio}
\else\def\\{ and }{\pnum\folio}\hfill\ifx\theshortauthors\relax\theauthors
\else\theshortauthors\fi\fi\fi}\vss}}
\footline{\vbox to 0pt{\vglue 0mm\line{\small\pfoot\ifnum\count0=\startpage
\copyright\ \gtp\hfill\else
\gt, Volume \thevolumenumber\ (\thevolumeyear)\hfill\fi}\vss
}}
\else
\makeatletter
\def\@oddhead{{\small\lhead\ifnum\count0=\startpage ISSN 1364-0380 (on line)
1465-3060 (printed) \hfill {\lnum\number\count0}\else\ifodd\count0
\def\\{ }\ifx\theshorttitle\relax \thetitle \else\theshorttitle\fi\hfill
{\lnum\number\count0}\else\def\\{ and }{\lnum\number\count0}
\hfill\ifx\theshortauthors\relax 
\theauthors\else\theshortauthors\fi\fi\fi}}\def\@evenhead{\@oddhead}
\def\@oddfoot{\small\lfoot\ifnum\count0=\startpage\copyright\ \gtp\hfill\else
\gt, Volume \thevolumenumber\ (\thevolumeyear)\hfill\fi}
\def\@evenfoot{\@oddfoot}
\makeatother
\fi

\newwrite\gtoutfile
\long\gdef\makeheadfile{  
{\def\\{, }\def\s{ }
\immediate\openout\gtoutfile head.xxx
\immediate\write\gtoutfile{To: math@arxiv.org}
\immediate\write\gtoutfile{Subject: put or rep NNNNN:pppp}
\immediate\write\gtoutfile{--text follows this line--}
\immediate\write\gtoutfile{Proxy-for: \ifx\theasciiauthors\relax
\theauthors\else\theasciiauthors\fi\s<\ifx\theasciiemail\relax\theemail\else\theasciiemail\fi>}
\immediate\write\gtoutfile{\noexpand\\}
\immediate\write\gtoutfile{Authors: \ifx\theasciiauthors\relax
\theauthors\else\theasciiauthors\fi}
\immediate\write\gtoutfile{Title: \ifx\theasciititle\relax
\thetitle\else\theasciititle\fi}
\immediate\write\gtoutfile{Subj-class: GT or SG or MG etc}
\immediate\write\gtoutfile{MSC-class: \theprimaryclass\ifx\thesecondaryclass\relax\else, \thesecondaryclass\fi}
\immediate\write\gtoutfile{Journal-ref: Geom. Topol. \thevolumenumber
(\thevolumeyear) \startpage-\finishpage}
\immediate\write\gtoutfile{Comments: Published by Geometry and Topology at}
\immediate\write\gtoutfile{\s\s http://www.maths.warwick.ac.uk/gt/GTVol\thevolumenumber/paper\thepapernumber.abs.html}
\immediate\write\gtoutfile{\noexpand\\}
\immediate\write\gtoutfile{}
\ifx\theasciiabstract\relax
\immediate\write\gtoutfile{\theabstract}\else
\immediate\write\gtoutfile{\theasciiabstract}\fi
\immediate\write\gtoutfile{}
\immediate\write\gtoutfile{\noexpand\\}
\immediate\write\gtoutfile{}
\immediate\closeout\gtoutfile}}  

\def\maketitlepage{\maketitlep\makeheadfile}
\let\maketitle\maketitlepage